\documentclass{agtart_a}
\pdfoutput=1

\title{Vortices and a TQFT for Lefschetz fibrations on 4--manifolds}
\author{Michael Usher}
\givenname{Michael}
\surname{Usher}
\address{Department of Mathematics\\Princeton Universtity\\\newline
Fine Hall\\Washington Road\\Princeton, NJ  08544\\USA}
\email{musher@math.princeton.edu}
\urladdr{}

\volumenumber{6}
\issuenumber{}
\publicationyear{2006}
\papernumber{60}
\startpage{1677}
\endpage{1743}

\doi{}
\MR{}
\Zbl{}


\received{10 July 2006}
\revised{}
\accepted{29 August 2006}
\published{21 October 2006}
\publishedonline{21 October 2006}
\proposed{}
\seconded{}
\corresponding{}
\editor{}
\version{}

\arxivreference{math.SG/0603128}

\AtBeginDocument{\let\bar\wbar\let\tilde\wtilde\let\hat\what
\def\notin{\not\in}}
\makeop{Mor}
\makeop{Diff}
\makeop{Hom}
\makeop{Nov}
\makeop{fiber}
\makeop{Symp}
\makeop{symp}
\makeop{Spin}
\makeop{torsion}
\makeop{sec}
\makeop{vol}
\makeop{sign}
\def\ECH{\hbox{\textit{ECH}}}

\makeatletter
\def\cnewtheorem#1[#2]#3{\newtheorem{#1}{#3}[section]
\expandafter\let\csname c@#1\endcsname\c@theorem}
\makeatother

\newtheorem{theorem}{Theorem}[section]
\cnewtheorem{prop}[theorem]{Proposition}
\cnewtheorem{lemma}[theorem]{Lemma}
\cnewtheorem{remark}[theorem]{Remark}
\cnewtheorem{question}[theorem]{Question}
\cnewtheorem{definition}[theorem]{Definition}
\cnewtheorem{cor}[theorem]{Corollary}
\def\delbar{\bar{\partial}}

\begin{document}

\begin{asciiabstract}
Adapting a construction of D Salamon involving the U(1) vortex
equations, we explore the properties of a Floer theory for
3-manifolds that fiber over S^1 which exhibits several parallels
with monopole Floer homology, and in all likelihood coincides with
it. The theory fits into a restricted analogue of a TQFT in which
the cobordisms are required to be equipped with Lefschetz
fibrations, and has connections to the dynamics of surface
symplectomorphisms.
\end{asciiabstract}

\begin{htmlabstract}
Adapting a construction of D Salamon involving the U(1) vortex
equations, we explore the properties of a Floer theory for
3&ndash;manifolds that fiber over S<sup>1</sup> which exhibits several
parallels with monopole Floer homology, and in all likelihood
coincides with it. The theory fits into a restricted analogue of a
TQFT in which the cobordisms are required to be equipped with
Lefschetz fibrations, and has connections to the dynamics of surface
symplectomorphisms.
\end{htmlabstract}

\begin{abstract}
Adapting a construction of D Salamon involving the $U(1)$ vortex
equations, we explore the properties of a Floer theory for
3--manifolds that fiber over $S^1$ which exhibits several parallels
with monopole Floer homology, and in all likelihood coincides with
it. The theory fits into a restricted analogue of a TQFT in which
the cobordisms are required to be equipped with Lefschetz
fibrations, and has connections to the dynamics of surface
symplectomorphisms.
\end{abstract}

\maketitle
\section{Background and summary of results}\label{sec.1}

For some time it has been known that two of the most important
invariants of smooth closed 4--manifolds, the Donaldson and
Seiberg--Witten invariants, can each be expressed in terms of
$(3+1)$--dimensional topological quantum field theories (see Donaldson
\cite{Don-Floer}, Marcolli and Wang \cite{MW}, Kronheimer and Mrowka
\cite{KM}).  In such a ``TQFT,'' to each oriented 3--manifold $Y$
(perhaps equipped with additional data, such as a
spin$^c$--structure), one associates canonically a group $V(Y)$
satisfying, among several other conditions, the property that a
cobordism $X$ from $Y_1$ to $Y_2$ functorially induces a homomorphism
$F_X\co V(Y_1)\to V(Y_2)$.  If $X$ is a smooth closed oriented
4--manifold, divided into two pieces as $X=X_1\cup_Y X_2$ with
$b^+(X_i)>0$, one views $X_1$ as a cobordism from the empty set
$\varnothing$ to $Y$ and $X_2$ as a cobordism from $\varnothing$ to
$-Y$ (ie, $Y$ with its orientation reversed).  One has a natural
identification $V(-Y)\cong V(Y)^*$, and the $4$--dimensional invariant
$I_X$ is obtained by a natural calculation in $V(Y)$ involving the
images of the maps $F_{X_1}$ and $F_{X_2}$; $I_X$ is independent of
the choice of splitting of $X$ into the two pieces $X_1$ and $X_2$.

In the presence of a symplectic structure $\omega$ on the spin$^c$
4--manifold $(X,s)$, the famous work of C Taubes collected in
\cite{Tbook} shows that the Seiberg--Witten invariant $SW_X(s)$
agrees with a ``Gromov invariant'' $Gr_{(X,\omega)}(\alpha_s)$ which
counts pseudoholomorphic submanifolds of $X$ representing a homology
class $\alpha_s$ corresponding to $s$.  Kronheimer and Mrowka's work
\cite{KM} (see \cite{KMOS} for a summary) lays the full foundations
for the TQFT underlying $SW_{X}(s)$, in which the role of the group
$V(Y)$ in the above description is played by $HM(Y,s,\eta)$, where
$s$ is a spin$^c$ structure and $\eta\in H^2(Y;\mathbb{R})$ is the
cohomology class of the perturbation used in the Seiberg--Witten
equations. Given the correspondence between $SW$ and $Gr$, it is
natural to expect that there might be a TQFT underlying $Gr$ which
corresponds to Kronheimer--Mrowka's field theory for $HM$. Progress
in this direction has been made by M Hutchings and his
collaborators, who introduce groups $\ECH(Y)$ (embedded contact
homology) \cite{ECH} and $HP(Y)$ (periodic Floer homology)
\cite{PFH}, \cite{Dehn} in the case that $Y$ is, respectively, a
contact manifold or a mapping torus. These groups are conjectured to
agree with $HM$ or the (conjecturally equivalent) Heegaard Floer
homology $HF^+$ under suitable hypotheses, and do so in each of the
several cases that have been computed. However, at this writing a
number of foundational questions (such as independence of the choice of
almost complex structure) remain to be settled for $\ECH$ and $HP$,
and there do not presently exist full-blown TQFT's incorporating
either one of them.

According to results of S Donaldson \cite{Don-Lef} and R Gompf
\cite{Gompf}, a smooth oriented 4--manifold $X$ admits a symplectic
structure if and only if, possibly after blowing $X$ up at finitely
many points, there is a \emph{Lefschetz fibration} $f\co X\to S^2$
whose fibers are homologically essential.    Recall here that a
Lefschetz fibration on an oriented 4--manifold is a map to a
2--manifold which is a submersion except at its only finitely many
critical points, near each of which there are orientation preserving
complex coordinates in terms of which the map has the form
$(z,w)\mapsto zw$.  As such, the fibers of $f$ are all complex
curves of some fixed arithmetic genus, all but finitely many of
which are smooth, with the singular fibers having at worst nodal
singularities.  In the presence of a Lefschetz fibration $f\co X\to
S^2$ (satisfying certain properties that can always be achieved
using the constructions of \cite{Don-Lef}), Donaldson and I Smith
introduced in \cite{DS} an invariant $DS_{(X,f)}(\alpha)$ (for
$\alpha\in H_2(X;\mathbb{Z}))$ which counts pseudoholomorphic
sections of a bundle of symmetric products constructed from $f$.  In
\cite{U} it was shown that this Donaldson--Smith invariant coincides
with Taubes' invariant $Gr$, and hence also with the Seiberg--Witten
invariant under the appropriate identification of
$H_2(X;\mathbb{Z})$ with the set of spin$^c$ structures on $X$.

The present paper concerns what might be described as a restricted
TQFT which underlies the $4$--dimensional invariant $DS$.  We view
this TQFT as a covariant functor \emph{to} the category of modules
over a certain ring $A$ \emph{from} a category whose objects are
closed oriented 3--manifolds $Y$ equipped with fibrations $f\co Y\to
S^1$ (along with some additional structure indicated below) with
fiber genus at least 2,\footnote{Throughout this paper, the genera
of the fibers of all surface fibrations will be implicitly assumed
to be at least two, unless indicated otherwise.} with
$\Hom((Y_-,f_-),(Y_+,f_+))$ consisting of Lefschetz fibrations $f\co
X\to B$ over a base with two boundary components $\partial_-B$ and
$\partial_+B$ such that $f^{-1}(\partial_{\pm}B)=Y_{\pm}$ and
$f|_{Y_{\pm}}=f_{\pm}$.  The ``additional structure'' alluded to
earlier on an object $(Y,f)$ consists of a homology class \[h\in
H_1(Y;\mathbb{Z}),\] a cohomology class \[ c\in H^2(Y;\mathbb{R})\]
which evaluates positively on the fibers of $f$, and \[\mbox{a
 real number }\tau\in (2\pi h\cap[\fiber],+\infty).
\]  To each such tuple $(Y,f,h,c,\tau)$ and suitable ring $A$
(often, $A$ will be a Novikov ring) we associate an $A$--module
$HF(Y,f,h,c,\tau; A)$. These groups have appeared in the literature
before: in a paper of D Salamon \cite{Sal} they were conjectured to
agree with the (at the time not-yet-rigorously-defined) monopole
Floer groups of $Y$. As will be explained in more detail below, the
fibration $f\co Y\to S^1$ singles out a canonical spin$^c$
structure, which provides an identification of $\Spin^c(Y)$ with
$H_1(Y;\mathbb{Z})$.  Let $s_h$ be the spin$^c$ structure
corresponding to $h\in H_1(Y;\mathbb{Z})$ under this identification.
Salamon's conjecture can then be restated as saying that,
\[ HF(Y,f,h,c,\tau; A)\cong HM(Y,s_h,\eta(h,c,\tau); A),\]
where, in an appropriate normalization, $\eta(h,c,\tau)=4\tau c+2\pi
c_1(s_h)$. (The normalization on $\eta$ in this formula is such that
$c=0$ (if it were allowed) would correspond to a ``balanced
perturbation'' as in Kronheimer--Mrowka \cite{KM}, ie, a perturbation as in the
hypothesis of Conjecture 1.1 of Lee \cite{LHF}; our requirement that $c$
pair positively with the fiber thus ensures that reducible solutions
to the Seiberg--Witten equations will not enter the picture (so that
$\overline{HM}(Y,s_h,\eta(h,c,\tau))=0$ and there is just one
nontrivial monopole Floer group corresponding to the perturbation
$\eta(h,c,\tau)$, making the notation $HM(Y,s_h,\eta(h,c,\tau); A)$
unambiguous) and that, in case $b_1(Y)=1$, all allowed values of $c$
and $\tau$ will put $\eta(h,c,\tau)$ on the same side of the
``wall'' familiar from Seiberg--Witten theory.) Note here the
general principle that the choice of $h$ in $HF$ corresponds to the
choice of a spin$^c$ structure in $HM$, while (given $h$) the choice
of $c$ corresponds to the choice of a cohomology class of
perturbation 2--forms in $HM$ (up to a scale factor determined by
$\tau$).  We should caution that the fact that $\eta(h,c,\tau)$ is
not a balanced perturbation means that in general $HF$ is not
conjectured to coincide with the Heegaard Floer group $HF^+$;
rather, they should be related by a change of coefficients as
detailed in  \cite[Chapter VIII]{KM}.

The following four subsections summarize our results concerning
these groups $HF$.  \fullref{Define} contains the explicit construction of
the groups, after which the results of subsections \ref{sec.1.1},
\ref{sec.1.2}, \ref{sec.1.3} and \ref{sec.1.4} are proven in 
Sections \ref{sec.3}, \ref{sec.4}, \ref{sec.5} and \ref{sec.6}, respectively.

\subsection{Coefficient rings}\label{sec.1.1}  
Typically, the natural choice for
the coefficient ring $A$ in $HF(Y,f,h,c,\tau;A)$ will be a Novikov
ring $\tilde{\Lambda}_{h,c}$ or $\Lambda_{h,c}$ which (as the
notation indicates) depends on the choices of $h\in
H_1(Y;\mathbb{Z})$ and $c\in H^2(Y;\mathbb{R})$. Let $R$ be a ring
(usually $\mathbb{Z}$, $\mathbb{Z}/2$, or $\mathbb{Q}$), $G$ an
abelian group, and $N\co G\to \mathbb{R}$ a homomorphism. Following
the notation of Lee \cite{L1}, the Novikov ring $\Nov(G,N;R)$ is defined
to be the set of formal sums $\sum_{g\in G}a_g\cdot g$ $(a_g\in R)$
satisfying the property that for every $C\in \mathbb{R}$ we have
$\#\{g|a_g\neq 0 \mbox{ and } N(g)<C\}<\infty$, with addition and
multiplication in $\Nov(G,N;\mathbb{R})$ defined as the obvious
extensions of the corresponding operations on the group ring $R[G]$.
For us the most common Novikov rings will be, aside from the
universal Novikov ring mentioned below,
\begin{equation} \tilde{\Lambda}_{h,c}=\Nov\left(\ker \langle
c_1(s_h),\cdot\rangle,\langle c,\cdot\rangle;R\right)\end{equation}
and
\begin{equation} \label{novring} \Lambda_{h,c}=\Nov\left(\frac{\ker
\langle c_1(s_h),\cdot\rangle}{\ker \langle
c_1(s_h),\cdot\rangle\cap\ker \langle c,\cdot\rangle}, \langle
c,\cdot\rangle;R\right)\end{equation} where $\langle
c_1(s_h),\cdot\rangle$ and $\langle c,\cdot\rangle$ are the
evaluation homomorphisms $H_2(Y;\mathbb{Z})\to\mathbb{R}$.

Obviously, multiplying $c$ by a positive constant leaves
$\Lambda_{h,c}$ unchanged.  Recall also that $HM(Y,s,\eta)$
naturally has coefficients in \[ \Nov\left(\ker \langle
c_1(s),\cdot\rangle, \langle 2\pi(\pi
c_1(s)-\eta),\cdot\rangle;R\right),\] which is the same as
$\tilde{\Lambda}_{h,c}$ in the event that $s=s_h$ and
$\eta=\eta(h,c,\tau)$, consistently with Salamon's conjecture.

One checks easily that $\Lambda_{h,c}$ as defined above embeds via
the ring homomorphism \[ \sum a_g g\mapsto \sum a_g T^{\langle
c,g\rangle}\] as a subring of the universal Novikov ring \[
\Lambda_{\Nov}^{R}=\left\{\sum_{i} a_i T^{\lambda_i}|a_i\in R,\,
(\forall C>0)(\#\{i|\lambda_i<C\}<\infty)\right\}.\]

We obtain our groups $HF$ as the homology of a chain complex
$CF(Y,f,h,c,\tau)$ which naturally has its coefficients in
$\tilde{\Lambda}_{h,c}$. If $A$ is any algebra over
$\tilde{\Lambda}_{h,c}$, then $HF(Y,f,h,c,\tau;A)$ is by definition the
homology of $CF(Y,f,h,c,\tau)\otimes_{\tilde{\Lambda}_{h,c}} A$.  Of
course, as a particular example of this, the projection \[ \ker
\langle c_1(s_h),\cdot\rangle\to \frac{\ker \langle
c_1(s_h),\cdot\rangle}{\ker \langle c_1(s_h),\cdot\rangle\cap\ker
\langle c,\cdot\rangle}\] makes $\Lambda_{h,c}$ into an algebra over
$\tilde{\Lambda}_{h,c}$.

Observe that when $c=\pm c_1(s_h)$, $\Lambda_{h,c}$ is just the ring
$R$ over which we are working while
$\tilde{\Lambda}_{h,c}=R[\ker\langle c_1(s_h),\cdot\rangle]$, so no
Novikov ring is needed in this case. This choice of $c$ corresponds
to the choice made in the construction of periodic Floer homology;
see \cite{Dehn}.  In this case, for any other $\tilde{c}\in
H^2(X;\mathbb{R})$ evaluating positively on the fibers of $f$,
$\tilde{\Lambda}_{h,\tilde{c}}$ and $\Lambda_{h,\tilde{c}}$ are
obviously algebras over $R[\ker\langle
c_1(s_h),\cdot\rangle]=\tilde{\Lambda}_{h,\pm c_1(s_h)}$, so if
$\langle c_1(s_h),\fiber\rangle\neq 0$ we have a well defined group
$HF(Y,F,h,\pm c_1(s_h),\tau;\tilde{\Lambda}_{h,\tilde{c}})$, where
the sign at the front of $\pm c_1(s_h)$ is chosen to make its
evaluation on the fiber positive.

We mention here that, while the presence of a Novikov ring such as
$\tilde{\Lambda}_{h,c}$ as the natural coefficient ring is a
standard aspect of Floer theory, the fact that this Novikov ring is
described directly in terms of the homology of $Y$ and is (crucially
for the invariance theorem below) independent of $\tau$ is more
subtle.  This fact follows from two basic ingredients: a formula
from Perutz \cite{P} for a certain cohomology class on the vortex moduli
space; and an expression for the evaluation of one of the terms in
that formula on certain cycles, derived below as Equation
\ref{dsquared}, which enables us to choose the forms
$\omega_{c,h,\tau}$ at the start of Section 2 in such a way as to
arrange that $\tilde{\Lambda}_{h,c}$ be the appropriate Novikov
ring.

 The following theorem, whose proof uses the $\tau$-independence of
 the Novikov ring along with the
difficult bifurcation analysis carried out by Y-J Lee in
\cite{L1,L2}, shows that at least in the majority of cases,
the groups $HF$ are independent of $\tau$, and depend on $c$ only to
the extent that $c$ determines the appropriate coefficient ring.

\begin{theorem} Writing $d=\langle PD(h),\fiber\rangle$ and letting $g$ be the genus of the
fibers of $f\co Y\to S^1$, assume that either $d\geq g-1$ or
$d<(g+1)/2$.  Then: \label{coeffs} \begin{itemize}\item[(i)]
$CF(Y,f,h,c,\tau)$ and $CF(Y,f,h,c,\tau')$ are canonically chain
homotopy equivalent whenever $\tau,\tau'\in \mathbb{R}$ are such
that both chain complexes are defined.
\item[(ii)]  Assume that $\tilde{c}$ and $\pm c_1(s_h)$ both
evaluate positively on the fiber.  Then for any $\tau>2\pi d$,
$CF(Y,f,h,\tilde{c},\tau)$ is chain homotopy equivalent to\[
CF(Y,f,h,\pm c_1(s_h),\tau)\otimes_{R[\ker\langle
c_1(s_h),\cdot\rangle]} \tilde{\Lambda}_{h,\tilde{c}}.\]
\end{itemize}
\end{theorem}

A similar result holds for the dependence of $HM(Y,s,\eta)$ on
$\eta$, as is shown in Section 31 of the current draft version of
\cite{KM}.

\subsection{Grading, module structure, duality, and local coefficients}
\label{sec.1.2} 
The groups $HF(Y,f,h,c,\tau;A)$ share additional algebraic structure
with the monopole Floer groups, as has been observed independently in
the thesis of T Perutz \cite{P}.  Throughout the paper, to a
$3$--manifold $Y$ associate the graded ring \[
\mathbb{A}(Y)=\mathbb{Z}[U]\otimes \Lambda^*(H_1(Y)/\torsion),\] where
$U$ is a formal variable of degree 2 and elements of $H_1(Y)$ have
degree $1$.  Then:
\begin{prop} \label{alg}
 $HF(Y,f,h,c,\tau; A)$ is a naturally $\mathbb{Z}/2$-graded,
relatively $\mathbb{Z}/\frak{d}(s_h)$-graded module over
$\mathbb{A}(Y)$, where
\[ \frak{d}(s_h)=\gcd_{T\in H_2(Y;\mathbb{Z})}\langle
c_1(s_h),T\rangle\]  and the action of an element of degree $p$ of
$\mathbb{A}(Y)$ on $HF(Y,f,h,c,\tau;A)$ decreases the relative
grading by $p$.\end{prop}

Since $HM$ likewise enjoys these properties (as seen, for instance,
in  \cite[sections 3.1 and 3.2]{KM}), it is natural to embellish
Salamon's conjecture to state that the (conjectural) isomorphism
between $HF$ and $HM$ is an isomorphism of graded modules.

For the ``Poincar\'e duality'' property (which, like \fullref{alg}, also appears in \cite{P}), given $(Y,f,h,c,\tau)$, let
$(-Y,\bar{f},h,c,\tau)$ be obtained by reversing the orientation of
$Y$ and composing $f$ with complex conjugation.

\begin{prop} \label{pairing} There is a perfect pairing \[
\langle\cdot,\cdot\rangle\co CF(Y,f,h,c,\tau;A)\otimes
CF(-Y,\bar{f},-h,c,\tau;A)\to A\] which satisfies \[ \langle
\partial_Y a,b\rangle=\langle a,\partial_{-Y} b\rangle\] and hence
descends to a pairing which identifies $HF(-Y,\bar{f},-h,c,\tau;A)$
with the dual of $HF(Y,f,h,c,\tau;A)$.\end{prop}

A handy device in monopole Floer theory is the use of ``local
coefficients'' in $HM$, in which a singular $1$--cycle $\gamma$ in
$Y$ gives rise to a twisted version $HM(Y,s,\eta;\Gamma_{\gamma})$
arising from a twisted coefficient system $\Gamma_{\gamma}$ on
configuration space associated to $\gamma$; homologous $1$-cycles
yield isomorphic coefficient systems and hence isomorphic Floer
groups, but a homology between the cycles must be specified in order
to make these isomorphisms canonical.  A parallel situation, which
seems most naturally expressed in terms of closed $2$--forms on $Y$
rather than their dual $1$--cycles, exists for our Floer groups, and
we expect the resulting twisted groups to be isomorphic to their
monopole Floer counterparts:

\begin{prop}  To each closed $2$--form $\theta\in \Omega^2(Y)$, we
may associate a local coefficient system $\Gamma_{\theta}$ and hence
Floer groups $HF(Y,f,h,c,\tau;\Gamma_{\theta})$.  If
$\theta_1,\theta_2$ are closed 2--forms, to each
$\zeta\in\Omega^1(Y)$ such that $d\zeta=\theta_2-\theta_1$ there is
associated a canonical isomorphism $\phi_{\zeta}\co
\Gamma_{\theta_1}\to \Gamma_{\theta_2}$, which then induces an
isomorphism of the associated Floer groups.\end{prop}

\subsection{Cobordisms}\label{sec.1.3} 

We define a category $\mathsf{FCOB}$ (for ``fibered cobordism'') as
follows.  An object $\frak{o}$ of  $\mathsf{FCOB}$ is a quintuple
$\frak{o}=(Y,f,h,c,\tau)$ where, as before, $Y$ is an oriented
3--manifold (we allow $Y$ to be empty, in which case $f,h,c,\tau$
need not be specified), $f\co Y\to S^1$ is a fibration, $h\in
H_1(Y;\mathbb{Z})$, $c\in H^2(Y;\mathbb{R})$, and $\tau\in(2\pi
h\cap[\fiber],+\infty)$ (when $h\cap[\fiber]$ is outside the interval
$[(g+1)/2,g-1)$ \fullref{coeffs} ensures that the Floer homology
$HF(\frak{o})$ associated to the object will be independent of
$\tau$ up to canonical isomorphism; these isomorphisms will commute
with the homomorphisms decribed below). A morphism
$\frak{m}=(X,\tilde{f},\tau)$ in $\Mor(\frak{o_-},\frak{o_+})$
consists of a Lefschetz fibration $\tilde{f}\co X\to B$ defined on a
4--manifold $X$ with oriented boundary $\partial X=(-Y_-)\coprod
Y_+$, with two-dimensional image $B$ having boundary components
$\partial_-B$, $\partial_+B$; here $\partial_{\pm} B=S^1$ if
$Y_{\pm}$ is nonempty, and otherwise $\partial_{\pm}B=\varnothing$.
We require $\tilde{f}^{-1}(\partial_{\pm}B)=Y_{\pm}$ and
$\tilde{f}|_{Y_{\pm}}=f_{\pm}$, and $\tau=\tau_+=\tau_-$.
Furthermore, where we denote by $\partial_{\pm}\co H_2(X,\partial
X;\mathbb{Z})\to H_1(Y_{\pm};\mathbb{Z})$ the obvious maps induced
by restriction to the boundary, we require that the sets
\begin{align*} C_{c_-,c_+}=\{\tilde{c}\in
H^2(X;\mathbb{R})|\tilde{c}|_{Y_{\pm}}&=c_{\pm} \}\end{align*} and
\[ H_{h_-,h_+}=\{\tilde{h}\in H_2(X,\partial
X;\mathbb{Z})|\partial_{\pm}\tilde{h}=h_{\pm}\}\] both be nonempty.

We also place the following additional structure on $B$:
\begin{definition} A ``starred surface with boundary'' is an
oriented surface $B$ with (say) genus $g$ and $n$ boundary
components equipped with distinguished points, arcs, and
parametrized loops as follows:
\begin{itemize}
\item[(0)] The distinguished points comprise one ``interior
base point'' $b$, $s\geq 0$ `interior special points''
$p_1,\ldots,p_s$, and one ``boundary base point'' $q_j$ on each of
the $n$ boundary components, and
\item[(1)] There are $n+s$ distinguished arcs, namely
one from $b$ to $p_i$ for each $i$ and one from $b$ to $q_j$ for
each $j$, one distinguished loop (namely the boundary component)
based at each $q_j$, and $2g$ distinguished interior loops
$\alpha_1,\ldots,\alpha_g,$\break $\beta_1,\ldots,\beta_g$, such that
$\alpha_1,\ldots,\alpha_g$ are based at $b$, are linearly
independent in $H_1(B;\mathbb{Z})$, and represent homotopy classes
having zero geometric intersection number, and
$\beta_1,\ldots,\beta_g$ are disjoint, with $\#(\alpha_i\cap
\beta_j)=\delta_{ij}$.\end{itemize} We require also that the
distinguished arcs and loops have no intersections other than the
ones implied by the above conditions.\end{definition}

We equip the base $B$ of the Lefschetz fibration $\tilde{f}\co X\to
B$ with the structure of a starred surface with boundary, with
interior special points comprising precisely the critical values of
$\tilde{f}$, and with the boundary loop at $\partial_- B$ (resp.\
$\partial_+B$) negatively (resp.\ positively) oriented. Finally, for
technical reasons that shall appear in \fullref{relfib} and its
proof, in the event that $g(B)>0$ we take as given in the data of
the morphism $\frak{m}$ a set of cohomology classes $b_i\in
H^2(\tilde{f}^{-1}(\beta_i);\mathbb{R})$ with the property that
there are $\tilde{c}\in C_{c_-,c_+},\tilde{h}\in H_{h_-,h_+}$ such
that $b_i=(\tilde{c}+2\pi
PD(\tilde{h})/\tau)|_{\tilde{f}^{-1}(\beta_i)}$ for
$i=1,\ldots,g(B)$. For the sake of conciseness, we shall nonetheless
generally denote morphisms with the notation
$\frak{m}=(X,\tilde{f},\tau)$, suppressing our marking of $B$ and
our choice of classes in the
$H^2(\tilde{f}^{-1}(\beta_i);\mathbb{R}).$

 Note that any two
starred surfaces $B_1,B_2$ with interior base points $b_1,b_2$
nonempty boundary may be glued along any of their common boundary
components $S_1, S_2$ (if necessary after reversing the
parametrization of one of the $S_i$) to obtain a new starred surface
$B_1\,_{S_1}\sharp_{S_2} B_2$ with (possibly empty) boundary as
follows. We join the corresponding boundary components at their
corresponding boundary base points; this in particular yields a path
$\gamma$ from $b_1$ to $b_2$.  To get the new distinguished arcs,
delete the loop at the former boundary base point, and extend all
the paths (and loops) based at $b_2$ in $B_2$ to paths (and loops)
based at $b_1$ in $B_1\sharp B_2$ by adding on the path $\gamma$
from $b_1$ to $b_2$. This construction applies equally well when
$B_1=B_2=B$ (as long as $S_1\neq S_2$); in this case
$B_{S_1}\sharp_{S_2}B$ will have genus one larger than $B$, and the
path $\gamma$ appears as a new $\alpha$-curve while the loop
resulting from the fusing of the old boundary components becomes a
new $\beta$-curve. In particular, a starred Riemann surface $B$ with
genus $g$ and $n$ boundary components can be cut along its
$\beta$-curves to obtain a starred Riemann surface $B^0$ with genus
$0$ and $n+2g$ boundary components, so that $B$ is recovered from
$B^0$  by applying the gluing construction $g$ times.\footnote{One
might naturally ask whether it is really necessary to incorporate
the complication of starred surfaces into the theory; the reason we
have done so is that it will help us construct certain $2$--forms on
symmetric product bundles in an essentially canonical way. If we
were only considering a single cobordism, say from $\frak{o}_-$ to
$\frak{o}_+$, it would suffice for the construction of the desired
$2$--form to replace the starred surface structure with an auxiliary
choice of a cohomology class from the set $C_{c_-,c_+}$. However, we
wish to compose our cobordisms, and the set of possible choices of
cohomology class on the composed cobordism having the appropriate
restrictions to the pieces is generally a positive-dimensional
affine space, so that there is no canonical way to glue cohomology
classes.}

Given an object $\frak{o}=(Y,f,h,c,\tau)$ of $\mathsf{FCOB}$, we
have a well-defined group\break $HF(\frak{o};\Lambda_{\Nov}^{R})$, where
again $\Lambda_{\Nov}^{R}$ is the universal Novikov ring over the
ring $R$; more generally if $\theta\in\Omega^2(Y)$ is closed we may
consider the twisted Floer homology group
$HF(\frak{o};\Lambda_{\Nov}^{R}\otimes \Gamma_{\theta})$ . A crucial
property of $HF$, a similar version of which was independently
discovered in \cite{P}, is that it is a functor from $\mathsf{FCOB}$
to the category $\mathsf{MOD}_{\Lambda_{\Nov}^{R}}$ of modules over
$\Lambda_{\Nov}^{R}$.

Denote by $\varnothing$ the object of $\mathsf{FCOB}$ whose
underlying 3--manifold is the empty set.  By way of definition, we
set $HF(\varnothing;\Lambda_{\Nov}^{R})=\Lambda_{\Nov}^{R}$.

Given morphisms $\frak{m}_0\in \Mor(\frak{o}_0,\frak{o}_1)$ and
$\frak{m}_1\in \Mor(\frak{o}_1,\frak{o}_2)$, we define the composite
morphism $\frak{m}_1\circ \frak{m}_0\in \Mor(\frak{o}_0,\frak{o}_2)$
by the obvious procedure of gluing the total spaces $X_0$ and $X_1$
of the Lefschetz fibrations $\tilde{f}_0,\tilde{f}_1$ underlying
$\frak{m}_0$, $\frak{m}_1$ along their common boundary component
$Y_1$ to obtain a new Lefschetz fibration $f\co X\to B$; as noted
earlier $B$ inherits the structure of a starred surface with
boundary from those of the bases of the $\tilde{f}_i$. If we choose
$\tilde{c}_0\in C_{c_0,c_1}$ and $\tilde{c}_1\in C_{c_1,c_2}$, so
that in particular $\tilde{c}_0$ and $\tilde{c}_1$ have the same
restriction to $Y_1$, then the Mayer--Vietoris sequence reveals that
the set of $\tilde{c}\in H^2(X;\mathbb{R})$ such that
$\tilde{c}|_{X_i}=\tilde{c}_i$ ($i=0,1$) is nonempty and is an
affine space over the image of the boundary map $\delta\co
H^1(Y;\mathbb{R})\to H^2(X;\mathbb{R})$. Poincar\'e duality implies
a parallel statement for the $h_i$.  In particular
$\frak{m}_1\circ\frak{m_0}$ is a morphism (for $C_{c_0,c_2}$ and
$H_{h_0,h_2}$ are nonempty, and where relevant the $\beta$-curves on
$B$ are just those on the $B_i$, so we can use the same cohomology
classes on the preimages of the $\beta$-curves as were used on the
$\frak{m}_i$).

\begin{theorem} \label{cobmaps} To each morphism $\frak{m}=(X,\tilde{f},\tau)$ from $\frak{o}_-$ to $\frak{o}_+$, where
$\frak{o}_{\pm}=(Y_{\pm},f_{\pm},h_{\pm},c_{\pm},\tau)$, and to each
closed form $\theta\in\Omega^2(X)$ vanishing near the critical
points of $f\co X\to B$, we may associate a homomorphism
\[ F_{\frak{m},\theta}\co \mathbb{A}(X)\otimes
HF(\frak{o}_-;\Lambda_{\Nov}^{R}\otimes\Gamma_{\theta|_{Y_-}})\to
HF(\frak{o}_+;\Lambda_{\Nov}^{R}\otimes\Gamma_{\theta|_{Y_+}})\]
where $\mathbb{A}(X)=\mathbb{Z}[U]\otimes
\Lambda^{*}(H_1(X;\mathbb{Z})/\torsion)$; in fact, each map
$F_{\frak{m},\theta}$ decomposes naturally as a sum
\[ F_{\frak{m},\theta}=\sum_{\tilde{h}\in H_{h_-,h_+}}
F_{\frak{m},\theta,\tilde{h}},\] and these maps enjoy the following
properties:\begin{itemize}
\item[(i)] For morphisms $\frak{m}_0=(X_0,f_0,\tau)$ from
  $\frak{o}_0=(Y_0,f_0,h_0,c_0,\tau)$ to $\frak{o}_1=$\break $(Y_1,f_1,h_1,c_1,\tau)$ and $\frak{m}_1=(X_1,f_1,\tau)$ from $\frak{o}_1$ to $\frak{o}_2=(Y_2,f_2,h_2,c_2,\tau)$, for $\theta$ a closed $2$-form on the total space $X=X_0\cup_YX_1$ of the Lefschetz fibration underlying $\frak{m}_1\circ\frak{m}_0$, for $\tilde{h}_0\in H_{h_0,h_1},\tilde{h}_1\in H_{h_1,h_2}$, and for $v\in
HF(\frak{o}_0;\Lambda_{\Nov}^{R}\otimes\Gamma_{\theta|_{Y_0}})$ we
have
\begin{align*} \sum_{\begin{array}{c}\scriptstyle\tilde{h}\in H_{h_0,h_2}:\\
\scriptstyle\tilde{h}|_{X_0}=\tilde{h}_0,\tilde{h}|_{X_1}=\tilde{h}_1\end{array}}
F_{\frak{m}_1\circ\frak{m}_0,\theta,\tilde{h}}&(U^{k+l}\otimes
1\otimes v)=\\
&F_{\frak{m}_1,\theta|_{X_1},\tilde{h}_0}(U^{k}\otimes 1\otimes
F_{\frak{m}_0,\theta|_{X_0},\tilde{h}_1}(U^l\otimes 1\otimes
v)).\end{align*}
\item[(ii)]
Where $i_-\co \mathbb{A}(Y_-)\to \mathbb{A}(X)$ is the map induced
by the action of the inclusion $Y_-\subset X$ on $H_1$,
$F_{\frak{m},\theta,\tilde{h}}$ is compatible with the module
structure of $HF$ in the sense that, for $\lambda\in
\mathbb{A}(Y_-)$ and $v\in
HF(\frak{o}_-;\Lambda_{\Nov}^{R}\otimes\Gamma_{\theta})$,\[
F_{\frak{m},\theta,\tilde{h}}(1\otimes \lambda\cdot
v)=F_{\frak{m},\theta,\tilde{h}}(i_-(\lambda)\otimes v).\]
\item[(iii)]  Where $-\frak{m}\in \Mor(\frak{o}_+,\frak{o}_-)$ is the
morphism obtained by reversing the orientation of the boundary
components of the base of $\tilde{f}$, with respect to the pairing
in \fullref{pairing} we have \[ \langle
F_{\frak{m},\theta,\tilde{h}}(v),w\rangle_{\frak{o}_+}=\langle
v,F_{-\frak{m},\theta,\tilde{h}}(w)\rangle_{\frak{o}_-}.\]
\item[(iv)]  Suppose that $\tilde{f}\co X\to \Sigma$ is a Lefschetz fibration on
the \emph{closed} manifold $X$ over the closed surface $\Sigma$, so
that $\frak{o}_-=\frak{o}_+=\varnothing$, and
$H_{h_-,h_+}=H_2(X;\mathbb{Z})$. Then  for a certain homomorphism
$A\co H_2(X;\mathbb{Z})\to\mathbb{R}$ and for $\frak{m}=(X,f,\tau)$
and $\theta\in \Omega^2(X)$ representing $[\theta]\in
H^2(X;\mathbb{R})$ we have,
\begin{align}\label{invtform} F_{\frak{m},\theta,\tilde{h}}(U^{r}\otimes
(\eta_1\wedge\cdots&\wedge\eta_k)\otimes 1)=
T^{A(\tilde{h})}e^{\langle[\theta],\tilde{h}\rangle}DS_{(X,\tilde{f})}(\tilde{h};pt^r,\eta_1,\ldots,\eta_k).\end{align}
Here $DS_{(X,\tilde{f})}(\tilde{h};pt^r,\eta_1,\ldots,\eta_k)$ is
the obvious extension of the Donaldson--Smith invariant \cite{DS} to
an invariant counting sections of the relative Hilbert scheme of $f$
which correspond to surfaces in $X$ representing $\tilde{h}\in
H_2(X;\mathbb{Z})$ and passing through $r$ generic points and
through generic cycles representing
$\eta_1,\ldots,\eta_k$.\end{itemize}\end{theorem}

 Our TQFT thus contains the
Donaldson--Smith invariant, which, thanks to \cite{Tbook} and
\cite{U}, is known to agree with the four-dimensional
Seiberg--Witten invariant under a natural (given the symplectic or
Lefschetz fibration structure on $X$) correspondence between
$H_2(X;\mathbb{Z})$ and $\Spin^c(X)$.  As is well-known (and shown in
detail in Chapter VII of \cite{KM}), there is a TQFT in
Seiberg--Witten theory whose cobordism maps enjoy properties exactly
parallel to those of \fullref{cobmaps}, with the Seiberg--Witten
invariant appearing in place of the (equivalent, by \cite{Tbook} and
\cite{U}) Donaldson--Smith invariant in part (iv). As such, we may
further embellish Salamon's conjecture to state that the conjectural
isomorphisms between $HF$ and $HM$ commute with the cobordism maps;
the agreement of $DS$ with $SW$ would then be a shadow of this
relationship.

In \cite{P}, Perutz uses and extends constructions similar to this
in order to construct a ``Lagrangian matching invariant'' for the
singular Lefschetz fibrations constructed in Auroux--Donaldson--Katzarkov
\cite{ADK} (which exist
on blowups for any 4--manifold with $b^+>0$, though it is not known
whether Perutz's invariant is independent of the choice of singular
Lefschetz fibration on a given 4--manifold), and conjectures that
this new invariant, too, agrees with the Seiberg--Witten invariant.

\subsection{Relation to dynamics of surface symplectomorphisms}\label{sec.1.4} 

Given a symplectomorphism $\phi\co (\Sigma,\omega)\to
(\Sigma,\omega)$ of a symplectic 2--manifold of genus at least 2,
form the mapping torus $Y_{\phi}=\mathbb{R}\times\Sigma/(t+1,x)\sim
(t,\phi(x))$; this fibers over $S^1$ and carries a fiberwise
symplectic form $\omega_{\phi}$ obtained by pushing forward the
obvious form induced by $\omega$ on $\mathbb{R}\times \Sigma$ via
the projection $\mathbb{R}\times\Sigma\to Y_{\phi}$.  Let $j$ be an
almost complex structure on the fibers of $f\co Y_{\phi}\to S^1$ and
$h\in H_1(Y_{\phi};\mathbb{Z})$; where $e$ is the Euler class of the
vertical tangent bundle of $Y_{\phi}\to S^1$, under several
assumptions on $\phi$ and $j$, including that $[\omega_{\phi}]\in
H^2(Y_{\phi};\mathbb{R})$ is proportional to $c_1(s_h)=e+2PD(h)$,
the \emph{periodic Floer homology} $HP(\phi,h,j)$ is defined in
\cite{PFH} to be the homology of a chain complex $CP(\phi,h)$ whose
generators are ``admissible orbit sets'' $\alpha=\{(\alpha_i,m_i)\}$
such that $\sum_{i}m_i[\alpha_i]=h$. Here the $\alpha_i$ are
periodic orbits for $\phi$ (which then give rise naturally to loops
in $Y$), and the $m_i$ are positive integers such that $m_i=1$ if
$\alpha_i$ is hyperbolic.  The matrix element $\langle\partial
\alpha,\beta\rangle$ for the boundary operator $\partial$ of the
chain complex counts certain embedded holomorphic curves $C$ in
$\mathbb{R}\times Y$ such that $C\cap(\{t\}\times Y)$ is asymptotic
to $\alpha$ (resp.\ $\beta$) as $t\to -\infty$ (resp.\ $t\to
+\infty$).  Note that $CP(\phi,h)$ is independent of $j$ as a graded
group; the same is expected to be true of $HP(\phi,h,j)$, but this
has not yet been proven.

$HP(\phi,h,j)$ is defined over the coefficient ring $\mathbb{Z}$;
more generally, if $[\omega_{\phi}]=c\in H^2(X;\mathbb{R})$, one
could define a periodic Floer homology $HP(\phi,h,c,j)$ over the
same Novikov ring $\tilde{\Lambda}_{h,c}$ as in \eqref{novring}.

Let $d=h\cap[\fiber]$.  $\phi$ then induces a continuous (but usually
not differentiable) map $S^d\phi\co S^d\Sigma\to S^d\Sigma$, where
$S^d\Sigma$ is the $d$th symmetric product of $\Sigma$.  Letting
$Y_{S^d\phi}$ denote the mapping torus of $S^d\phi$, $h\in
H_1(Y_{\phi};\mathbb{Z})$ naturally determines a homotopy class
$p_h$ of sections of $Y_{S^d\phi}$, and any generator
$\alpha=\{(\alpha_i,m_i)\}$ of $CP(\phi,h)$ corresponds to a fixed
point $D_{\alpha}$ of $S^d\phi$ such that the ``constant section''
of $Y_{S^d\phi}$ at $D_{\alpha}$ represents $p_h$ (in such a
situation we say $D_{\alpha}$ is ``in the $p_h$-sector''). The
admissibility condition on the hyperbolic orbits in $\alpha$
prevents this from being a one-to-one correspondence. Notably, at
least in some simple cases, one can perturb the (usually only
H\"older continuous) map $S^d\phi$ to a smooth map $\Phi\co
S^d\Sigma\to S^d\Sigma$ such that the fixed points of $\Phi$ in the
$p_h$-sector are precisely the $D_{\alpha}$ for \emph{admissible}
$\alpha$; the non-admissible fixed points disappear on this
perturbation to a smooth map.  For instance, if $\phi$ is given in
local holomorphic coordinates near one of its hyperbolic fixed
points by $x+iy\mapsto \lambda x+i\lambda^{-1}y$ where $\lambda > 1$
(ie, $z\mapsto az+\sqrt{a^2-1}\bar{z}$ where
$a=(\lambda+\lambda^{-1})/2 > 1$), then in terms of the natural
holomorphic coordinates $\sigma_1=z_1+z_2$, $\sigma_2=z_1z_2$ near
$\{0,0\}$ on $S^2\Sigma$, $S^2\phi$ is given by \begin{align*}
(\sigma_1,&\sigma_2)\mapsto \left(a\sigma_1+\sqrt{a^2-1}\overline{\sigma_1},a^2\sigma_2+(a^2-1)\overline{\sigma_2}+a\sqrt{a^2-1}(z_1\overline{z_2}+\overline{z_1}z_2)\right)\\
&=\left(a\sigma_1+\sqrt{a^2-1}\overline{\sigma_1},a^2\sigma_2+(a^2-1)\overline{\sigma_2}+\frac{a}{2}\sqrt{a^2-1}(|\sigma_1|^2-|\sigma_{1}^{2}-4\sigma_2|)\right).
\end{align*}
One easily checks that leaving the first component of this function unchanged and adding an appropriate small imaginary-valued function
supported near the origin to the second component results in a smooth function with no fixed points in the coordinate neighborhood
under consideration.

With this in mind, we state a basic property of our groups $HF$
which suggests a connection to $HP$.

\begin{theorem} \label{pfhmain} For any symplectomorphism $\phi$,
$HF(Y_{\phi},f,h,[\omega_{\phi}],\tau)$ arises as the homology of a
chain complex whose generators are the fixed points in the
$p_h$-sector of a smooth map $\Phi_{\tau}\co S^d\Sigma\to
S^d\Sigma$, where $\Phi_{\tau}\to S^d\phi$ in $C^0$-norm as
$\tau\to\infty$.\end{theorem}

Note that the fixed points of $\Phi_{\tau}$ will, for large $\tau$,
all be close to fixed points of $S^d\phi$; one would like to
conclude that they will all be close to the fixed points coming from
the admissible orbit sets that generate $HP$, but it is not clear
that this is the case.  We hope that further analysis of the maps
$\Phi_{\tau}$ might make it possible to establish a correspondence
between the generators and flowlines for $HP$ and those for $HF$
when $\tau$ is large enough and hence to equate the two groups, but
this seems out of reach at present.

While we cannot go so far as to prove that $HF$ and $HP$ are
equivalent, our results do suffice to imply the existence of
periodic points with certain periods for a certain class of surface
symplectomorphisms. Recall from Seidel \cite{Sei} that a symplectomorphism
$\phi\co\Sigma\to\Sigma$ is called \emph{monotone} provided that,
where $Y_{\phi}$ is the mapping torus of $\phi$ and $\omega_{\phi}$
is the form on $Y_{\phi}$ induced by the symplectic form on
$\Sigma$, the cohomology classes $e(T^{vt}Y_{\phi})$ and
[$\omega_{\phi}]$ are proportional.  There are monotone
symplectomorphisms in every mapping class, and if $\phi$ is monotone
then so is any $\psi\circ\phi$ where $\psi$ is the flow of any
possibly-time-dependent vector field symplectically dual to a
(possibly-time-dependent) closed one-form representing an element of
the image of $id-\phi^*\co H^1(\Sigma;\mathbb{R})\to
H^1(\Sigma;\mathbb{R})$, so in particular the intersection of the
space of monotone symplectomorphisms of $\Sigma$ with any component
of the space of orientation-preserving diffeomorphisms of $\Sigma$
is infinite-dimensional.  The following (which is a special case of
known results previously established by rather different means) is
the most general statement we can make about the dynamics of such
symplectomorphisms. \begin{cor}\label{fixedpoint}  Let
$\phi\co\Sigma\to\Sigma$ be a monotone symplectomorphism of a
surface of genus $g\geq 2$.  Then the induced map $S^{2g-2}\phi\co
S^{2g-2}\Sigma\to S^{2g-2}\Sigma$ on the $(2g-2)$th symmetric
product of $\Sigma$ has a fixed point.
\end{cor}

As will be clear from our proof, for $\phi$ belonging to particular
mapping classes the number $2g-2$ can be lowered depending on the
properties of the Seiberg--Witten basic classes of the total spaces
of Lefschetz fibrations having $\phi$ as monodromy (in particular,
if $\phi$ is the monodromy around a loop of a Lefschetz fibration
obtained via Donaldson's construction by blowing up a Lefschetz
pencil, the basic class corresponding to a section of square $-1$
forces $\phi$ itself to have a fixed point). This connection seems
to deserve further study.

Note that a fixed point of $S^{2g-2}\phi$ is equivalent to, for some
partition $2g-2=\sum_{i=1}^{m} n_id_i$, periodic orbits
$o_1,\ldots,o_m$ of $\phi$ with minimal periods $d_1,\ldots,d_m$
respectively.  For $g>2$, \fullref{fixedpoint} can also be
deduced via elementary methods: by considering the relationship of
(what are now called) the Lefschetz numbers of the iterates of
$\phi$ to the characteristic polynomial of the action of $\phi$ on
$H^1(\Sigma;\mathbb{Z})$, Nielsen showed in \cite{N} that an
orientation-preserving homeomorphism $\phi$ must have a periodic
point of period at most $2g-2$ and that this estimate is best
possible; by examining his argument more carefully one can show that
it implies that one of the Lefschetz numbers $L(\phi),L(S^2\phi),
L(S^{2g-2}\phi)$ is nonzero, so that in any event $S^{2g-2}\phi$ has
a fixed point.  Since Nielsen's argument does not work for the case
$g=2$, he asked in \cite{N} whether orientation-preserving
homeomorphisms of surfaces of genus 2 always have points of period
at most 2; this question remained open for decades before eventually
being answered affirmatively by Dicks and Llibre \cite{DL}, using
methods quite different from those we use in the special case
considered here.

\subsection{Remarks
and Acknowledgements}\label{sec.1.5} 
During  the development of this work, I
learned of the thesis \cite{P} of T Perutz, which (among several
other things) develops ideas parallel to some of those discussed in
this paper. By and large, due to slight technical differences in our
formulations, the proofs given here are independent of similar
results in \cite{P} (a notable exception is that I appeal to
Perutz's calculation of a certain cohomology class, which is needed
in the proof of \fullref{coeffs} (which has no analogue in
\cite{P}) and simplifies some other arguments); this leads to some
redundancy between some results in \cite{P} and parts of Sections 4
and 5 of this paper, which seems justified because it ensures that
the relevant results are in the particular form that we need, and
hopefully makes this paper more readable than it otherwise would be.
I am very grateful to T Perutz for sending me his thesis and for
some interesting conversations. Thanks are also due to C Taubes for
advice regarding vortices which played an essential role in the
proof of \fullref{pfhmain}; to M Hutchings, Y-J Lee, and T
Mrowka for answering some questions regarding their work; and to the
anonymous referee for his or her detailed comments and useful
suggestions. This work was partially supported by an NSF
postdoctoral fellowship.

\section{Defining $HF$}\label{Define}

Let $Y,f,h,c,\tau$ be as in the previous section.  Our groups
$HF(Y,f,h,c,\tau)$ are defined, adapting Salamon \cite{Sal}, as the
Floer homologies of certain symplectomorphisms that the monodromy of
$f$ induces on the symmetric products of the fibers of $f$.  Note that
whereas \cite{Sal} begins with an explicit presentation of $Y$ as a
mapping torus, we do not begin with such data; since the fibration
$f\co Y\to S^1$ only specifies the monodromy of $f$ up to smooth
isotopy, we will need additional data to describe its symplectic
behavior. The triple $(h,c,\tau)$ (principally, just $c$) provides
these data: recall that $h\in H_1(X;\mathbb{Z})$ was arbitrary; $c\in
H^2(X;\mathbb{R})$ was a class having positive pairing with the fiber
of $f$, and $\tau>2\pi h\cap [\fiber]$.

 Set $d=h\cap
[\fiber]$. First, if $d<0$, set $HF(Y,f,h,c,\tau;A)=0$; this is
consistent with adjunction relations in monopole Floer theory.  If
$d=0$, set $HF(Y,f,h,c,\tau;A)=A$ if $h=0$ and
$HF(Y,f,h,c,\tau;A)=0$ otherwise. Restrict attention now to the case
$d>0$.

We then have $\langle c+\frac{2\pi}{\tau}PD(h),[\fiber]\rangle>0$,
which enables us to use the Thurston trick to find closed  forms
$\omega=\omega_{c,h,\tau}\in \Omega^2(Y)$ representing
$c+\frac{2\pi}{\tau}PD(h)$ which restrict symplectically to every
fiber of $f$.  The $\omega$--orthogonal complement of $\ker
f_{*}\subset TY$ then defines a horizontal subspace of $TY$, and so
picking a basepoint in $S^1$ and flowing along the horizontal lift
of the vector field $\partial_{\theta}$ on $S^1$ defines a
symplectomorphism $\phi_{\omega}\co\Sigma\to \Sigma$, where $\Sigma$
denotes the fiber of $f\co Y\to S^1$ over the chosen basepoint in
$S^1$.  We now show that the Hamiltonian isotopy class of the
symplectomorphism $\phi_{\omega}$ depends only on the class
$[\omega]=c+\frac{2\pi}{\tau}PD(h)\in H^2(Y;\mathbb{R})$. (This fact
is assuredly well-known and we indicate the proof only for
completeness.)
\begin{prop}\label{cvx} Let $f\co Y\to S^1$ be a fibration, and
$\omega_0,\omega_1$ two closed forms on $Y$ which restrict to each
fiber as volume forms and represent the same class in
$H^2(Y;\mathbb{R})$.  Then $\{\phi_{s\omega_1+(1-s)\omega_0}\}_{s\in
[0,1]}$ is a Hamiltonian isotopy from $\phi_{\omega_0}$ to
$\phi_{\omega_1}$.
\end{prop}
\begin{proof} Note first that, for any $t\in S^1$, $\omega_0$ and
$\omega_1$ are nonvanishing 2--forms on the 2--manifold $f^{-1}(t)$
which induce the same orientation, so the same statement applies to
$\omega_s:=s\omega_1+(1-s)\omega_0$ for each $s$; in particular the
symplectomorphism $\{\phi_{s\omega_1+(1-s)\omega_0}\}$ is
well-defined.

Since $\omega_0$ and $\omega_1$ are cohomologous, write
$\omega_1=\omega_0+d\alpha$.  Now consider the fibration
\[\pi\co [0,1]\times Y\to [0,1]\times S^1\]\[(s,y)\mapsto
(s,f(y));\] the form $\Omega=\omega_0+d(s\alpha)$ is then closed and
restricts to each $Y=\pi^{-1}(\{s\}\times S^1)$ as $\omega_s$. If we
let $\gamma_s$ be the loop in $[0,1]\times S^1$ obtained by
juxtaposing the paths (each with domain $[0,1]$ for the parameter $t$)
$t\mapsto (0,e^{2\pi it})$, $t\mapsto (st,1)$, $t\mapsto (s,e^{-2\pi
it})$, and $t\mapsto (s(1-t),1)$, then the monodromy around $\gamma_s$
(using the horizontal lift of $\dot{\gamma_s}$ given by the
$\Omega$--orthogonal complement of $\ker d\pi$) is
$\phi_{\omega_s}\circ \phi_{\omega_0}^{-1}$ (modulo the identification
of $\pi^{-1}(0,1)$ with $\pi^{-1}(s,1)$ via horizontal translation).
But according to Proposition 6.31 of McDuff and Salamon \cite{McDS},
since $\Omega\in \Omega^2([0,s]\times Y)$ is closed and $\gamma_s$ is
contractible, the $\Omega$--monodromy around $\gamma_s$ is
Hamiltonian.  Thus each $\phi_s$ differs from $\phi_0$ by a
Hamiltonian isotopy, and the proposition is proven. \end{proof} Thus,
specifying the cohomology class $c+\frac{2\pi}{\tau}PD(h)\in
H^2(Y;\mathbb{R})$ specifies the monodromy $\phi_{\omega}$ of the
fibration $f\co Y\to S^1$ up to Hamiltonian symplectomorphism; this
makes it reasonable to expect basic symplectic properties (eg, Floer
homology) of $\phi$ to depend only on $c$, $h$, and $\tau$.  The
reader might wonder at this point why we are using the formula
$c+\frac{2\pi}{\tau}PD(h)$ to refer to the cohomology class of the
form on $Y$ rather than just, say, $c$ (as is done in \cite{P}); the
reason is that with this choice we ensure that the coefficient ring
$\tilde{\Lambda}_{h,c}$ over which the still-to-be defined group
$HF(Y,f,h,c,\tau)$ is naturally a module will depend only on $h$ and
$c$ and not on $\tau$ (indeed, as will be seen in the proof of
\fullref{coeffs}, under a technical assumption on $h\cap [\fiber]$ this
choice ensures that the groups themselves are independent of $\tau$).

We now explain the definition of the groups $HF$ when $d>0$.   Given
the data $(Y,f,h,c,\tau)$, we choose a closed, fiberwise symplectic
form $\omega$ representing the class $c+\frac{2\pi}{\tau}PD(h).$  As
above, choosing a basepoint $1\in S^1$ and using the horizontal
distribution induced by $\omega$, we obtain a fiberwise
symplectomorphism \[(Y,\omega)\cong \left(\frac{\mathbb{R}\times
\Sigma}{(t+1,x)\sim (t,\phi_{\omega}(x))},\omega\right),\] where
$\Sigma=f^{-1}(1)$ and the symplectomorphism $\phi_{\omega}\co
(\Sigma,\omega|_{\Sigma})\to (\Sigma,\omega|_{\Sigma})$ is the
monodromy of $f$.
 We shall follow \cite{Sal} to define a family of
symplectic forms $\Omega_{d,\omega,\tau}$ on the symmetric products
$S^d\Sigma$, and a family of maps $\Phi_{d,\omega,\tau}\co
S^d\Sigma\to S^d\Sigma$ which are
$\Omega_{d,\omega,\tau}$--symplectomorphisms.

 These symplectomorphisms are obtained from a
construction involving the \emph{$U(1)$ vortex equations} on the
fiber $\Sigma$.  This construction is carried out in Sections 4 and
5 of \cite{Sal}; we recall it here. Let $L$ be a Hermitian line
bundle on $\Sigma$ of degree $d$. The configuration space
$\mathcal{C}_{L}$ is defined as the space of pairs $(A,\theta)$
where $A$ is a Hermitian connection on $L$ and $\theta$ is a section
of $L$ (in practice these should be viewed as elements of Sobolev
spaces $L^{2}_{k}$ for $k>2$, but for the most part we shall
suppress these standard details; incidentally, our $\theta$ will be
$(2\tau)^{-1/2}$ times the section $\Theta_0$ used in \cite{Sal}).
$J$ will denote an $\omega$--compatible complex structure on
$\Sigma$. The $U(1)$ vortex equations for a pair
$(A,\theta)\in\mathcal{C}_L$ are then
\begin{align}\label{vortex} \bar{\partial}_{J,A}\theta&=0\nonumber \\ iF_A&=\tau(1-|\theta|^2)\omega.\end{align}

Now the tangent space to $\mathcal{C}_L$ is given by
$T_{(A,\theta)}\mathcal{C}_L=\Omega^1(\Sigma;i\mathbb{R})\times
\Omega^0(E)$.  There exists a \emph{universal symplectic form}
$\tilde{\Theta}\in \Omega^2(\mathcal{C}_L)$ given by \[
\tilde{\Theta}((\alpha_1,\theta_1),(\alpha_2,\theta_2))=-\int_{\Sigma}
\alpha_1\wedge\alpha_2+2\tau \int_{\Sigma}\Im \langle
\theta_1,\theta_2\rangle \omega.\]  Where again $J$ is an arbitrary
complex structure on $\Sigma$ compatible with $\omega$, \[
\mathcal{X}_J=\{(A,\theta)\in
\mathcal{C}_L|\bar{\partial}_{J,A}\theta=0\}\] is a symplectic
submanifold of $(\mathcal{C}_L,\tilde{\Theta})$.  A K\"ahler
structure is induced on both $\mathcal{X}_J$ and $\mathcal{C}_L$ by
$\tilde{\Theta}$ together with the complex structure
$(\alpha,\theta)\mapsto (\ast \alpha,i\theta)$.

One has an action of the gauge group
$\mathcal{G}=L^{2}_{k+1}(\Sigma,S^1)$ by \[ u\cdot
(A,\theta)=(A-u^{-1}du,u\theta);\] in fact this is a Hamiltonian
action on $\mathcal{C}_L$, with moment map \[ \mu\co
(A,\theta)\mapsto \ast iF_A+\tau|\theta|^2 \] (where $\ast
\omega=1$). The set $\tilde{\mathcal{M}}_{\Sigma,d}(J,\tau)$ of
solutions to the vortex equations is thus $\mathcal{X}_J\cap
\mu^{-1}(\tau)$, and the set \[
\mathcal{M}_{\Sigma,d}(J,\tau)=\tilde{\mathcal{M}}_{\Sigma,d}(J,\tau)/\mathcal{G}\]
of gauge equivalence classes of solutions to the vortex equations is
the symplectic reduction of $\mathcal{X}_J$ by the action of
$\mathcal{G}$.

An extension by O Garcia-Prada \cite{GP} of a theorem of C Taubes
implies that the map
\[\mathcal{M}_{\Sigma,d}(J,\tau)\to S^d(\Sigma,J)\] which
sends an equivalence class $[A,\theta]$ to the zero set of $\theta$
is an isomorphism of complex manifolds. (We write $S^d(\Sigma,J)$
here to emphasize that the complex structure, and indeed even the
$C^{\infty}$ charts, on $S^d\Sigma$ depend on the complex structure
on $\Sigma$.) The form $\tilde{\Theta}$ on $\mathcal{C}_L$ descends
to a symplectic form $\Theta_{J,\tau}$ on each
$\mathcal{M}_{\Sigma,d}(J,\tau)$.

Let $\mathcal{J}(\Sigma)$ be the space of (almost) complex
structures on $\Sigma$, and consider the space \[
\tilde{\mathcal{X}}=\{(J,A,\theta)|J\in\mathcal{J}(\Sigma),
(A,\theta)\in \mathcal{X}_J\}.\]  $\tilde{\mathcal{X}}$ obviously
fibers over $\mathcal{J}(\Sigma)$ with fiber $\mathcal{X}_J$, and
$\tilde{\Theta}$ defines a closed, fiberwise symplectic form on this
fibration.  Carrying out the above symplectic reduction process
fiberwise, we obtain a closed, fiberwise symplectic form
$\hat{\Theta}_{\tau}$ on the fibration \[
\hat{\mathcal{X}}=\{(J,D)|D\in S^d(\Sigma,J)\}\to
\mathcal{J}(\Sigma).\]  \begin{align*} \gamma\co [0,1]\tag*{\hbox{If}}&\to
\mathcal{J}(\Sigma)\\ t&\mapsto J_t\end{align*} is a smooth path of
almost complex structures on $\Sigma$, we then obtain a symplectic
fibration on $\gamma^{*}\hat{\mathcal{X}}$ with closed fiberwise
symplectic $2$--form $\gamma^{*}\hat{\Theta}_{\tau}$.  Using the
parallel translation given by the
$(\gamma^{*}\hat{\Theta}_{\tau})$--orthogonal complement of
$T^{vt}\gamma^*\hat{\mathcal{X}}$ then gives a symplectomorphism
\begin{equation} \label{transport} F_{\{J_t\}}\co (S^{d}(\Sigma,J_0),\Theta_{J_0,\tau})\to
(S^{d}(\Sigma,J_1),\Theta_{J_1,\tau}).\end{equation}

With this preparation in hand, let us return to the data consisting
of a fibration $f\co Y\to S^1$ along with $h\in H_1(X;\mathbb{Z})$,
$c\in H^2(X;\mathbb{R})$, and $\tau\in \mathbb{R}$; as earlier, we
set $d=h\cap(\fiber)$ and choose a fiberwise symplectic $2$--form
$\omega=\omega_{c,h,\tau}$ representing $c+\frac{2\pi}{\tau}PD(h)$;
letting $\phi_{\omega}$ be the resulting monodromy map $\Sigma\to
\Sigma$, the pair $(f\co Y\to S^1,\omega)$ is identified as a
symplectic fibration with the mapping torus of $\phi_{\omega}$.  Now
let $\bar{J}$ be an $\omega$--compatible almost complex structure on
$T^{vt}Y$; the datum of $\bar{J}$ amounts to a path $J_t$ of almost
complex structures such that $J_1=\phi_{\omega}^{*}J_0$.  Now if
$\phi$ is any diffeomorphism of $\Sigma$ and $J$ any complex
structure on $\Sigma$ there is a tautological K\"ahler isomorphism
\[ S^d\phi\co S^d(\Sigma,\phi^*J)\to S^d(\Sigma,J)\]
(we should caution here that if we were to instead view $S^d\phi$ as
a map $S^d(\Sigma,J)\to S^d(\Sigma,J)$ it typically would not even
be differentiable).  We then set
\[ \Phi_{d,\omega,\tau}=S^d\phi_{\omega}\circ F_{\{J_t\}}\co
(S^{d}(\Sigma,J_0),\Theta_{J_0,\tau})\to
(S^{d}(\Sigma,J_0),\Theta_{J_0,\tau}).\] $\Phi_{d,\omega,\tau}$ is
the composition of two symplectomorphisms and hence is a
symplectomorphism.

\begin{prop}\label{Ham}  If $\omega_0$ and $\omega_1$ are two representatives of
$c+\frac{2\pi}{\tau}PD(h)$ as above, and $J_{0,t}$ (resp.\ $J_{1,t}$)
are $\omega_0$-- (resp.\ $\omega_1$--) compatible almost complex
structures, with $J_{s,t}$ a family of
$(s\omega_1+(1-s)\omega_0)$--compatible almost complex structures
interpolating between them, then $\Phi_{d,\omega_0,\tau}$ is
Hamiltonian equivalent to
$F_{\{J_{s,0}\}}^{-1}\circ\Phi_{d,\omega_1,\tau}\circ
F_{\{J_{s,0}\}}$.\end{prop}

\begin{proof}  First note that, by the proof of \fullref{cvx}, the forms $\omega_s=s\omega_1+(1-s)\omega_0$ are all
fiberwise symplectic on $Y$; in light of this the family of almost
complex structures $J_{s,t}$ exists by the contractibility of
$\mathcal{J}(\Sigma)$.  Defining $\Gamma\co [0,1]^2\to
\mathcal{J}(\Sigma)$, we have a closed fiberwise symplectic
$2$--form $\Gamma^{*}\hat{\Theta}$ on the fibration \[
\Gamma^*\mathcal{\hat{X}}=\{(s,t,D)|D\in S^d(\Sigma,J_{s,t})\}\to
[0,1]^2,\] and since each
\[ S^d\phi_{\omega_s}\co (S^d(\Sigma,J_{s,1}),\Theta_{J_{s,1},\tau})\to
(S^d(\Sigma,J_{s,0}),\Theta_{J_{s,0},\tau})\] is a
symplectomorphism, $\Gamma^*\hat{\Theta}$ descends to the mapping
torus \[ \Gamma^*\mathcal{\hat{X}}/(s,1,D)\sim (s,0,S^d\phi (D)),\]
making this fibration over the cylinder a locally Hamiltonian
fibration. Now the monodromy of this fibration around a loop similar
to that in the proof of \fullref{cvx} (using, as usual, the
form induced by $\Gamma^*\hat{\Theta}$ to determine the horizontal
distribution) is
$F_{\{J_{s,0}\}}^{-1}\circ\Phi_{d,\omega_1,\tau}\circ
F_{\{J_{s,0}\}}\circ \Phi_{d,\omega_0,\tau}^{-1}$.  But since the
loop is contractible this monodromy is Hamiltonian by 
\cite[Proposition 6.31]{McDS}. \end{proof}

We can now finally define $HF(Y,f,h,c,\tau)$. If $\Phi\co
(X,\omega)\to (X,\omega)$ is a symplectomorphism of a symplectic
manifold $X$ with nondegenerate fixed points, let $Y_{\Phi}\to S^1$
denote the mapping torus of $\Phi$ and $\omega_{\Phi}$ denote the
$2$--form on $Y_{\Phi}$ induced by the pullback of $\omega$ to
$\mathbb{R}\times Y$. Recall then that, for $\mathcal{P}\in \pi_1
(\Gamma(Y_{\Phi}))$, the Floer homology
$HF^{\symp}(\Phi,\mathcal{P})$ is obtained, na\"ively, by
Floer--Morse theory on the subset of $\Gamma(Y_{\Phi})$ consisting
of sections representing the homotopy class $\mathcal{P}$ using the
action $1$--form \begin{equation}\label{actionform}
\mathcal{Y}_{\gamma}(\xi)=-\int_{0}^{1}\omega_{\Phi}(\dot{\gamma}(t),\xi(t))\,dt,\end{equation}
for $\gamma\in\mathcal{P}$ and $\xi\in
T_{\gamma}\mathcal{P}=\gamma^*T^{vt}Y_{\phi}$ (if $\Phi$ has
degenerate fixed points, a perturbation is used; see Section 3 of
\cite{L1}, which is the most thorough reference for this subject).
The complex $CF^{\symp}(\Phi,\mathcal{P})$ is then generated by fixed
points $x$ of $\Phi$ with the property that the ``constant section''
of $Y_{\Phi}$ at $x$ belongs to the homotopy class $\mathcal{P}$;
the boundary operator counts holomorphic cylinders in
$\mathbb{R}\times Y_{\Phi}$.  $CF^{\symp}$ naturally has its
coefficients in, depending on convention, either the Novikov ring \[
\Nov\left(\frac{\ker\langle
c_1(T^{vt}Y_{\phi}),\cdot\rangle}{\ker\langle
c_1(T^{vt}Y_{\phi}),\cdot\rangle\cap\ker\langle
\omega,\cdot\rangle},\langle\omega,\cdot\rangle;R\right)\] or the
larger Novikov ring \[ \Nov\left(\ker\langle
c_1(T^{vt}Y_{\phi}),\cdot\rangle,\langle\omega,\cdot\rangle;R\right)\]
for a ring $R$ (in our context, for $g/2+1<d<g-1$, the virtual
moduli methods of Liu and Tian \cite{LT} are required, and so $R$ will need to be
a field of characteristic zero).

We consider the case $X=S^d(\Sigma,J)$ and
$\Phi=\Phi_{d,\omega,\tau}$.  As is shown in 
\cite[Section 7]{Sal}, there is a one-to-one correspondence
$\mathcal{P}_h\leftrightarrow h$ between
\[ \pi_1(\Gamma(Y_{\Phi_{d,\omega,\tau}}))\mbox{ and } \{h\in
H_1(Y_{\phi_{\omega}};\mathbb{Z})|h\cap(\fiber)=d\}\] (roughly
speaking, a homotopy class $\mathcal{P}$ of sections of the bundle
$Y_{\Phi_{d,\omega,\tau}}\to S^1$ of symmetric products corresponds
to the homology class represented by the union of points appearing
in the divisors represented by some section in $\mathcal{P}$; we'll
be clearer about this later).

As such, we can set \[
HF(Y,f,h,c,\tau)=HF^{\symp}(\Phi_{d,\omega,\tau},\mathcal{P}_h),\]
where, once again, $\omega$ is a fiberwise symplectic representative
of $c+\frac{2\pi}{\tau}PD(h)$.  By \fullref{Ham} and the
standard fact that $HF^{\symp}$ is invariant under conjugation by
symplectomorphisms and under Hamiltonian isotopy, we see immediately
that $HF(Y,f,h,c,\tau)$ does not depend on the choice of $\omega$.
We shall soon verify that the Novikov ring over which it is defined
is as promised in the introduction, and that it is independent of
$\tau$ at least for $d$ outside a certain range.  First, however, a
digression regarding the topology of $Y_{\Phi}$ is in order.

\section{Basic properties of $HF$}\label{sec.3}
\subsection{Topology of
(relative) symmetric products} We review here some basic facts
regarding the cohomology of symmetric products.  A standard
reference for some of this material is Macdonald \cite{Mac}; for the relative
versions see also the appendices and Section 2.1 of \cite{P}.

First note that if $\Sigma$ is a Riemann surface we obtain a natural
map \[ \uparrow\co H^*(\Sigma;\mathbb{Z})\to
H^*(S^d\Sigma;\mathbb{Z})\] as follows.  Inside the product
$\Sigma\times S^d\Sigma$ we have a divisor $\mathcal{D}=\{(p,D)\in
\Sigma\times S^d\Sigma|p\in D\}$, where we view $D\in S^d\Sigma$ as
a set of points in $\Sigma$ (in other references, such as \cite{P},
$\mathcal{D}$ is called $\Delta$, but we prefer to use $\Delta$ to
denote the diagonal stratum in the  symmetric product).  Letting
$\pi_1$ and $\pi_2$ denote the projections of $\Sigma\times
S^d\Sigma$ onto either factor, the map $\uparrow$ is defined by \[
\uparrow\! c=(\pi_2)_!(\frak{d}\cup\pi_{1}^{*}c),\] where
$\frak{d}=PD[\mathcal{D}]\in H^2(\Sigma\times
S^d\Sigma;\mathbb{Z})$.  Note the simple geometric interpretation of
this map: if $c\in H^*(\Sigma;\mathbb{Z})$, let $A$ be a cycle
Poincar\'e dual to $c$; then the Poincar\'e dual of $\uparrow\! c$
is represented by a cycle whose image is the set $\{D\in
S^d\Sigma|D\cap A\neq \varnothing\}$ of degree $d$ effective
divisors on $\Sigma$ which contain a point of $A$.

Dually, there is a map \begin{align*} \downarrow\co
H_*(S^d\Sigma;\mathbb{Z})&\to H_{*}(\Sigma;\mathbb{Z})\\
A&\mapsto (\pi_1)_*(\frak{d}\cap \pi_{2}^{!}A);\end{align*} again we
may intuitively visualize this as sending $A\in
H_*(S^d\Sigma;\mathbb{Z})$ to the homology class in $\Sigma$
represented by union of all the points in $\Sigma$ which appear in
the set of divisors which is the image of some chain representing
$A$.  Note that for $c\in H^*(\Sigma;\mathbb{Z})$ and $A\in
H_*(S^d\Sigma;\mathbb{Z})$, we have \[ \langle\uparrow\!
c,A\rangle=\langle
(\pi_2)_!(\frak{d}\cup\pi_{1}^{*}c),A\rangle=\langle
c,(\pi_1)_*(\frak{d}\cap \pi_{2}^{!}A)\rangle=\langle c,\downarrow
\! A\rangle.\]

We can now describe the cohomology of $S^d\Sigma$; for proofs (in a
somewhat different language) see \cite{Mac}.
\begin{prop}\label{basicsd}\begin{itemize}\item[(i)] The map $\uparrow\co
H^1(\Sigma;\mathbb{Z})\to H^1(S^d\Sigma;\mathbb{Z})$ is an
isomorphism. \item[(ii)] Identifying $H^1(\Sigma;\mathbb{Z})$ with
$H^1(S^d\Sigma;\mathbb{Z})$ by $\uparrow$, where $\omega$ is a
positive generator of $H^2(\Sigma;\mathbb{Z})$ and
$U=\uparrow\omega$, one has \[
H^{*}(S^d\Sigma;\mathbb{Z})=\frac{\mathbb{Z}[U]\otimes_{\mathbb{Z}}\Lambda^{*}H^1(\Sigma;\mathbb{Z})}{\langle
U^{i}\otimes(\gamma_1\wedge\cdots\wedge\gamma_j)|i,j\geq 0,\,i+j>d,
\gamma_i\in H^1(\Sigma;\mathbb{Z}))\rangle}\] as graded rings.  In
particular $H^*(S^d\Sigma;\mathbb{Z})$ is naturally a module over
$\mathbb{Z}[U]\otimes_{\mathbb{Z}}\Lambda^{*}H^1(\Sigma;\mathbb{Z})$.\end{itemize}\end{prop}

We will also be interested in \emph{relative} symmetric products
associated to surface fibrations.  To wit, let $\pi\co X\to B$ be a
fibration with fiber $\Sigma$ a closed surface of genus $g\geq 2$
($B$ is a compact manifold, possibly with boundary, in which case
$\partial X=\pi^{-1}(\partial B)$). By choosing an almost complex
structure $J$ on $T^{vt}X$ and appealing to the parametrized Riemann
mapping theorem to obtain ``restricted charts'' which are smooth
horizontally and holomorphic vertically (see, eg,
\cite{DS}), one can construct a fibration \[ \Pi\co X_d(\pi)\to B\]
carrying a vertical almost complex structure $\tilde{J}$ such that
each fiber $\Pi^{-1}(b)$ is identified as a complex manifold with
$(S^d\pi^{-1}(b),J|_{\pi^{-1}(b)})$.  The maps $\uparrow$ and
$\downarrow$ extend to the relative context: inside the fiber
product \[ X_{\,\pi}\!\times_{\Pi}X_d(\pi) \] we have a
codimension--2 submanifold \[ \mathcal{D}=\{(b,p,D)\in B\times
X\times X_d(\pi)|p\in \pi^{-1}(b),D\in\Pi^{-1}(b), p\in D\}\]
determining a class \[[\mathcal{D}]\in
H_{2d+2}(X_{\,\pi}\!\times_{\Pi}X_d(\pi),\partial(X_{\,\pi}\!\times_{\Pi}X_d(\pi)
)).\] Again let $\frak{d}=PD[\mathcal{D}]$ and define
\begin{align*}
\uparrow\co H^*(X;\mathbb{Z})&\to H^*(X_d(\pi);\mathbb{Z})\\
c&\mapsto (\pi_2)_!(\frak{d}\cup\pi_{1}^{*}c),\end{align*} and
similarly for $\downarrow$, where $\pi_1$ and $\pi_2$ are the
projections from $X_{\,\pi}\!\times_{\Pi}X_d(\pi)$ to $X$ and
$X_d(\pi)$, respectively.

The following formula, proven using the family Atiyah--Singer
theorem, expresses the first Chern class of the vertical tangent
bundle $T^{vt}X_{d}(\pi)$ (with its induced almost complex
structure) in terms of the Euler class of $T^{vt}X$.
\begin{lemma}\label{euler}{\rm\cite[Lemma 2.1.1]{P}}\qua
Assume that $B$ is closed.  Then
\[c_1(T^{vt}X_{d}(\pi))=\frac{1}{2}\left(\uparrow\!(e(T^{vt}X))+(\pi_2)_{!}(\frak{d}^2)\right).\]
\end{lemma}

We now give a geometric interpretation of the class
$(\pi_2)_{!}(\frak{d}^2)$ appearing in \fullref{euler}.
$\Delta\subset X_d(f)$ will denote the real--codimension 2
subvariety of $X_d(f)$ consisting of divisors $D$ having one or more
points repeated.  $\Delta$ is easily seen to represent an element
$[\Delta]\in H_{2r-2+\dim B}(X_d(f),\partial X_d(f))$.

\begin{prop}\label{dtw}  If $B$ is closed, then
 $(\pi_2)_{!}(\frak{d}^2)=PD[\Delta]+\uparrow\!(e(T^{vt}X)).$
\end{prop}
\begin{proof}
Let $v\in \Gamma(T^{vt}X)$ be a transversally-vanishing
vertically-valued vector field on $X$. $v$ then induces a
vertically-valued vector field $V$ on $X_{\pi}\times_{\Pi}X_d(\pi)$
by, with respect to the splitting
$T_{(b,p,D)}^{vt}X_{\pi}\times_{\Pi}X_d(\pi)=T_{(b,p,D)}^{vt}X\oplus
T_{(b,p,D)}^{vt}X_d(\pi)$, setting $V(b,p,D)=(v(b,p),0)$.  So
$\frak{d}^2$ is  represented by the Poincar\'e dual of the vanishing
locus $\mathcal{V}$ of the projection of $V|_{\mathcal{D}}$ to the
normal bundle $N\mathcal{D}$. This latter is easily seen to
be\begin{align*}
 \{(b,p,D)\in X_{\,\pi}\times_{\Pi}X_d(f)|D=2p+D'\mbox{
for some } D'\in X_{d-2}(\pi)\}\\ \cup\{(b,p,D)\in
X_{\,\pi}\times_{\Pi}X_d(f)|v(p)=0\}.\end{align*}  Hence
$(\pi_2)_{!}(\frak{d}^2)=PD(\pi_2)_*[\mathcal{V}]$ is represented by
the Poincar\'e dual of the homology class represented by the union
\[\Delta\cup\{D\in X_d(f)|D\cap v^{-1}(0)\neq\varnothing.\}\]  But
the second set in this union represents $PD(\uparrow\!(e(T^{vt}X)))$
since $v$ is a transversally vanishing section of $T^{vt}X$.  The
conclusion is then immediate.\end{proof}
\subsection{Novikov rings}  We shall now recall how
exactly Novikov rings enter into the general picture of symplectic
Floer theory; \cite{L1} is a good reference for those seeking
further details.  After this, we shall be prepared to verify that
the Novikov ring over which $HF(Y,f,h,c,\tau)$ is defined is indeed
the ring $\tilde{\Lambda}_{h,c}$ of the introduction.  As many
references exist for technicalities relating to the relevant
Fredholm theory and compactness, our treatment shall be essentially
formal, but our conclusion will remain valid in full generality,
with the exception that the ring $R$ below needs to be a field of
characteristic zero when the virtual moduli methods of \cite{LT} are
required. Returning to the notation of \fullref{Define}, assume
that $\Phi:(X,\omega)\to (X,\omega)$ is a symplectomorphism. The
configuration space for $HF^{\symp}(\Phi,\mathcal{P})$ is then the
space of sections of the mapping torus $Y_{\Phi}$ of $\Phi$
belonging to the homotopy class $\mathcal{P}$, and generators for
the Floer complex will be those sections of $Y_{\Phi}$ corresponding
to fixed points of $\Phi$.  We have a natural evaluation map \[
ev_{\mathcal{P}}\co H_1(\mathcal{P};\mathbb{Z})\to
H_2(Y_{\Phi};\mathbb{Z}).
\]  Let $A_0\in H_1(\mathcal{P};\mathbb{Z})$ be such that  $\langle c_1(T^{vt}Y_{\Phi}),ev(A_0)\rangle$
is a positive generator for $Im(\langle
c_1(T^{vt}Y_{\Phi}),ev(\cdot)\rangle\co
H_1(\mathcal{P};\mathbb{Z})\to \mathbb{Z})$. Let
$p\co\tilde{\mathcal{P}}\to \mathcal{P}$ be the universal abelian
cover of the configuration space $\mathcal{P}$, with covering group
$H_1(\mathcal{P};\mathbb{Z})$. One then has a natural relative
$\mathbb{Z}$--grading $\tilde{gr}(\tilde{x},\tilde{y})$ for any
nondegenerate $\tilde{x},\tilde{y}\in\tilde{\mathcal{P}}$ given by
the Maslov index, such that
\begin{equation}\label{gr}\tilde{gr}(\tilde{x},\tilde{y})+\tilde{gr}(\tilde{y},\tilde{z})=\tilde{gr}(\tilde{x},\tilde{z})\mbox{
and } \tilde{gr}(\tilde{x},A\cdot\tilde{x})=\langle
 c_1(T^{vt}Y_{\Phi}),ev(A)\rangle.\end{equation}

Where $\mathcal{F}\subset\mathcal{P}$ denotes the set of critical
points for the action $1$--form $\mathcal{Y}$ of \eqref{actionform},
for all $x\in\mathcal{F}$ choose and fix lifts $\tilde{x}$ in
$\tilde{\mathcal{P}}$ with the property that \[
|\tilde{gr}(\tilde{x},\tilde{y})|<\langle
 c_1(T^{vt}Y_{\Phi}),ev(A_0)\rangle\] for each $x,y$, as is possible
 using the second property in \eqref{gr}.

Where $J$ is a generic $\mathbb{R}$--invariant almost complex
structure on $\mathbb{R}\times T^{vt}(Y_{\Phi})$, the boundary
operator for the Floer complex counts finite energy solutions
\begin{align*} U\co \mathbb{R}\times S^1&\to \mathbb{R}\times
Y_{\Phi}\\(s,t)&\mapsto (s,u(s,t))\end{align*} to a perturbed
Cauchy--Riemann equation $\bar{\partial}_J U=X_H(U)$; each of these
is asymptotic as $s\to\pm\infty$ to generators
$x^{\pm}\in\mathcal{F}$, and $s\mapsto u(s,\cdot)$ determines a path
in $\mathcal{P}$ from $x^-$ to $x^+$.  This path then lifts to a
unique path $\tilde{u}$ in $\tilde{\mathcal{P}}$ from $\tilde{x}^-$
to $A\cdot\tilde{x}^+$ for some $A\in H_1(\mathcal{P};\mathbb{Z})$.
Accordingly, given $A$, let $\mathcal{M}_{J,H}(x^-,x^+;A)$ be the
set of those $U$ which, as above, satisfy $\bar{\partial}_J=X_H(U)$
and determine a path $s\mapsto \tilde{u}(s,\cdot)$ in
$\tilde{\mathcal{P}}$ from $\tilde{x}^-$ to $A\cdot\tilde{x}^+$. For
generic Hamiltonian perturbations $X_H$, this will be a manifold of
dimension $\tilde{gr}(\tilde{x}^-,A\cdot\tilde{x}^+)$ with a free
$\mathbb{R}$--action. Note that the choice of $A_0$ determines a
splitting $H_1(\mathcal{P};\mathbb{Z})=\mathbb{Z}A_0\oplus
\ker(\langle c_1(T^{vt}Y),ev(\cdot)\rangle)$; let $p_2$ denote the
projection onto the second summand in this splitting.

The Floer complex $CF^{\symp}(\Phi,\mathcal{P})$ is then the free
module generated by the elements of $\mathcal{P}$ over the Novikov
ring \begin{equation}\label{gennov} \Nov\left(\ker\langle
c_1(T^{vt}Y_{\Phi}),ev(\cdot)\rangle,\langle
[\omega_{\Phi}],ev(\cdot)\rangle;R\right)\end{equation} where $R$ is
an arbitrary ring (the notation is as in the introduction); the
differential is given by the formula \[
\partial
\langle
x^-\rangle=\sum_{\tilde{gr}(\tilde{x}^-,A\cdot\tilde{x}^+)=1}\#(\mathcal{M}(x,y;A)/\mathbb{R})[p_2(A)]\langle
x^+\rangle.\]

Here $\#(\mathcal{M}(x,y;A)/\mathbb{R})$ refers to a signed count of
points in the indicated compact $0$--manifold, using coherent
orientations as in Floer and Hofer \cite{FH}.

Recall now that
$HF(Y,f,h,c,\tau)=HF^{\symp}(\Phi_{d,\omega,\tau},\mathcal{P}_{h})$
where $\omega$ is a fiberwise symplectic $2$--form on $Y$
representing $c+\frac{2\pi}{\tau}PD(h)$ and $d=h\cap(\fiber)$. We can
be somewhat clearer than before about the definition of the homotopy
class $\mathcal{P}_h$ of sections of $Y_{\Phi_{d,\omega,\tau}}$
corresponding to $h$; note that $Y_{\Phi_{d,\omega,\tau}}$ is the
relative symmetric product built from the fibration $f\co Y\to S^1$;
$\mathcal{P}_h$ is characterized by the property that, where \[
ev\co \pi_0(\Gamma(Y_{\Phi_{d,\omega,\tau}}))\to
H_1(Y_{\Phi_{d,\omega,\tau}};\mathbb{Z})\] is the obvious evaluation
map, we have $\downarrow\!(ev(\mathcal{P}_h))=h$, the map
$\downarrow$ having been defined in the last section.

\begin{lemma}\label{deltacorr} Let $\gamma\in H_1(\mathcal{P}_h)$, with
$ev(\gamma)=T\in H_2(Y_{\Phi_{d,\omega,\tau}};\mathbb{Z})$.  Then \[
\langle PD(\Delta),T\rangle_{Y_{\Phi_{d,\omega,\tau}}}=\langle
2PD(h)-e(T^{vt}Y),\downarrow\!T\rangle_{Y}\] where $\Delta$ is the
diagonal in the relative symmetric product
$Y_{\Phi_{d,\omega,\tau}}$.
\end{lemma}
\begin{proof}  Assume $\gamma$ is represented by a loop $\alpha\co
S^1\to \mathcal{P}_h$.  Define \begin{align*} \tilde{\alpha}\co
S^1\times S^1&\to S^1\times
Y_{\Phi_{d,\omega,\tau}}\\(\theta,t)&\mapsto
(\theta,\alpha(\theta,t))\end{align*} and let $\tilde{T}\in
H_2(S^1\times Y_{\Phi_{d,\omega,\tau}};\mathbb{Z})$ be the class
represented by $\tilde{\alpha}$.  By perturbing $\alpha$, we may
assume that $\tilde{\alpha}$ is transverse to the diagonal in
$S^1\times Y_{\Phi_{d,\omega,\tau}}$ (which is a relative symmetric
product over the torus, of course).  Where $i\co Y\cong \{1\}\times
Y\to S^1\times Y$ is the inclusion, we have
\begin{equation}\label{circle}
\downarrow\!\tilde{T}=i_*(\downarrow\!T)+[S^1]\times
h.\end{equation} Let $C=\{p\in S^1\times Y|\exists D\in
Im(\tilde{\alpha})\, p\in D\}$ be the image of the cycle in
$S^1\times Y$ representing $\downarrow\! T$, so that $C$ is ``swept
out'' by $\tilde{\alpha}$. Where $v$ is a transversally vanishing
section of $T^{vt}(S^1\times Y)$ all of whose zeroes occur over
points in $S^1\times S^1$ which are not contained in the finite set
$\{(s,t)\in S^1\times S^1|\tilde{\alpha}(s,t)\in\Delta\}$ and
$\phi^t\co Y\to Y$ is its time--$t$ flow, set \[C^t=\{p\in S^1\times
Y|\phi^t(p)\in C\}.\] For $t$ small, one sees that there is one
intersection of $C$ with $C^t$ for each point of $v^{-1}(0)\cap C$,
and another intersection of $C$ with $C^t$ for each intersection of
the image of $\tilde{\alpha}$ with $\Delta$, and that moreover these
all occur with correctly--corresponding signs.  This shows that \[
\langle
PD(\downarrow\!\tilde{T}),\downarrow\!\tilde{T}\rangle_{S^1\times
Y}=\langle PD(\Delta),\tilde{T}\rangle_{S^1\times
Y_{\Phi_{d,\omega,\tau}}}+\langle e(T^{vt}(S^1\times
Y)),\downarrow\!\tilde{T}\rangle_{S^1\times Y}.\] But from
\eqref{circle} we see \[ \langle
PD(\downarrow\!\tilde{T}),\downarrow\!\tilde{T}\rangle_{S^1\times
Y}=2\langle PD(h),\downarrow\! T\rangle_Y,\] while it is
straightforward to see that
\[
\langle PD(\Delta),\tilde{T}\rangle_{S^1\times
Y_{\Phi_{d,\omega,\tau}}}=\langle
PD(\Delta),T\rangle_{Y_{\Phi_{d,\omega,\tau}}}\] and
\[e(T^{vt}(S^1\times Y),\downarrow\!\tilde{T}\rangle_{S^1\times Y}=\langle
e(T^{vt}Y),\downarrow\!T\rangle_Y,\] proving the lemma.
\end{proof}

We now invoke the following very useful calculation of Perutz.

\begin{theorem}\label{cohclass}{\rm \cite[page 70]{P}}\qua
\[[(\Omega_{d,\omega,\tau})_{\Phi_{d,\omega,\tau}}]=2\pi\left(\tau\uparrow[\omega]-\pi(\pi_2)_{!}(\frak{d}^2)\right)\in
H^2(Y_{\Phi_{d,\omega,\tau}};\mathbb{R}).\]\end{theorem}

With this in hand we can identify the Novikov ring over which
$HF(Y,f,c,h,\tau)$ is defined.

Combining \fullref{dtw} and \fullref{deltacorr} and using
several times the duality between $\uparrow$ and $\downarrow$, we
see that, for $\gamma\in H_1(\mathcal{P}_h;\mathbb{Z})$, we have
\begin{equation}\label{dsquared} \langle(\pi_2)_{!}(\frak{d}^2),ev(\gamma)\rangle=\langle
PD[\Delta]+\uparrow\!e(T^{vt}Y),ev(\gamma)\rangle=\langle
2PD(h),\downarrow \,(ev(\gamma))\rangle,\end{equation} so by \fullref{euler}
\begin{equation}\label{c-one} \langle
c_1(T^{vt}Y_{\Phi_{d,\omega,\tau}}),ev(\gamma)\rangle=\frac{1}{2}\langle
e(T^{vt}Y)+2PD(h),\downarrow\! ev(\gamma)\rangle.\end{equation}

Meanwhile since we are choosing $\omega$ as a representative of
$c+\frac{2\pi}{\tau}PD(h)$, we have \begin{align}\label{eval}
\langle[(\Omega_{d,\omega,\tau})_{\Phi_{d,\omega,\tau}}],ev(\gamma)\rangle&=2\pi\langle\tau
(c+2\pi PD(h)/\tau)-\pi(2PD(h)),\downarrow\!
ev(\gamma)\rangle\nonumber\\&=2\pi\tau\langle c,\downarrow\!
ev(\gamma)\rangle.\end{align}

Therefore, in light of \eqref{gennov}, $HF(Y,f,h,c,\tau)$ is defined
over the Novikov ring \[ \Nov\left(\ker\langle
e(T^{vt}Y)+2PD(h),\cdot\rangle,\langle c,\cdot\rangle;R\right).\] In
the introduction, allusions were made to ``a spin$^c$ structure
$s_h$ corresponding to $h\in H_1(Y;\mathbb{Z})$;'' we clarify that
slightly here: in \cite[Section 8]{Sal}, Salamon defines a
canonical spin$^c$ structure on the total space of $f\co Y\to S^1$
which has (rank 2, Hermitian) spinor bundle
$S=\underline{\mathbb{C}}\oplus T^{vt}Y$ (the reader who prefers to
think of spin$^c$ structures as given by nonvanishing vector fields
can identify this as the structure specified by a vector field which
is transverse to the fibers of $f$). $s_h$ is then defined by
tensoring $S$ with a line bundle $L$ such that $c_1(L)=PD(h)$ and
extending Clifford multiplication trivially.  Evidently, then, \[
c_1(s_h)=c_1(L\oplus(T^{vt}Y\otimes L))=e(T^{vt}Y)+2PD(h),\] proving
that the Novikov ring we are considering is precisely the ring
$\tilde{\Lambda}_{h,c}$ specified in Section 1.1. (If one instead
uses the convention that the natural coefficient ring for
$CF^{\symp}(\Phi,\mathcal{P})$ is \[ \Nov\left(\frac{\ker\langle
c_1(T^{vt}Y_{\Phi}),ev(\cdot)\rangle}{\ker\langle
c_1(T^{vt}Y_{\Phi}),ev(\cdot)\rangle\cap\langle[\omega_{\Phi}],ev(\cdot)\rangle},\langle[\omega_{\Phi}],ev(\cdot)\rangle;R\right)\]
one evidently instead obtains the smaller ring $\Lambda_{h,c}$
here.)

\subsection{Invariance}  In this subsection, we shall prove \fullref{coeffs}.  We remark first of all that for invariance results
such as \fullref{coeffs} which equate the homologies of two
Floer chain complexes $CF_-$ and $CF_+$ which depend on different
auxiliary data, the usual technique of proof has for some time been
the ``method of continuation,'' wherein one defines a chain map
$CF_-\to CF_+$ by counting finite-energy solutions to some modified
version of the Cauchy--Riemann equations on the cylinder
$\mathbb{R}_t\times S^1$ which has the property that such solutions
are asymptotic to generators of $CF_{\pm}$ as $t\to \pm \infty$.  In
our setting, in which we consider the effect on $HF(Y,f,h,c,\tau)$
of varying the parameter $\tau$, na\"ive attempts to use this method
do not appear to work.  Indeed, the method of continuation would
suggest that, to equate $HF(Y,f,h,c,\tau_-)$ with
$HF(Y,f,h,c,\tau_+)$ for (say) $\tau_-<\tau_+$, we should consider
maps \[ u\co \mathbb{R}\times S^1\to \mathbb{R}\times Y_d(f) \]
which satisfy a perturbed Cauchy--Riemann equation for an almost
complex structure on $\mathbb{R}\times Y_{d}(f)$, which is
compatible with a form $\Omega$ on $\mathbb{R}\times Y_d(f)$ which
agrees with the form induced by $\Omega_{d,\omega,\tau_{\pm}}$ as
the $\mathbb{R}$ parameter approaches $\pm\infty$.  But since
$\Omega_{d,\omega,\tau_+}$ and $\Omega_{d,\omega,\tau_-}$ are not
cohomologous, such a form $\Omega$ could not be closed, and this
would prevent us from obtaining the energy bounds on these maps $u$
which are needed to show the continuation method validly defines a
map between the two chain complexes.

Instead of the continuation method, we make use of recent work of
Lee \cite{L1,L2}) which will enable us to understand in
fairly explicit terms how the chain complexes $CF(Y,f,h,c,\tau)$
vary as $\tau$ increases from $\tau_-$ to $\tau_+$.  It is
interesting to note that A Floer himself used a similar method in
his original paper \cite{Floer} on Lagrangian Floer homology, though
some details needed to justify this approach did not appear until
Lee's work.  Lee  was concerned with torsion invariants rather than
homology in her work, and so did not explicitly prove an invariance
theorem for homology, even though such a theorem readily follows
from her analysis  (Lee also did not consider the effect of smoothly
varying the symplectic form; however this does not complicate the
analysis as long as the ``$H^1$--codirectionality'' hypothesis
discussed below is maintained). Since this result may be useful in
other contexts, we state it in general form here and give an outline
of the algebraic arguments needed in the proof, referring readers to
\cite{L1} and \cite{L2} for the subtle geometric and analytic
arguments necessary to show that the Floer complexes behave as we
claim.  Recall from \cite[Definition 2.3.1]{L1} that a
$2n$-dimensional symplectic manifold $(X,\omega)$ is $w^+$--monotone
provided that every sphere with Chern number strictly between $0$
and $n-1$ has positive symplectic area.

\begin{theorem}[Lee]\label{Lee}  Let $\{\omega_t\}_{t\in [0,1]}$ be a smooth family of symplectic forms on a manifold $X$ such that
the symplectic manifolds $(X,\omega_t)$ are $w^+$--monotone, let
$\phi_t\co X\to X$ be a smooth family of diffeomorphisms such that
$\phi_{t}^{*}\omega_t=\omega_t$, and let $\mathcal{P}\in
\pi_0(\Gamma(Y_{\phi_t}))$ be such that the action functionals
$\mathcal{Y}^t$ for the Floer complexes
$CF^{\symp}(\phi_t,\mathcal{P})$ are $H^1$--codirectional, in the
sense that, where $K\leq \pi_1(\Gamma(Y_{\phi_t}))$ denotes the
kernel of the spectral flow homomorphism we have
$[\mathcal{Y}_t]|_K=f(t)[\mathcal{Y}_0]|_K$ where $f$ is a
nonnegative continuous function; thus where $\Lambda_t$ is the
Novikov ring over which $CF^{\symp}(\phi_t,\mathcal{P})$ is naturally
defined $\Lambda_0$ is a module over $\Lambda_t$, with
$\Lambda_0=\Lambda_t$ when $f(t)\neq 0$.  Then \[
CF^{\symp}(\phi_0,\mathcal{P})\mbox{ is chain homotopy equivalent to
} CF^{\symp}(\phi_1,\mathcal{P})\otimes_{\Lambda_1}\Lambda_0.\]
\end{theorem}

We first explain how \fullref{coeffs} follows from this result.
In \fullref{coeffs}, we consider two cases where we allow one
entry from the standard data $(Y,f,h,c,\tau)$ to vary in a
one-parameter family; namely, we either: \begin{itemize}
\item[(i)] Let $\tau$ vary from $\tau_-$ to $\tau_+$ (say $\tau_t=(1-t)\tau_-+t\tau_+$), fixing
$(Y,f,h,c)$, or
\item[(ii)] Let $c$ vary from $\tilde{c}$ to $\alpha c_1(s_h)$ (say
$c_t=t\alpha c_1(s_h)+(1-t)\tilde{c}$; here $\alpha\in
\mathbb{R}\setminus\{0\}$), fixing $(Y,f,h,\tau)$.
\end{itemize}

Recall that in the statement of \fullref{coeffs} we have assumed
that either $d\geq g-1$ or $d<(g+1)/2$; since $\pi_2(S^d\Sigma)$ is
an infinite cyclic group generated by a sphere on which the
symplectic form is positive and $c_1$ evaluates as $d-g+1$
\cite{Mac} this is the assumption needed to ensure that $S^d\Sigma$
is $w^+$--monotone and so for \fullref{Lee} to apply. For the
range $(g+1)/2\leq d<g-1$ it seems likely that the methods of
\cite{LT} could be used to prove similar results to the ones we cite
here, but this does not appear to be at all straightforward.

In either case, let $\Omega_t$ be the symplectic form obtained by
Salamon's construction using the data at time $t$.  The relevant
action one--form, defined on the space of sections $\gamma$ of
$Y_d(f)$ representing the homotopy class $\mathcal{P}_h$, is then \[
\mathcal{Y}^{t}_{\gamma}(\xi)=-\int_{0}^{1}\Omega_t
(\dot{\gamma}(s),\xi(s))\,ds;\] the action of this $1$--form on
loops in $\mathcal{P}_h$ is given by Equation \ref{eval}.  Meanwhile
the action of $c_1(T^{vt}Y_d(f))$ (which plays the role here of
Lee's spectral flow homomorphism $\psi\co \pi_1(\mathcal{P}_h)\to
\mathbb{Z}$) is given by $\langle c_1(s_h),\downarrow\!
ev(\cdot)\rangle.$  So in case (i) above, the classes of the forms
$\mathcal{Y}^t$ in $H^1(\mathcal{P}_h;\mathbb{R})$ satisfy
\[ [\mathcal{Y}^t]|_{\ker \psi}=\frac{\tau_t}{\tau_0}
[\mathcal{Y}^0]|_{\ker \psi}\] while in case (ii)  \[
[\mathcal{Y}^t]|_{\ker \psi}=(1-t) [\mathcal{Y}^1]|_{\ker \psi}.\]
Thus in both cases, our path of $1$--forms $\mathcal{Y}^t$ on the
infinite-dimensional space $\mathcal{P}_h$ is $H^1$--codirectional.
This fact makes \fullref{Lee} relevant to our situation.

We now briefly outline the facts from \cite{L1} which enter into the proof of \fullref{Lee}.
Lee's work implies the existence of a ``regular homotopy of
Floer systems'' (RHFS) between the (partially-defined) flows on the
space $\mathcal{P}\subset Y_{\phi_t}$ which underlie complexes
$CF^{\symp}(\phi_0,\mathcal{P})$ and $CF(\phi_1,\mathcal{P})$.
Namely, there is a path $(J_{t},H_{t})_{t\in [0,1]}$ of
$\omega_t$--compatible almost complex structures and Hamiltonian
perturbations such that as $t$ varies the chain complexes
$CF^{\symp}(\phi_t,\mathcal{P})$ change only at certain values of $t$
corresponding to ``handleslides'' and (just finitely many)
``death-births,'' all in the complement of a set $S_{reg}$ of second
category in $[0,1]$; the Floer complex corresponding to $(J_t,H_t)$
is well-defined for each $t\in S_{reg}$.

If $[t_0,t_1]\subset [0,1]$ is an interval in which the complex
changes only by handleslides (that is, Floer flow lines between
generators of equal index; there may be infinitely many of these,
but only finitely many with energy below any given bound) and
$t_0,t_1\in S_{reg}$, the generators for the chain complex remain
unchanged throughout the interval, while the differentials
$\partial_0,\partial_1$ of the chain complexes at times $t_0,t_1$
are related by
\begin{align}\label{hs} \langle \partial_1 x,y[p_2(A)]\rangle=&\langle\partial_0 x,y[p_2(A)]
\rangle+\\&\sum_{z,B}\sum_{s\in
\mathcal{M}^{[t_0,t_1]}_{hs}(\mathbf{x},\mathbf{z};B)}\epsilon(s)\langle\partial_0
z,y[p_2(B-A)]\rangle+\nonumber
\\&\sum_{z,C}\sum_{s\in
\mathcal{M}^{[t_0,t_1]}_{hs}(\mathbf{z},\mathbf{y};C)}-\epsilon(s)\langle\partial_1
x,z[p_2(C-A)]\rangle \nonumber \end{align} where in general
$\mathcal{M}^{[t_0,t_1]}_{hs}(\mathbf{u},\mathbf{v};B)$ denotes the
moduli space of handleslides between the generators $\mathbf{u}$ and
$\mathbf{v}$ having Novikov ring weight $B$, and for each
handleslide $s$ $\epsilon(s)=\pm 1$ is a sign determined by coherent
orientations.

If we then set \[ Tx=x+\sum_{z,B}\sum_{s\in
\mathcal{M}^{I}_{hs}(\mathbf{x},\mathbf{z};B)}\epsilon(s)
z[p_2(B)],\] then the matrix element $\langle
\partial_0Tx,y[p_2(A)]\rangle$ is the sum of the first two terms on
the right hand side of \eqref{hs}, while $\langle T\partial_1
x,y[p_2(A)]\rangle$ is the difference of the term on the left hand
side and the last term on the right. Hence $T$ defines a chain map
$CF^{\symp}(\phi_1,\mathcal{P})\otimes \Lambda_0\to
CF^{\symp}(\phi_0,\mathcal{P})$ (note that this is well-defined over
$\Lambda_0$ as a result of the finiteness condition on handleslides
mentioned earlier).  But $T$ is an invertible map (if we write
$T=I+U$ where $U$ is the identity, $\sum_{n=0}^{\infty}(-U)^n$ will
be a well-defined endomorphism over the Novikov ring and will
obviously be inverse to $T$), so it actually defines an isomorphism
of chain complexes.

This reduces the invariance problem to showing that the Floer
homology is unchanged on intervals containing death-births, of which
there are only finitely many.  Let $I\subset [0,1]$ be an open
interval containing a single death or birth, say at $\bar{t}$ (we'll
assume it's a birth; the death case may be obtained by reversing
various arrows and inequality signs in the discussion below).  For
$t<t'<\bar{t}$ such that $t,t'\in S_{reg}$ we have isomorphisms of
chain complexes
\[ T_{t,t'}\co CF^{\symp}(\phi_t,\mathcal{P})\otimes \Lambda_0\to CF^{\symp}(\phi_{t'},\mathcal{P})\otimes \Lambda_0 \]
as above; these form a directed system indexed by a dense subset of
$\{t\in I|t<t'\}$, so we may let $(CF^{-}_{\bar{t}},\partial^-)$ be
the direct limit of the $CF^{\symp}(\phi_t,\mathcal{P})$ under this
directed system. Likewise let $(CF^{+}_{\bar{t}},\partial^+)$ be the
inverse limit of the directed system indexed by the subset $\{t\in
I\cap S_{reg}|t>\bar{t}\}$ of $I$, given by the isomorphisms
$T_{t,t'}$. We pass to these limits in order to allow ourselves to
ignore handleslides in the following discussion, since $\bar{t}$
might not be contained in any open interval over which there are no
handleslides.

Lee's axioms for an RHFS in \cite[Section 4]{L1} (see also
Hutchings \cite{Htors} for a more explicit description in the finite
dimensional context) imply the following description of the
relationship of $(CF^{+}_{\bar{t}},\partial^+)$ to
$(CF^{-}_{\bar{t}},\partial^-)$.  We have \[
CF^{+}_{\bar{t}}=CF^{-}_{\bar{t}}\oplus \Lambda_{0}\langle
x^+,x^-\rangle \] for some two generators $x^{\pm}$ which differ in
relative grading by $1$; these new generators are ``born'' from a
degenerate critical point $x^0$ of the action form
$\mathcal{Y}^{\bar{t}}$ which appears at $t=\bar{t}$.

Using Lee's axioms RHFS2 and RHFS2c, one can deduce the following
description for the relationship between the Floer boundary
operators $\partial^-$ and $\partial^+$.  There are maps $v\co
\Lambda_0\to CF^{-}_{\bar{t}}$, $\mu\co
CF^{-}_{\bar{t}}\to\Lambda_0$ and an invertible element $\alpha\in
\Lambda_0$ such that, with respect to an ordered basis for
$CF^{+}_{\bar{t}}$ consisting of an ordered basis for
$CF^{-}_{\bar{t}}$ followed by $(x^+,x^-)$, the differential
$\partial^+$ may be written in block form as \[
\partial^+=\left(\begin{array}{ccc}\partial^-+\alpha^{-1}v\circ\mu &
v & 0\\0 & 0 & 0\\ \mu & \alpha & 0\end{array}\right).\]  Note that
the fact that $(\partial^+)^2=0$ then implies that $\mu\circ
v=\mu\circ\partial^-=\partial^-\circ v=0$.  Define a map $i\co
CF^{-}_{\bar{t}}\to CF^{+}_{\bar{t}}$ by the block matrix
\[ i=\left(\begin{array}{c}Id\\ -\alpha^{-1}\mu \\
0\end{array}\right)
\] and a map $p\co CF^{+}_{\bar{t}}\to CF^{-}_{\bar{t}}$ by \[
p=\left(\begin{array}{ccc} Id & 0 &
-\alpha^{-1}v\end{array}\right).\]  The fact that $\mu\circ
\partial^-=0$ implies that $i$ is a chain map, while the fact that
$\partial^-\circ v=0$ implies that $p$ is a chain map.  Now
obviously $p\circ i$ is the identity, while defining $K\co
CF^{+}_{\bar{t}}\to CF^{+}_{\bar{t}}$ by
\[K=\left(\begin{array}{ccc}0&0&0\\0&0&\alpha^{-1}\\0&0&0\end{array}\right),\]
one easily computes (using the fact that $\mu\circ v=0$) that \[
\partial^+K+K\partial^+=1-i\circ
p=\left(\begin{array}{ccc}0&0&\alpha^{-1}v\\\alpha^{-1}\mu&1&0\\0&0&1\end{array}\right).\]
$i$ and $p$ thus put $CF^{-}_{\bar{t}}$ and $CF^{+}_{t}$ into chain
homotopy equivalence.

\fullref{Lee} is immediate from this, for the interval $[0,1]$
contains just finitely many values of $t$ (say $t_1,\ldots,t_n$) at
which death-births occur; write $t_0=0,t_{n+1}=1$). Our treatment of
handleslides shows that the chain complexes
$CF^{\symp}(\phi_t,\mathcal{P})\otimes \Lambda_0$ are mutually
isomorphic for all $t\in S_{reg}\cap [t_i,t_{i+1}]$, and our
treatment of births shows that the direct limit
$\lim_{\rightarrow}^{t<t_i}CF^{\symp}(\phi_t,\mathcal{P})\otimes\Lambda_0$
under these isomorphisms is chain homotopy equivalent to
$\lim_{\leftarrow}^{t>t_i}CF^{\symp}(\phi_t,\mathcal{P})\otimes
\Lambda_0$, implying that the various
$CF^{\symp}(\phi_t,\mathcal{P})\otimes\Lambda_0$ for all regular $t$
in the entire interval are chain homotopy equivalent. This completes
our account of the proof of \fullref{Lee}; as explained earlier,
the assumption that $d\notin [(g+1)/2,g-1)$ makes \fullref{coeffs} a special case of this more general result.

\section{Further algebraic properties}\label{sec.4}

With $HF(Y,f,h,c,\tau;A)$ defined as the Floer homology of the
symplectomorphism $\Phi_{d,\omega,\tau}$, the algebraic properties
alluded to in Section 1.2 now follow quickly from standard
properties of Floer homology.  First, as always in Floer theory, we
have a absolute $\mathbb{Z}/2$ grading provided here by the
Lefschetz index: a generator of $CF(Y,f,h,c,\tau;A)$ is a fixed
point $p$ of (possibly a Hamiltonian perturbation of)
$\Phi_{d,\omega,\tau}$, and the absolute grading of $p$ is just \[
\sign \det(Id-(d\Phi_{d,\omega,\tau})_p).\]  The relative grading
also follows from the standard setup and equation \eqref{c-one};
quite generally the ambiguity in the relative grading of a generator
of the Floer homology of $HF^{\symp}(\psi;\mathcal{P})$ of a
symplectomorphism $\psi$ in the fixed point class $\mathcal{P}$ is
given by \[ 2\gcd_{\gamma\in H_1(\mathcal{P})}\langle
c_1(T^{vt}Y_{\psi}),ev(\gamma)\rangle,\] and in our case with
$\psi=\Phi_{d,\omega,\tau}$ and $\mathcal{P}=\mathcal{P}_h$, this is
precisely the divisibility of $c_1(s_h)=e(T^{vt}Y_{\phi})+2PD(h)$,
as stated in \fullref{alg}.

Poincar\'e duality for $HF$ is still another simple consequence of
the setup: replacing the tuple $\frak{o}=(Y,f,h,c,\tau)$ with
$\bar{\frak{o}}=(-Y,\bar{f},-h,c,\tau)$ of course preserves the
orientations of the fibers, while we have $PD_{-Y}(-h)=PD_Y(h)$, so
the same fiberwise symplectic form $\omega$ representing
$c+\frac{2\pi}{\tau}PD_{\pm Y}(\pm h)$ can be used for both
$\frak{o}$ and $\bar{\frak{o}}$ and the same fiberwise complex
structures $J_t$ can be used on $Y\to S^1$ and $-Y\to S^1$.  Using
these same auxiliary data, the horizontal vector field whose flow we
use to define the monodromy $\Phi^{\bar{\frak{o}}}_{d,\omega,\tau}$
in Salamon's construction for $\bar{\frak{o}}$ will then be
precisely the opposite of the vector field which generates the
monodromy $\Phi^{\frak{o}}_{d,\omega,\tau}$.  Hence
$\Phi^{\bar{\frak{o}}}_{d,\omega,\tau}=(\Phi^{\frak{o}}_{d,\omega,\tau})^{-1}$.
Now it is quite generally the case that, for Floer chain complexes
$CF^{\symp}(\psi,\mathcal{P})$ of symplectomorphisms $\psi$,
$CF^{\symp}(\psi^{-1},\bar{\mathcal{P}})$ is naturally the dual
complex to $CF^{\symp}(\psi,\mathcal{P})$ under an appropriate
identification $\mathcal{P}\leftrightarrow\bar{\mathcal{P}}$ of
fixed point classes: a fixed point of $\psi$ is of course also a
fixed point of $\psi^{-1}$ and,  tautologically, two such are
Nielsen--equivalent for $\psi$ iff they are for $\psi^{-1}$, so as
groups
$CF^{\symp}(\psi^{-1},\bar{\mathcal{P}})=CF^{\symp}(\psi,\mathcal{P})$.
The differentials are related by the observation that since the
mapping torus fibration $Y_{\psi^{-1}}\to S^1$ is (up to equivalence
of fibrations) the conjugate of $Y_{\psi}\to S^1$, a cylinder in
$\mathbb{R}\times Y_{\psi}$ which serves as a flowline from the
generator $a$ to the generator $b$ in $CF^{\symp}(\psi)$ is the same
thing as a cylinder in $\mathbb{R}\times Y_{\psi^{-1}}$ which serves
as a flowline from $b$ to $a$ in $CF^{\symp}(\psi^{-1})$.  In our
context, the fixed point class $\mathcal{P}_h$ for
$\Phi^{\frak{o}}_{d,\omega,\tau}$ corresponds tautologically to the
fixed point class $\mathcal{P}_{-h}$ for
$\Phi^{\bar{\frak{o}}}_{d,\omega,\tau}$, and so $CF(\frak{o})$ and
$CF(\bar{\frak{o}})$ are the same as groups and have dual
differentials, which proves \fullref{pairing} (the pairing
between $CF(\frak{o})$ and $CF(\bar{\frak{o}})$ promised therein is
of course  obtained by, for generators $\mathbf{x}$ and $\mathbf{y}$
of the identical groups, setting $\langle
\mathbf{x},\mathbf{y}\rangle$ to be $1$ if $\mathbf{x}=\mathbf{y}$
and $0$ otherwise and then extending linearly).

The only remaining item from \fullref{alg} is the structure
of $HF(Y,h,c,\tau;A)$ as a module over\[
\mathbb{A}(Y)=\mathbb{Z}[U]\otimes
\Lambda^{*}(H_1(Y;\mathbb{Z})/\torsion).\] We obtain this by
considering the quantum cap product structure in Floer theory, which
we describe here in the case that virtual cycle methods are not
needed.  Quite generally, the boundary operator $\partial\co CF\to
CF$ in a Floer theory with configuration space $\mathcal{C}$ counts
paths $\gamma\co [-\infty,\infty]\to \mathcal{C}$ with prescribed
endpoints which are (formally) gradient lines for some Morse
function on $\mathcal{C}$. Letting
$\mathcal{M}(\mathbf{x},\mathbf{y};A)$ denote the moduli space of
flowlines from a generator $\mathbf{x}$ of $CF$ to a generator
$\mathbf{y}$ having relative homotopy class $A$, if $k$ denotes the
common index of these flowlines then evaluation of the flowline at
time zero determines a $k$--dimensional chain
$ev_*\mathcal{M}(\mathbf{x},\mathbf{y};A)$ in $\mathcal{C}$. So if
$a\in C^{k}(\mathcal{C};\mathbb{Z})$ we get a degree--$(-k)$ map
\[ a\cdot\co CF\to CF\] by setting \[ a\cdot \mathbf{x}=\sum_{\mathbf{y},A}\langle
a,ev_*\mathcal{M}(\mathbf{x},\mathbf{y};A)\rangle[p_2(A)]\mathbf{y};\]
considering the boundary components of the
$\mathcal{M}_k(\mathbf{x},\mathbf{y};A)$ reveals that $a\cdot$ is a
chain map.  Further, the map that it induces on the Floer homology
$HF$ depends only on the cohomology class of $a$, and the resulting
map $H^{*}(\mathcal{C};\mathbb{Z})\times HF\to HF$ makes $HF$ into a
module over $H^*(\mathcal{C};\mathbb{Z})$.  For more details on this
see Viterbo \cite{Vit} and Liu--Tian \cite{LTprod}, in the latter of which it is shown
that the quantum cap product can be made compatible with virtual
cycle machinery.

In our context, the configuration space $\mathcal{C}$ is a homotopy
class $\mathcal{P}_h$ of sections of the degree--$d$ relative
symmetric product $Y_d(f)\to S^1$ of a fibered 3--manifold $f\co
Y\to S^1$.  Thus the quantum cap product construction makes
$HF(Y,f,h,c,\tau;A)$ into a module over
$H^{*}(\mathcal{P}_h;\mathbb{Z})$.  To obtain the asserted module
structure over $\mathbb{A}(Y)$ (and so exhibit still another
parallel with the monopole and Heegaard Floer theories), we thus
just need to exhibit a natural graded ring homomorphism \[
\mathbb{Z}[U]\otimes \Lambda^{*}(H_1(Y;\mathbb{Z})/\torsion)\to
H^*(\mathcal{P}_h;\mathbb{Z}),\] where $U$ has degree 2. In this
direction, note that if $g\co K\to \mathcal{P}_h$ is a chain of
dimension $k$, we get a dimension--$(k+1)$ chain $e(g)\co S^1\times
K\to Y_d(f)$ by setting $e(g)(t,k)=(g(k))(t)$; $e$ evidently defines
a degree--$1$ chain map
\[ e\co C_*(\mathcal{P}_h;\mathbb{Z})\to
C_{*+1}(Y_d(f);\mathbb{Z});\] dualizing this and passing to
cohomology yields a homomorphism \[ e^*\co H^*(Y_d(f);\mathbb{Z})\to
H^{*-1}(\mathcal{P}_h;\mathbb{Z}).\]  Recalling the map $\uparrow\co
H^*(Y;\mathbb{Z})\to H^*(Y_d(f);\mathbb{Z})$, we now set \[
U=e^*\left(\uparrow(PD[pt])\right)\in
H^2(\mathcal{P}_h;\mathbb{Z})\] and use the homomorphism
\begin{align*} H_1(Y;\mathbb{Z})&\to H^1(\mathcal{P}_h;\mathbb{Z})\\
\gamma&\mapsto e^*\uparrow PD(\gamma).\end{align*}  Note that the
image of any torsion elements of $H_1(Y;\mathbb{Z})$ will be
trivial, and so the above map factors through a map \[
H_1(Y;\mathbb{Z})/\torsion\to H^1(\mathcal{P}_h;\mathbb{Z}).\]

  Now that we have
chosen an element $U\in H^2(\mathcal{P}_h;\mathbb{Z})$ and a
homomorphism $$H_1(Y;\mathbb{Z})/\torsion\to
H^1(\mathcal{P}_h;\mathbb{Z}),$$ a unique ring homomorphism \[
\mathbb{A}(Y)=\mathbb{Z}[U]\otimes
\Lambda^{*}(H_1(Y;\mathbb{Z})/\torsion)\to
H^*(\mathcal{P}_h;\mathbb{Z})\] is forced on us by the graded ring
structure of $H^*(\mathcal{P}_h;\mathbb{Z})$.  This completes the
proof of the existence of the module structure over the promised
ring.

More geometrically, the map on $HF(Y,f,h,c,\tau;A)$ induced by the
element \[ U^r\otimes \gamma_1\wedge\cdots\wedge\gamma_k\in
\mathbb{Z}[U]\otimes\Lambda^*(H_1(Y;\mathbb{Z})/\torsion)\] counts
holomorphic sections \begin{align*} \mathbb{R}\times S^1&\to
\mathbb{R}\times Y_d(f)\\(s,t)&\mapsto u(s,t)\end{align*} with the
property that the cycle obtained by taking the union of the points
in the divisors $u(0,t)$ for the various $t\in S^1$ contains a
generic set of $r$ points and meets generic representatives of the
classes $\gamma_1,\ldots,\gamma_k$.

Finally, we explain the construction of the local coefficient
systems $\Gamma_{\theta}$ for closed forms $\theta\in \Omega^2(Y)$.
First, we choose once and for all a de Rham representative $\delta$
of the class $PD(\mathcal{D})\in H^2(Y\times_{S^1}
Y_d(f);\mathbb{R})$. Now from our initial data set $(Y,f,h,c,\tau)$
we have constructed a closed form $\omega\in \Omega^2(Y)$
restricting to each fiber of $f\co Y\to S^1$ as a volume form;
pulling back this form by $\pi_1\co Y\times_{S^1} Y_d(f)\to Y$ gives
a closed $2$--form on $Y\times_{S^1} Y_d(f)$ which restricts as a
volume form to each fiber of $\pi_2\co Y\times_{S^1} Y_d(f)\to
Y_d(f)$; this gives rise via integration down the fibers of $\pi_2$
to a form--level extension
\[ (\pi_2)_{!}\co \Omega^*(Y\times_{S^1} Y_d(f))\to
\Omega^{*-2}(Y_d(f))\] of the Gysin map, and so to a form--level
extension \begin{align}\label{formlevel} \uparrow_{\delta}\co
\Omega^*(Y)&\to \Omega^*(Y_d(f)) \\
\theta&\mapsto (\pi_2)_{!}(\delta\wedge\pi_{1}^{*}\theta)\nonumber
\end{align} of $\uparrow$ which is a cochain map and so takes closed forms to
closed forms.  This yields a local system $\Gamma_{\theta}$ on the
configuration space $\mathcal{P}_h$ as follows: to each
$\mathbf{x}\in \mathcal{P}_h$ take
$(\Gamma_{\theta})_{\mathbf{x}}=\mathbb{R}$ as the fiber over
$\mathbf{x}$; since an element of $\mathcal{P}_h$ gives rise via
evaluation to a loop in $Y_d(f)$, a path $\gamma$ from $\mathbf{x}$
to $\mathbf{y}$ in $\mathcal{P}_h$ gives rise  via evaluation to a
cylinder $C_{\gamma}\subset Y_d(f)$ with boundary components
prescribed by $\mathbf{x}$ and $\mathbf{y}$, and we define the
isomorphism \[ \phi_{\theta}([\gamma])\co
(\Gamma_{\theta})_{\mathbf{x}}\to (\Gamma_{\theta})_{\mathbf{y}}\]
to be multiplication by \[
\exp\left(\int_{C_{\gamma}}\uparrow_{\delta}\theta\right);\] that
$\uparrow_{\delta}\theta$ is closed of course ensures that this
depends only on the homotopy class of $\gamma$ relative to its
boundary. The resulting twisted Floer groups are then obtained by
letting the twisted Floer chain complex be the direct sum over the
fixed points $x$ of $\Phi_{d,\omega,\tau}$ of the groups
$\Gamma_{\theta}(\mathbf{x})\otimes \tilde{\Lambda}_{h,c}$ (where
$\mathbf{x}\in\mathcal{P}_h$ is the section of $Y$ corresponding to
$x$) and weighting each term in the standard Floer differential
(which corresponds to a path $\gamma$ in $\mathcal{P}_h$) by
$\phi_{\theta}([\gamma])$; see \cite[Section 2.7]{KM} for the
analogous construction in Morse theory.  If $\zeta\in \Omega_1(Y)$
satisfies $d\zeta=\theta_2-\theta_1$, then if $\gamma$ is a path in
$\mathcal{P}_h$ from $\mathbf{x}$ to $\mathbf{y}$, we have \[
\frac{\phi_{\theta_2}([\gamma])}{\phi_{\theta_1}([\gamma])}=\exp\left(\int_{C_{\gamma}}\uparrow_{\delta}
d\zeta\right)=\exp\left(\int_{C_{\gamma}}d\uparrow_{\delta}
\zeta\right)=\frac{\exp(\int_{\mathbf{y}}\uparrow_{\delta}\zeta)}{\exp(\int_{\mathbf{x}}\uparrow_{\delta}\zeta)},\]
\begin{align*} \Gamma_{\theta_1}(\mathbf{x})\tag*{\hbox{and so}} &\to
\Gamma_{\theta_2}(\mathbf{x})\\ z&\mapsto
z\exp\left(\int_{\mathbf{x}}\uparrow_{\delta}\zeta\right)
\end{align*} defines an isomorphism of local systems
$\Gamma_{\theta_1}\cong\Gamma_{\theta_2}$, as claimed in the
introduction.

The construction of the $\Gamma_{\theta}$ depends on a choice of de
Rham representatives of the classes $PD(\Delta)\in
H^2(Y\times_{S^1}Y_d(f);\mathbb{R})$ and
$c+\frac{2\pi}{\tau}PD(h)\in H^2(Y;\mathbb{R})$; it is easy to see
that different choices of these representatives give rise to local
systems which are isomorphic, with the isomorphism depending on
cohomologies between the different choices.

\section{Maps from cobordisms}\label{sec.5}

The basic ingredient in the construction of the maps on $HF$
obtained from ``fibered cobordisms'' between objects
$(Y,f,h,c,\tau)$ in the category $\mathsf{FCOB}$ is the fact that,
given a Lefschetz fibration $f\co X\to B$ with closed fibers, a
fiberwise positive class $c\in H^2(X;\mathbb{R})$, and a class $h\in
H_2(X,\partial X;\mathbb{Z})$ with $h\cap [\fiber]=d>0$, one can
construct a \emph{relative Hilbert scheme} $F\co X_d(f)\to B$ such
that for regular values $t$ of $f$ one has
$F^{-1}(t)=Sym^df^{-1}(t)$, and which carries a symplectic form
$\Omega$ such that for each component $S$ of $\partial B$,
$(F^{-1}(S),\Omega|_{F^{-1}(S)})$ is isomorphic as a locally
Hamiltonian fibration to the mapping torus
$Y_{\Phi_{d,\omega,\tau,S}}$, $\Phi_{d,\omega,\tau,S}$ being the
Salamon monodromy map defining the Floer group
$HF(Y,f|_{f^{-1}(S)},\partial_S h,c|_{f^{-1}(S)},\tau)$.  In the
case that $\partial B$ has two components and $X$ is a morphism in
$\mathsf{FCOB}$ from $\frak{o}_-$ to $\frak{o}_+$, $(X_d(f),\Omega)$
then defines a symplectic cobordism from the mapping torus used to
define $HF(\frak{o}_-)$ to that used to define $HF(\frak{o}_+)$, and
the map between the Floer groups is obtained by adding half-infinite
cylindrical ends to this symplectic cobordism and then counting
pseudoholomorphic cylinders with prescribed asymptotics.

The construction of the relative Hilbert scheme $X_d(f)$, including
the proof of the crucial fact that it is smooth in spite of the
presence of singular fibers in the Lefschetz fibration $f\co X\to
B$, is carried out in detail in Smith \cite{Smith}.  The existence of a
natural deformation equivalence class of symplectic structures on
$X_d(f)$ is also proven in \cite{Smith}, but in order to obtain the
proper behavior of the symplectic form on $\partial X_d(f)$ we shall
need somewhat more refined results, which we now set about proving.

Recall that in Section 2, given an object $(Y,f,h,c,\tau)$ in
$\mathsf{FCOB}$ and a closed form $\omega\in\Omega^2(Y)$
representing $c+\frac{2\pi}{\tau}PD(h)$, when $d>0$ we have obtained
a symplectic form $\Omega_{d,\omega,\tau}$ on $S^d\Sigma$ and a map
$\Phi_{d,\omega,\tau}\co S^d\Sigma\to S^d\Sigma$ which preserves
$\Omega_{d,\omega,\tau}$. (When $d\leq 0$ we have set $HF=0$ except
when $h=0$, in which case $HF$ is the coefficient ring $A$; in all
such cases we hereby define the cobordism maps of this section to be
the identity, and so restrict attention to the case $d>0$
hereinafter.) The mapping torus $Y_{\Phi_{d,\omega,\tau}}$ then
carries a closed fiberwise symplectic form
$(\Omega_{d,\omega,\tau})_{\Phi_{d,\omega,\tau}}$; note that
$Y_{\Phi_{d,\omega,\tau}}\cong Y_d(f)$ as fibrations over $S^1$.
Recall also the definition of a starred  surface $B$ in the
introduction, which in particular specifies a set
$\{\beta_1,\ldots,\beta_{g(B)}\}$ of disjoint curves on $B$.

\begin{lemma}\label{relfib} Fix $\tau>2\pi d$, and let $f\co X\to B$ be a fibration by closed surfaces
of genus $g\geq 2$ over a compact connected starred surface $B$
whose boundary decomposes into connected components as $\partial
B=\partial_1 B\cup\cdots\cup\partial_n B.$ Write
$Y_i=f^{-1}(\partial_i B)$, and let $c_i\in H^2(f^{-1}(\partial_i
B);\mathbb{R})$, $h_i\in H^1(f^{-1}(\partial_i B);\mathbb{Z})$, and
$b_j\in H^2(f^{-1}(\beta_j);\mathbb{R})$ be such that there exist
$\tilde{c}\in H^2(X;\mathbb{R})$, $PD(\tilde{h})\in
H^2(X;\mathbb{Z})$ both evaluating positively on the fibers (say
$\langle PD(\tilde{h}), [\fiber]\rangle=d$) with
$c_i=\tilde{c}|_{Y_i}$, $PD(h_i)=PD(\tilde{h})|_{Y_i}$, and
$b_j=(\tilde{c}+2\pi PD(\tilde{h})/\tau)|_{f^{-1}(\beta_j)}$. Let
$\omega_i\in \Omega^2(Y_i)$ be closed fiberwise symplectic forms
representing $c_i+2\pi PD(h_i)/\tau=(c+2\pi PD(h)/\tau)|_{Y_i}\in
H^2(Y_i;\mathbb{R})$. Then there is a fiberwise symplectic form
$\tilde{\Omega}=\Omega_{\vec{c},\vec{h},\vec{b}}\in
\Omega^2(X_d(f))$, determined canonically up to isotopy by the
$c_i,h_i$, and $b_j$ and independent of the choices of $\tilde{c}$
and $\tilde{h}$, such that
\[
[\tilde{\Omega}|_{(Y_i)_d(f)}]=[(\Omega_{d,\omega_i,\tau})_{\Phi_{d,\omega_i,\tau}}]\]
for each $i$.

\end{lemma}

\begin{proof} First note that if $g\co Z\to S^1$ is a fibration and
if $\Omega_-,\Omega_+$ are cohomologous closed fiberwise symplectic
forms on $Z$ (say $\Omega_+=\Omega_-+d\alpha$) which tame a common
fiberwise almost complex structure, then on $[-1,1]_t\times Z$, if
we set $\tilde{\Omega}=\Omega_-+d(\rho(t)\alpha)$ where $\rho\co
[-1,1]\to [0,1]$ is a smooth function which is identically zero near
$-1$ and identically $1$ near $1$, then $\tilde{\Omega}$ is a closed
fiberwise symplectic form on $Z\times [-1,1]\to S^1\times [-1,1]$
which, for small $\epsilon$, restricts to $Z\times[-1,-1+\epsilon]$
as the pullback of $\Omega_-$ and to $Z\times[1-\epsilon,1]$ as the
pullback of $\Omega_+$.  Using this device, in the context of the
lemma it becomes straightforward to construct a closed, fiberwise
K\"ahler form $\tilde{\Omega}$ over all of $B$ satisfying the
desired properties as soon as we have constructed it over the
surface obtained by cutting $B$ along its $\beta$-curves in such a
way that $\tilde{\Omega}$ restricts over $\beta_j$ to
$\Omega_{d,\omega_j,\tau}$, $\omega_j$ being a fiberwise symplectic
representative of $b_j\in H^2(f^{-1}(\beta_j);\mathbb{R})$. This
observation reduces us to the case that $B$ has genus zero.

Having made this reduction, consider first the case that $B=S^2$,
which is simpler in that here there are no boundary conditions for
the form to satisfy; the only datum given to us is the fibration map
$f\co X\to B$, which is a bundle whose structure group is the
identity component $\Diff_0(\Sigma)$ of the diffeomorphism group of
$\Sigma$. Now it follows from Moser stability \cite{Moser} that
$\Diff_0(\Sigma)$ retracts to $\Symp(\Sigma,\omega)$, while since
$g\geq 2$ the Earle--Eells theorem \cite{EE} states that
$\Diff_0(\Sigma)$ is contractible; hence the structure group of the
fibration reduces to $\Symp(\Sigma,\omega)$, which is contractible.
 The fibration thus admits a symplectic trivialization which is unique up to isotopy, giving us a canonical symplectic identification of $X$ with
$S^2\times\Sigma$ with $f$ being the projection to the first factor,
and this induces an identification $X_d(f)$ with $S^2\times
S^d\Sigma$. So we take for $\tilde{\Omega}$ the pullback of the
Salamon form $\Omega_{d,\omega,\tau}$ on any of the individual
fibers.  Having dispensed with this trivial case, assume that $B$
has at least one boundary component.  Let $\tilde{\omega}$ be a
fiberwise symplectic form on the total space of $f\co X\to B$ which
represents the class $\tilde{c}+2\pi PD(\tilde{h})/\tau\in
H^2(X;\mathbb{R})$ for an arbitrary choice of $\tilde{c},\tilde{h}$
as in the statement of the lemma (such $\tilde{\omega}$ exists by
the Thurston trick). After trivializing $f$ over a tubular
neighborhood $S\subset B$ of the contractible ``star'' consisting of
the arcs from the interior base point $b$ to the boundary basepoints
that are part of the data of the starred surface $B$ and modifying
$\tilde{\omega}$ by an exact form, we can assume that
$\tilde{\omega}|_{f^{-1}(S)}$ is just the pullback to
$f^{-1}(S)=S\times \Sigma$ of a volume form $\omega$ on $\Sigma$.
Let $U^{(1)}$ be the union of $S$ with tubular neighborhoods of the
boundary components (so $U^{(1)}$ is a regular neighborhood of a
$1$-skeleton for $B$ since we've reduced to $g(B)=0$).

Using Salamon's construction we can obtain a closed fiberwise
K\"ahler form $\tilde{\Omega}^{(1)}$ over $U^{(1)}$ as follows.
Start with the trivial extension of Salamon's form
$\Omega=\Omega_{d,\omega,\tau}$ on $S^df^{-1}(0)$ to $S\times
S^df^{-1}(0)$, and then, over the strips $(-\epsilon,\epsilon)\times
\partial_i B\subset U^{(1)}$ corresponding to the boundary components attach strips $(-\epsilon,\epsilon)\times (Y_i)_d(f)$ to
$D^2\times S^df^{-1}(0)$. To describe $\tilde{\Omega}^{(1)}$ on
these latter strips, let $\Phi_i\co S^d f^{-1}(0)\to S^d f^{-1}(0)$
be the Salamon monodromy map associated to the mapping torus of the
monodromy of $\tilde{\omega}$ around these loops; the form on the
strip $(-\epsilon,\epsilon)\times (Y_i)_d(f)$ can then just be taken
to be the pullback of the form $(\Omega)_{\Phi_i}$ on $(Y_i)_d(f)$.

Now $U^{(1)}$ has $n+1$ boundary components $C,\partial_1 B,\ldots,
\partial_n B$ with $[C]=\sum [\partial_i B]$ in homology, and $B$ is
obtained by attaching a disc $D'$ to $C$.  Since $C$ is contractible
in $B$, the $\tilde{\omega}$--monodromy around it is Hamiltonian,
from which it follows as in \fullref{Ham} that the
$\tilde{\Omega}$--monodromy around $C$ (which is just the
composition of the monodromies around the $\partial_i B$) is
Hamiltonian as well. Write $F\co X_d(f)\to B$ for the map defining
the relative symmetric product. Then picking $p\in C=\partial D$,
$\tilde{\Omega}|_{F^{-1}(C)}$ is cohomologous to the pullback of
$\Omega|_{F^{-1}(p)}$ to $S^1\times F^{-1}(p)$ as a result of the
fact that they have Hamiltonian--equivalent monodromies, and so just
as in the first paragraph of this proof we can interpolate between
these two forms and so glue $\tilde{\Omega}|_{F^{-1}(C)}$ to the
pullback of $\tilde{\Omega}|_{F^{-1}(p)}$ to $D'\times F^{-1}(p')$,
thus producing the desired form $\tilde{\Omega}$ on all of
$F^{-1}(B)$. (The conclusion about $[\tilde{\Omega}|_{(Y_i)_d(f)}]$
follows from the fact that the monodromy of $\tilde{\Omega}$ around
$\partial_iB$ is Hamiltonian--equivalent to the Salamon monodromy
$\Phi_{d,\omega_i,\tau}$.)
\end{proof}

In order to extend the previous lemma from genuine fibrations to
Lefschetz fibrations, consider now an elementary Lefschetz fibration
$p\co E\to D^2$, that is, a Lefschetz fibration over the disc with
just one singular fiber $E_0$, which lies over the origin and
contains just one node. If $Y=\partial E$, the image of the
restriction map $H^2(E;\mathbb{Z})\to H^2(Y;\mathbb{Z})$ then has
rank one, generated by $(2-2g)^{-1}e(T^{vt}Y)$.  In particular,
letting $c\in H^2(E;\mathbb{R})$ and $h\in H_2(E,\partial
E;\mathbb{Z})$ meet the fibers positively, the class
$[\omega_Y]=c|_Y+\frac{2\pi}{\tau}PD(\partial h)\in
H^2(Y;\mathbb{R})$ of the form $\omega_Y$ on $Y$ used in the
construction of the Salamon monodromy map for the object
$(Y,p|_Y,\partial h,c,\tau)$ will be some negative (since $g\geq 2$)
multiple of $e(T^{vt}Y)$.

According to \cite[Lemma 3.15]{P2} (which for our present purposes
replaces an erroneous lemma in \cite{DS} which we had referred to in
earlier versions of this paper; we thank the referee for pointing
out this issue and suggesting a way of resolving it), where $\eta\in
H^2(E;\mathbb{Z})$ restricts to the singular fiber $E_0$ as the
orientation class (so $\eta$ is a positive multiple of
$c+\frac{2\pi}{\tau}PD(h)$), any cohomology class in
$H^2(E_d(p);\mathbb{Z})$ of form \[
s\uparrow\eta+t(d\uparrow\eta-(\pi_2)_!(\frak{d}^2)/2)\] with
$s,t>0$ is represented by K\"ahler forms.  In particular, for
$\tau>2\pi d$, there are K\"ahler forms on $E_d(p)$ whose
restrictions to the boundary $Y_d(p)$ represent the class
$[(\Omega_{d,\omega_Y,\tau})_{\Phi_{d,\omega_Y,\tau}}]=2\pi(\tau\uparrow[\omega_Y]-\pi(\pi_2)_!(\frak{d}^2))$
of \fullref{cohclass}.

This positions us
to state the result underlying the construction of the maps induced
by fibered cobordisms.

\begin{prop}\label{rhsform}  Let $f\co X\to B$ be a Lefschetz fibration on a
4--manifold $X$ over a starred surface $B$ with boundary $\partial
B=\partial_1 B\cup\cdots \cup \partial_n B$ such that the critical
values of $f$ are precisely  the interior special points of $B$.
Write $Y_i=f^{-1}(\partial_i B)$, and let $c_i\in
H^2(f^{-1}(\partial_i B);\mathbb{R})$, $h_i\in H^1(f^{-1}(\partial_i
B);\mathbb{Z})$, and $b_j\in H^2(f^{-1}(\beta_j);\mathbb{R})$ be
such that there exist $\tilde{c}\in H^2(X;\mathbb{R})$,
$PD(\tilde{h})\in H^2(X;\mathbb{Z})$ both evaluating positively on
the fibers (say $\langle PD(\tilde{h}), [\fiber]\rangle=d$) with
$c_i=\tilde{c}|_{Y_i}$, $PD(h_i)=PD(\tilde{h})|_{Y_i}$, and
$b_j=(\tilde{c}+2\pi PD(\tilde{h})/\tau)|_{f^{-1}(\beta_j)}$. Let
$\omega_i\in \Omega^2(Y_i)$ be closed fiberwise symplectic forms
representing $c_i+2\pi PD(h_i)/\tau=(c+2\pi PD(h)/\tau)|_{Y_i}\in
H^2(Y_i;\mathbb{R})$. Then for  $\tau>2\pi d$ there is a fiberwise
symplectic form $\tilde{\Omega}=\Omega_{\vec{c},\vec{h},\vec{b}}\in
\Omega^2(X_d(f))$, determined canonically up to isotopy by the
$c_i,h_i,\tau$, and $b_j$ and independent of the choices of
$\tilde{c}$ and $\tilde{h}$, such that
\[
[\tilde{\Omega}|_{(Y_i)_d(f)}]=[(\Omega_{d,\omega_i,\tau})_{\Phi_{d,\omega_i,\tau}}]\]
for each $i$.

\end{prop}

\begin{proof}
Let $p_1,\ldots,p_k$ be the critical values of $f$, so that since
the $p_i$ are the interior special points of $B$ our data include
arcs from the interior base point $b$ of $B$ to the $p_i$.  Let
$U_1,\ldots,U_k$ small disjoint neighborhoods of $p_i$, and take an
almost complex structure on $X$ with respect to which $f\co X\to B$
is pseudoholomorphic and whose restriction to each $f^{-1}(U_i)$ is
integrable.  Write $B^0=B\setminus\cup U_i$ and $F\co X_d(f)\to B$
for the map defining the relative Hilbert scheme.  $B^0$ then
inherits from $B$ the structure of a starred surface with boundary;
$B^0$ has no interior special points, and has new boundary
components $\partial U_i$ corresponding to the special points of
$B$, with the arcs from $b$ to the new boundary components just
given by the portions of the arcs from $b$ to the special points of
$B$ that lie in $B^0$.

We now apply \fullref{relfib} to $B^0$ (taking for the boundary
data $c,PD(h)$ on the new boundary components $f^{-1}(\partial U_i)$
appropriate multiples of $e(T^{vt}f^{-1}(\partial U_i))$). This
gives us a closed fiberwise symplectic form
$\tilde{\Omega}^{0}=\tilde{\Omega}^{0}_{\vec{c},\vec{h},\vec{b}}$ on
$F^{-1}(B^0)$, while the remarks before the proposition give us
K\"ahler forms $\tilde{\Omega}^i$ on each $F^{-1}(U_i)$; furthermore
the restrictions of $\tilde{\Omega}^{0}$ and $\tilde{\Omega}^i$ to
$F^{-1}(\partial U_i)$ are cohomologous and have restrictions which
are compatible with the complex structures on each fiber of
$F^{-1}(\partial U_i)\to S^1$. But then it is straightforward to
glue these together to obtain the desired form on the total space:
letting $t$ be the first coordinate on $[-1,1]\times F^{-1}(U_i)$,
$\rho\co [-1,1]\to [0,1]$ a monotone smooth function which is
identically $0$ near $-1$ and identically $1$ near $1$, and
$\alpha_i\in \Omega^1(F^{-1}(\partial U_i))$ such that
$\tilde{\Omega}^i=\tilde{\Omega}^{0}+d\alpha_i$, the form
$\tilde{\Omega}^{0}+d(\rho(t)\alpha_i)$ will be a closed form equal
to $\tilde{\Omega}^{0}$ for $t$ near $-1$ and to $\tilde{\Omega}^i$
for $t$ near $1$, and its restriction to each fiber will be a convex
combination of K\"ahler forms and so will be K\"ahler.
\end{proof}

With this preparation, the existence of the maps promised in \fullref{cobmaps} follows from an application of the fairly standard
idea that a symplectic cobordism gives rise to maps on Floer
homology groups. Let $\frak{m}=(X,\tilde{f},\tau)$ be a morphism
between objects $(Y_{\pm},f_{\pm},h_{\pm},c_{\pm},\tau)$ in
$\mathsf{FCOB}$, so that in particular $\tilde{f}\co X\to B$ is a
Lefschetz fibration with boundary components $Y_{\pm}$ (either of
which is allowed to be empty), and the sets \[
H_{h_-,h_+}=\{\tilde{h}\in H_2(X,\partial
X;\mathbb{Z})|\partial_{\pm}\tilde{h}=h_{\pm}\in
H_1(Y_{\pm};\mathbb{Z})\} \] $$ C_{c_-,c_+}=\{\tilde{c}\in
H^2(X;\mathbb{Z})|\tilde{c}|_{Y_{\pm}}=c_{\pm}\}\leqno{\hbox{and}}$$ 
are nonempty.
Further, we have, as part of the data of the morphism, classes
$b_j\in H^2(f^{-1}(\beta_j);\mathbb{R})$ as in \fullref{rhsform}. From that Proposition, we then obtain a fiberwise
symplectic form $\tilde{\Omega}$ on $X_d(\tilde{f})$ restricting to
a fiberwise K\"ahler form cohomologous to the Salamon form
$\Omega_{d,\omega_{\pm},\tau}$ on the boundary components, where
$\omega_{\pm}\in \Omega^2(Y_{\pm})$ is a fiberwise symplectic
representative of $c_{\pm}+\frac{2\pi}{\tau}PD(h_{\pm})$.  Let
$\bar{f}\co \bar{X}\to \bar{B}$ be the Lefschetz fibration obtained
by adding trivial cylindrical ends $(-\infty,-1]\times Y_-\to
(-\infty,-1]\times S^1$ and $[1,\infty)\times Y_+\to
[1,\infty)\times S^1$ to $X\to B$ (of course if one or both of
$Y_{\pm}$ is the empty object we don't add an end corresponding to
that object).  Then our usual device involving a cutoff function
gives us a closed, fiberwise symplectic form $\bar{\Omega}$ on
$\bar{X}_d(\bar{f})$ extending $\tilde{\Omega}\in
\Omega^2(X_d(\tilde{f}))$ and equal to the pullback of
$\Omega_{d,\omega_-,\tau}$ on $(-\infty,-2]\times (Y_-)_d(f)$ and to
the pullback of $\Omega_{d,\omega_+,\tau}$ on $[2,\infty)\times
(Y_+)_d(f)$.

Note that our construction of $\tilde{\Omega}$ is compatible with
the composition of morphisms, in the sense that if $\frak{m}_0\in
\Mor(\frak{o}_0,\frak{o}_1)$ and $\frak{m}_1\in
\Mor(\frak{o}_1,\frak{o}_2)$, then after isotoping the forms
$\tilde{\Omega}_i\in \Omega^2((X_i)_{d}(\tilde{f}_i))$ obtained from
\fullref{rhsform} to coincide near their common boundary
component $(Y_i)_d(f_i)$, the form $\tilde{\Omega}$ on the relative
Hilbert scheme associated to $\frak{m}_1\circ \frak{m}_0$ is (up to
isotopy) obtained by gluing together the forms $\tilde{\Omega_i}$
coming from the two pieces.  This compatibility property  is the
main motivation for the additional technical data that we have
included in our definition of a morphism.

Now if $\omega_{\bar{B}}$ is a volume form on the base and
$\bar{F}\co \bar{X}_d(\bar{f})\to \bar{B}$ is the map defining the
relative Hilbert scheme, then
$\bar{\Omega}+K\bar{F}^*\omega_{\bar{B}}$ will be a symplectic form
for large enough $K\in\mathbb{R}$; let $J$ be an almost complex
structure on $\bar{X}_d(\bar{f})$ which is compatible with this
symplectic structure, which makes $\bar{F}$ a pseudoholomorphic map,
and which agrees with the standard complex structure on the relative
Hilbert scheme on the preimages of small neighborhoods of each of
the critical points of $\bar{f}$.  We shall define our maps on the
Floer homology groups by counting certain $J$--holomorphic sections
of $\bar{F}$ with prescribed asymptotics in $(Y_{\pm})_d(f)$.

To be more specific about which sections are counted and how,
consider any class $\tilde{h}\in H_{h_-,h_+}\subset H_2(X,\partial
X;\mathbb{Z})$. Now the fiber product $X_d(f)\times_B X$ contains a
universal divisor $\mathcal{D}$ (for regular values $t\in B$ of
$X\to B$, $\mathcal{D}$ meets the fiber $S^df^{-1}(t)\times
f^{-1}(t)$ over $t$ in $\{(D,p)|p\in D\}$; the extension of
$\mathcal{D}$ over the critical values follows from the
algebro-geometric description of the relative Hilbert scheme of an
elementary Lefschetz fibration, the details of which will not be
relevant to us here).  Hence we get a map $\downarrow\co
H_2(X_d(f),\partial X_d(f);\mathbb{Z})\to H_2(X,\partial
X;\mathbb{Z})$ analogous to the map on relative symmetric products
considered earlier.  Where $\Gamma(X_d(f))$ denotes the space of
sections of $X_d(f)$, there is a natural evaluation map
$\pi_0(\Gamma(X_d(f)))\to H_2(X_d(f),\partial X_d(f);\mathbb{Z})$.
Composing this with $\downarrow$ gives a map
\begin{align*} \pi_0(\Gamma(X_d(f))&\to H_2(X,\partial
X;\mathbb{Z})\\\gamma&\mapsto h_{\gamma}; \end{align*} Lemma 4.1 of
\cite{Smith} shows that this map is injective.
\vspace{8pt}

Our maps $F_{\frak{m},\theta,\tilde{h}}$ in \fullref{cobmaps}
will be the maps induced on homology of certain chain maps
$\tilde{F}_{\frak{m},\theta,\tilde{h}}$ on the Floer chain
complexes. As always taking $d=\tilde{h}\cap[\fiber]$, if our class
$\tilde{h}\in H_2(X,\partial X;\mathbb{Z})$ is not in the image of
$(\gamma\mapsto h_{\gamma})$ we then define
$\tilde{F}_{\frak{m},\theta,\tilde{h}}$ to be zero; otherwise let
$\gamma_{\tilde{h}}$ denote the unique preimage of $\tilde{h}$ under
$(\gamma\mapsto h_{\gamma})$.
$\tilde{F}_{\frak{m},\theta,\tilde{h}}$ will then be constructed
using pseudoholomorphic sections of $\bar{F}\co \bar{X}_d(f)\to
\bar{B}$ in the homotopy class $\gamma_{\tilde{h}}$, with weights
depending on $\theta$ and $h$.
\vspace{8pt}

We now define the promised maps
\[\tilde{F}_{\frak{m},\theta,\tilde{h}}(U^r\otimes\eta_1\wedge\cdots\wedge\eta_k\otimes\cdot)\co
CF(\frak{o}_-;\Lambda_{\Nov}^{R}\otimes\Gamma_{\theta})\to
CF(\frak{o}_+;\Lambda_{\Nov}^{R}\otimes\Gamma_{\theta})\] where
$r\geq 0$ and $\eta_i\in H_1(X;\mathbb{Z})$. Assume first that
neither of $\frak{o}_{\pm}$ is the empty object. Then
$CF(\frak{o}_{\pm};\Lambda_{\Nov}^{R})$ is freely generated over
$\Lambda_{\Nov}^{R}$ by the ``constant sections'' $\mathbf{x}$ of
$(Y_{\pm})_d(f)$ corresponding to the fixed points of Salamon's
symplectomorphism $\Phi_{d,\omega_{\pm},\tau}$.  If $\mathbf{x}_-$
is a generator of $CF(\frak{o}_-;\Lambda_{\Nov}^{R})$ and
$\mathbf{x}_+$ is a generator of $CF(\frak{o}_+;\Lambda_{\Nov}^{R})$,
and if $J$ is an almost complex structure on $\bar{X}_d(\bar{f})$
which is compatible with the symplectic structure, which makes the
projection $\bar{X}_d(\bar{f})\to \bar{B}$ pseudoholomorphic for
some chosen complex structure on $\bar{B}$, and which agrees with
the K\"ahler structure near the singular fibers, then let
\[\mathcal{M}_{J,\frak{m},\tilde{h}}(\mathbf{x_-},\mathbf{x_+})\] denote the
moduli space of $J$--holomorphic sections of $\bar{X}_d(\bar{f})$
which \begin{itemize} \item[(i)] are asymptotic to a cylinder on
$\mathbf{x_-}$ on the $(Y_-)_d(f_-)$ end of
$\bar{X}_d(\bar{f})$,\item[(ii)] are asymptotic to a cylinder on
$\mathbf{x_+}$ on the $(Y_+)_d(f_+)$ end of $\bar{X}_d(\bar{f})$,
and \item[(iii)] represent the homotopy class
$\gamma_{\tilde{h}}$.\end{itemize}
\vspace{8pt}

Provided that $d\notin [(g+1)/2,g-1)$, so that the
$w^+$--monotonicity condition of \cite{L1} holds (implying that
moduli spaces will generically not contain bubble trees involving
multiply covered spheres of negative Chern
number\footnote{Superficially, one would also need to rule out
bubbles in the singular fibers of the relative Hilbert scheme;
however, \cite[Lemma 4.8]{DS} shows that any such bubble would be
homologous to a sphere in a smooth fiber, so such bubbles do not
complicate the analysis.}),
$\mathcal{M}_{J,\frak{m},\tilde{h}}(\mathbf{x_-},\mathbf{x_+})$
will, for generic $J$, be a smooth manifold of a certain dimension
$\delta(\tilde{h})$ which admits a compactification
$\bar{\mathcal{M}}_{J,\frak{m},\tilde{h}}(\mathbf{x_-},\mathbf{x_+})$
by
\begin{itemize} \item[(i)] ``broken sections'' consisting of a chain
of Floer flowlines in $(Y_-)_d(f_-)$, followed by a section of
$\bar{X}_d(\bar{f})$, followed by a chain of Floer flowlines in
$(Y_+)_d(f_+)$.  Each section in the entire chain begins where the
previous one ends, and the homotopy class of the whole sequence is
$\gamma_{\tilde{h}}$; \item[(ii)] ``cusp sections'' consisting of a
section of $\bar{X}_d(\bar{f})$, together with spherical bubbles in
various fibers of $\bar{X}_d(\bar{f})$; and \item[(iii)]
combinations of (i) and (ii).\end{itemize}  These additional strata
will have codimension at least 1, and the only codimension one
strata will be those made up of broken sections with just one Floer
flowline component.  When $d\in [(g-1)/2,g-1)$, so that virtual
moduli methods are required, we can still find rational chains in an
ambient space containing
$\mathcal{M}_{J,\frak{m},\tilde{h}}(\mathbf{x_-},\mathbf{x_+})$
which satisfy properties parallel to these; see \cite{LT}.
\vspace{8pt}

The dimension $\delta(\tilde{h})$ is obtained as follows (see
\cite[Theorem 3.3.11]{Schwarz}): take an arbitrary  section $s$ of
$X_d(f)$ representing the class $\gamma_{\tilde{h}}$.
$s^*T^{vt}X_d(f)$ is then a symplectic bundle over the open
$2$--manifold $B$ and so is trivial; choosing any trivialization
$\tau$ we can then compute the Conley--Zehnder indices
$\mu_{CZ}^{\tau}(\mathbf{x}_{\pm})$ of
$T^{vt}(Y_{\pm})_d(f_{\pm})|_{\mathbf{x_{\pm}}}$ with respect to the
trivialization; then
\[
\delta(\tilde{h})=\mu_{CZ}^{\tau}(\mathbf{x}_+)-\mu_{CZ}^{\tau}(\mathbf{x}_-)+2d\chi
(B)\] is independent of the various choices and is the desired
virtual dimension.

Now if $\alpha\in \Omega^2(\bar{X}_d(\bar{f}))$ is any closed form,
then the integral of $\alpha$ over any cylinder in
$\bar{X}_d(\bar{f})$ representing an element of
$\mathcal{M}_{J,\frak{m},\tilde{h}}(\mathbf{x_-},\mathbf{x_+})$
depends only on $\tilde{h},\mathbf{x_-},$ and $\mathbf{x_+}$; denote
this common value by $\alpha(\tilde{h},\mathbf{x_-},\mathbf{x_+})$.
This remark in particular applies to both the forms $\uparrow
\theta$ (which we can easily make sense of thanks to the fact that
we assumed in the statement of \fullref{cobmaps} that $\theta$
vanishes near the critical points of $f$) and to the form
$\bar{\Omega}$ introduced after \fullref{rhsform}.  For
$\mathbf{x}_-\in \mathcal{A}_-$, we set
\begin{align*} \tilde{F}_{\frak{m},\theta,\tilde{h}}&(U^k\otimes\eta_1\wedge\cdots\wedge\eta_l\otimes
\mathbf{x_-})
\\&=\!\sum_{\mathbf{x_+}\in\mathcal{A}_+}\!T^{\bar{\Omega}(\tilde{h},\mathbf{x_-},\mathbf{x_+})}e^{\uparrow\theta(\tilde{h},\mathbf{x_-},\mathbf{x_+})}
\langle\mathbf{x_-}|U^r\otimes\eta_1\wedge\cdots\wedge\eta_k|\mathbf{x_+}\rangle,\\
&\hspace{1in}\langle\mathbf{x_-}|U^r\otimes\eta_1\wedge\cdots\wedge\eta_k|\mathbf{x_+}\rangle\tag*{\hbox{where}}\end{align*}
is defined by: \begin{itemize} \item if $\delta(\tilde{h})\neq
2k+l$, then
$\langle\mathbf{x_-}|U^r\otimes\eta_1\wedge\cdots\wedge\eta_k|\mathbf{x_+}\rangle=0$.
\item if $\delta(\tilde{h})=2k+l$, let $A_1,\ldots,A_r$ be generic
representatives of the class $PD(\uparrow PD[pt])\in
H^{2d-2}(\bar{X}_d(\bar{f}))$, and let
$D_{\eta_1},\ldots,D_{\eta_k}$ be generic representatives of the
classes $PD(\uparrow PD(\eta_i))\in H^{2d-1}(\bar{X}_d(\bar{f}))$.
The set
\begin{align*}
\mathcal{M}_{J,\frak{m},\tilde{h}}(&A_1,\ldots,A_r,D_1,\ldots,D_k;\mathbf{x}_-,\mathbf{x}_+)=\\
&\{(u,x_1,\ldots,x_r,y_1,\ldots,y_k)\in
\mathcal{M}_{J,\frak{m},\tilde{h}}(\mathbf{x}_-,\mathbf{x}_+)\times
B^r\times B^k|\\&\qquad\qquad\qquad\qquad \qquad\qquad \qquad
u(x_i)\in A_i,u(y_j)\in D_j \,\forall i,j\}\end{align*} will have
virtual dimension $\delta(\tilde{h})+2(2k+l)-4k-3l=0$, and we let
$\langle\mathbf{x_-}|U^r\otimes\eta_1\wedge\cdots\wedge\eta_k|\mathbf{x_+}\rangle$
be the signed number of elements in this set for generic $J$
(counted virtually as in \cite{LT,LTprod} if necessary).
\end{itemize}

Since the chains $A_i$ have even codimension and the $D_i$ have odd
codimension, this construction induces a map \[
\tilde{F}_{\frak{m},\theta,h}\co \mathbb{A}(X) \otimes
CF(\frak{o}_-;\Lambda_{\Nov}^{R}\otimes\Gamma_{\theta})\to
CF(\frak{o}_+;\Lambda_{\Nov}^{R}\otimes\Gamma_{\theta}),\] which
might depend on the additional choices (the finiteness condition for
the Novikov ring is trivial, since there are only finitely many
pairs $(\mathbf{x}_-,\mathbf{x}_+)$ and so
$\bar{\Omega}(\tilde{h},\mathbf{x_-},\mathbf{x_+})$ takes just
finitely many values for a given choice of $\tilde{h}$).

\begin{lemma} \label{basiccob}  Given a morphism $\frak{m}=(X,\tilde{f},\tau)$
between the objects
$$\frak{o}_{\pm}=(Y_{\pm},f_{\pm},h_{\pm},c_{\pm},\tau)$$ and a closed
form $\theta\in\Omega^2(X)$ which vanishes near the critical points
of $X$,  for each $\tilde{h}\in H_{h_-,h_+}$, the map
\[ \tilde{F}_{\frak{m},\theta,\tilde{h}}\co
\mathbb{A}(X)\otimes
CF(\frak{o}_-;\Lambda_{\Nov}^{R}\otimes\Gamma_{\theta|_{Y_-}})\to
CF(\frak{o}_+;\Lambda_{\Nov}^{R}\otimes\Gamma_{\theta|_{Y_+}})\]  is
a chain map, and the induced map on homology \[
F_{\frak{m},\theta,\tilde{h}}\co \mathbb{A}(X)\otimes
HF(\frak{o}_-;\Lambda_{\Nov}^{R}\otimes\Gamma_{\theta|_{Y_-}})\to
HF(\frak{o}_+;\Lambda_{\Nov}^{R}\otimes\Gamma_{\theta|_{Y_+}}) \] is
independent of the choices  of chains $A_i,D_j$ and of the almost
complex structure $J$.  The sum
$\tilde{F}_{\frak{m},\theta}=\sum_{\tilde{h}\in
H_{h_-,h_+}}\tilde{F}_{\frak{m},\theta,\tilde{h}}$ is a well-defined
chain map $CF(\frak{o}_-;\Lambda_{\Nov}^{R})\to
CF(\frak{o}_+;\Lambda_{\Nov}^{R})$.   Furthermore if $X=[0,1]\times
Y_-$ is the trivial cobordism, then $F_{\frak{m},\theta,\tilde{h}}$
coincides with the map defining the $\mathbb{A}(Y)$--module
structure of
$HF(\frak{o}_-;\Lambda_{\Nov}^{R}\otimes\Gamma_{\theta_{Y_-}})$.
\end{lemma}

\begin{proof}  To see that $\sum_{\tilde{h}\in
H_{h_-,h_+}}\tilde{F}_{\frak{m},\theta,\tilde{h}}$ is well-defined
we just need to check the finiteness condition for the Novikov ring.
However this follows directly from Gromov compactness, which ensures
that for any $c$ the set \[ \{u\in
\cup_{\tilde{h},\mathbf{x}_-,\mathbf{x}_+}
\bar{\mathcal{M}}_{J,\frak{m},\tilde{h}}(\mathbf{x_-},\mathbf{x_+})|\int_B
u^*\bar{\Omega}<c\}\] is compact, together with the definition that
$\Omega(\tilde{h},\mathbf{x}_-,\mathbf{x}_+)=\int_B u^*\bar{\Omega}$
for any section $u$ in the homotopy class $\gamma_{\tilde{h}}$ which
is asymptotic to $\mathbf{x}_{\pm}$.
\vspace{8pt}

Now for the trivial cobordism, $\bar{\Omega}$ is just the pullback
of the Salamon form $\Omega_{d,\omega,\tau}$ on $Y_d(f)$ to
$\mathbb{R}\times Y_d(f)$.  Also, using the exactness of
$H_2(\partial I\times Y;\mathbb{Z})\to H_2(I\times Y;\mathbb{Z})\to
H_2(I\times Y,\partial I\times Y;\mathbb{Z})\to H_1(\partial I\times
Y;\mathbb{Z})$ and the fact that the first map is a surjection, the
difference between any two elements of $H_{h_-,h_+}$ lies in the
image of $H_2(I\times Y;\mathbb{Z})\to H_2(I\times Y,\partial
I\times Y;\mathbb{Z})$ and so is zero.  So in this case
$H_{h_-,h_+}$ is a singleton $\{\tilde{h}\}$ and
$\tilde{F}_{\frak{m},\theta}=\tilde{F}_{\frak{m},\theta,\tilde{h}}$.
The complex
$CF(Y,f,h,c,\tau;\Lambda_{\Nov}^{R}\otimes\Gamma_{\theta})$ is
generated over $\Lambda_{\Nov}^{R}$ by the fixed points $\mathbf{x}$,
and the differential just counts holomorphic sections $u$ of
$\mathbb{R}\times Y_d(f)$, weighted by $T^{\int_{\mathbb{R}\times
S^1}u^*\Omega_{d,\omega,\tau}}e^{\int_{\mathbb{R}\times
S^1}u^*(\uparrow\theta)}$.  In particular, \[
\tilde{F}_{\frak{m},\theta,\tilde{h}}(1\otimes 1\otimes\cdot)\co
CF(\frak{o};\Lambda_{\Nov}^{R}\otimes \Gamma_{\theta})\to
CF(\frak{o};\Lambda_{\Nov}^{R}\otimes \Gamma_{\theta})\] is none
other than the identity (corresponding to sections of form
$\mathbb{R}\times \{\mathbf{x}\}$) plus the Floer boundary operator.
Furthermore, by taking all the chains $A_i,D_j\in
C_*(\mathbb{R}\times Y_d(f);\mathbb{Z})$ to be contained in
$\{0\}\times Y_d(f)$, comparing with section 4 reveals that
\[ \tilde{F}_{\frak{m},\theta,\tilde{h}}\co
\mathbb{Z}[U]\otimes \Lambda^{*}(H_1(\mathbb{R}\times
Y;\mathbb{Z})/\torsion)\otimes
CF(\frak{o};\Lambda_{\Nov}^{R}\otimes\Gamma_{\theta})\to
CF(\frak{o}_+;\Lambda_{\Nov}^{R}\otimes\Gamma_{\theta})\] is a chain
map which induces the $\mathbb{Z}[U]\otimes
\Lambda^*(H_1(Y;\mathbb{Z})/\torsion)$--module structure on
$HF(\frak{o};\Lambda_{\Nov}^{R}\otimes \Gamma_{\theta})$.  This
proves the last statement of the lemma.
\vspace{8pt}

To see that, for general morphisms $\frak{m}$ and for each
$\tilde{h}\in H_{h_-,h_+}$,  $\tilde{F}_{\frak{m},\theta,\tilde{h}}$
defines a chain map, let $\partial_{\pm}$ be the Floer boundary
operators for $\frak{o}_{\pm}$. We can then write, for any
generators $\mathbf{x}_{\pm}$ of $\frak{o}_{\pm}$,
\[
\langle
\partial_+\tilde{F}_{\frak{m},\theta,\tilde{h}}(U^r\otimes\eta_1\wedge\cdots\wedge\eta_k\otimes\mathbf{x}_-),\mathbf{x}_+\rangle=\Theta_{\tilde{h}}(\mathbf{x}_-,\mathbf{x}_+)T^{\bar{\Omega}(\tilde{h},\mathbf{x}_-,\mathbf{x}_+)}
e^{\uparrow\theta(\tilde{h},\mathbf{x}_-,\mathbf{x}_+)}\]
and
\[ \langle \tilde{F}_{\frak{m},\theta,\tilde{h}}(U^r\otimes\eta_1\wedge\cdots\wedge\eta_k\otimes\partial_-\mathbf{x}_-),\mathbf{x}_+\rangle=\Psi_{\tilde{h}}(\mathbf{x}_-,\mathbf{x}_+)T^{\bar{\Omega}(\tilde{h},\mathbf{x}_-,\mathbf{x}_+)}
e^{\uparrow\theta(\tilde{h},\mathbf{x}_-,\mathbf{x}_+)}.\]

Here all
terms are zero except when $\delta({\tilde{h}})=2k+l+1$, in which
case $\Theta_{\tilde{h}}(\mathbf{x}_-,\mathbf{x}_+)$ is the signed
count of broken sections from $\mathbf{x}_-$ to $\mathbf{x}_+$
consisting of a section of $\bar{X}_d(\bar{f})$ in the relative
homotopy class  $\gamma_{\tilde{h}}$ followed by a flowline for
$CF(\frak{o}_+)$ and which satisfy incidence conditions
corresponding to $r$ and the $\eta_i$, while
$\Psi_{\tilde{h}}(\mathbf{x}_-,\mathbf{x}_+)$ is the signed count of
broken sections from $\mathbf{x}_-$ to $\mathbf{x}_+$ consisting of
a flowline for $CF(\frak{o}_-)$ followed by a section of
$\bar{X}_d(\bar{f})$ in the class $\gamma_{\tilde{h}}$ which satisfy
these same incidence conditions.  Note that if $s_0\#s_1$ is a
broken section from $\mathbf{x}_-$ to $\mathbf{x}_+$ consisting of a
section $s_0$ of $\bar{X}_d(\bar{f})$ asymptotic at its negative end
to $\mathbf{x}_-$ and at its positive end to a generator (say
$\mathbf{y}$) of $CF(\frak{o}_+)$, followed by a flowline $s_1$ for
$CF(\frak{o}_+)$ from $\mathbf{y}$ to $\mathbf{x_+}$ (the latter
which extends to a map $[-\infty,\infty]\times S^1\to
(Y_+)_d(f_+)$), then extending $s_0$ by concatenating it with
$s_1|_{[-\infty,T]\times S^1}$ for $T\in \mathbb{R}\cup
\{-\infty,\infty\}$ defines a homotopy rel $(Y_+)_d(f_+)$ between
$s_0$ and the broken section $s_0\# s_1$.  Hence for any such
$s_0\#s_1$, $s_0$ belongs to the relative homotopy class
$\gamma_{\tilde{h}}$ if and only if the broken section $s_0\#s_1$
does, and so $\Theta_{\tilde{h}}(\mathbf{x}_-,\mathbf{x}_+)$ may
equally well be described as the signed count of broken sections
from $\mathbf{x}_-$ to $\mathbf{x}_+$ consisting of a section of
$\bar{X}_d(\bar{f})$ in followed by a flowline for $CF(\frak{o}_+)$
such that the (concatenated) broken section belongs to
$\gamma_{\tilde{h}}$ and satisfies certain incidence conditions; a
similar description applies to
$\Psi_{\tilde{h}}(\mathbf{x}_-,\mathbf{x}_+)$. But then a standard
argument shows that
$\Psi_{\tilde{h}}(\mathbf{x}_-,\mathbf{x}_+)=\Theta_{\tilde{h}}(\mathbf{x}_-,\mathbf{x}_+)$,
for their difference counts the oriented number of boundary points
of the $1$-manifold
$\mathcal{M}_{J,\frak{m},\tilde{h}}(A_1,\ldots,A_r,D_1,\ldots,D_k;\mathbf{x}_-,\mathbf{x}_+)$.
\vspace{8pt}

Thus
$\partial_+\tilde{F}_{\frak{m},\theta,\tilde{h}}=\tilde{F}_{\frak{m},\theta,\tilde{h}}\partial_-$,
and, summing over $\tilde{h}\in H_{h_-,h_+}$,
$\partial_+\tilde{F}_{\frak{m},\theta}=\tilde{F}_{\frak{m},\theta}\partial_-$.
The fact that the induced maps $F_{\frak{m},\theta,\tilde{h}}$ on
homology are independent of $J$ and of the choices of chains giving
the incidence conditions now follows by standard cobordism
arguments: for example if $D_1,D'_1$ are both representatives of
$PD(\uparrow PD(\eta_1))$ (say $D_1-D'_1=\partial E$), then a
consideration of the boundary components of the moduli spaces
$\mathcal{M}_{J,\frak{m},\tilde{h}}(A_1,\ldots,A_r,E,D_2,\ldots,D_k)$
shows that replacing $D_1$ by $D'_1$ in the definition of
$\tilde{F}_{\frak{m},\theta,\tilde{h}}$ does not affect the induced
map $F_{\frak{m},\theta,\tilde{h}}$ on homology.
\end{proof}

We have been assuming that our morphism
$\frak{m}=(X,\tilde{f},\tau)$ is a morphism between nonempty objects
$\frak{o}_{\pm}=(Y_{\pm},f_{\pm},h_{\pm},c_{\pm},\tau)$; the
construction of $F_{\frak{m},\theta,\tilde{h}}$ in the case in which
one or both of the fibered 3--manifolds $Y_{\pm}$ underlying
$\frak{o}_{\pm}$ is empty is a simple modification of what we have
already done. If the ``incoming'' boundary component $Y_-$ is empty
but $Y_+\neq\varnothing$, then we are to define a map \[
F_{\frak{m},\theta,\tilde{h}}\co
\mathbb{A}(X)\otimes\Lambda_{\Nov}^{R}\to
HF(\frak{o}_+;\Lambda_{\Nov}^{R}\otimes\Gamma_{\theta});\] this is
done as before by first defining maps
$\tilde{F}_{\frak{m},\theta,\tilde{h}}$ to
$CF(\frak{o}_+;\Lambda_{\Nov}^{R}\otimes\Gamma_{\theta})$ for classes
$\tilde{h}\in H_2(X,Y_+;\mathbb{Z})$ such that
$\partial\tilde{h}=h_+$ and $\delta(\tilde{h})=2r+k$ by
\begin{align*}
\tilde{F}_{\frak{m},\theta,\tilde{h}}(U^r&\otimes\eta_1\wedge\cdots\wedge\eta_k\otimes
1)\\=&\sum_{\mathbf{x}_+\in\mathcal{A}_+}T^{\bar{\Omega}(\tilde{h},\mathbf{x}_+)}e^{\uparrow\theta(\tilde{h},\mathbf{x}_+)}\#\mathcal{M}_{J,\frak{m},\tilde{h}}(A_1,\ldots,A_r,D_1,\ldots,D_k;\mathbf{x}_+),\end{align*}
where the notation is as before; of course since $B$ now only has
one boundary component there is only one asymptotic condition to
specify in expressions such as
$\bar{\Omega}(\tilde{h},\mathbf{x}_+)$ and
$\#\mathcal{M}_{J,\frak{m},\tilde{h}}(A_1,\ldots,A_r,D_1,\ldots,D_k;\mathbf{x}_+)$.
Consideration of the boundary of the spaces
$\mathcal{M}_{J,\frak{m},\tilde{h}}(A_1,\ldots,A_r,D_1,\ldots,D_k;\mathbf{x}_+)$
when $\delta(\tilde{h})=1$ reveals that $\partial_+
{F}_{\frak{m},\theta,\tilde{h}}=0$; $F_{\frak{m},\theta,\tilde{h}}$
is then the map on homology induced by this sum, and is independent
of the additional choices exactly as in \fullref{basiccob}

Dually, if $Y_-\neq\varnothing$ but $Y_+=\varnothing$, and if
$\mathbf{x}_-$ is a generator of $CF(\frak{o}_-)$, we define
$\tilde{F}_{\frak{m},\theta,\tilde{h}}(U^r\otimes\eta_1\wedge\ldots\wedge\eta_k\otimes\mathbf{x}_-)$
by counting pseudoholomorphic sections of $\bar{X}_d(\bar{f})$ which
are asymptotic at the boundary to $\mathbf{x}_-$ and which satisfy
the usual incidence conditions corresponding to $r$,
$\eta_1,\ldots,\eta_k$, with the usual weights
$T^{\bar{\Omega}(\tilde{h},\mathbf{x}_-)}e^{\uparrow\theta(\tilde{h},\mathbf{x}_-)}$;
this is a chain map exactly as in the previous case, and
$F_{\frak{m},\theta,\tilde{h}}\co \mathbb{A}(X)\otimes
HF(\frak{o}_-;\mathbb{Z})\to\Lambda_{\Nov}^{R}$ is the induced map on
homology.

We pause to consider the simplest nontrivial case, where
$Y_-=\varnothing$ and $Y_+=S^1\times \Sigma$ with the obvious
fibration $\pi_1$, and $h_+=d[S^1\times pt]$ with $c$ proportional
to $PD(h)$, so that $\frak{m}=(D^2\times \Sigma,\pi_1,\tau)$ gives a
morphism from $\frak{o}_-=\varnothing$ to
$\frak{o}_+=(Y_+,\pi_1,h_+,c_+,\tau)$.  Now as in \fullref{basicsd}(ii) we have a natural surjection $\Pi\co
\mathbb{Z}[U]\otimes\Lambda^*H_1(D^2\times\Sigma;\mathbb{Z})\twoheadrightarrow
H^*(S^d\Sigma;\mathbb{Z})$.  Meanwhile the monodromy of a fiberwise
symplectic form on $S^1\times \Sigma$ in the class $c_++2\pi
PD(h_+)/\tau$ is Hamiltonian, so the Salamon monodromy
$\Phi_{d,\omega,\tau}\co S^d\Sigma\to S^d\Sigma$ is also
Hamiltonian, so that the Floer homology
$HF(\frak{o}_+;\Lambda_{\Nov}^{R})$ is isomorphic to
$H^*(S^d\Sigma;\Lambda_{\Nov}^{R})$.  It follows directly from the
relevant definitions that our map
\[F_{\frak{m},0}\co
\mathbb{Z}[U]\otimes\Lambda^*H_1(D^2\times\Sigma;\mathbb{Z})\otimes\Lambda_{\Nov}^{R}\to
HF(\frak{o}_+;\Lambda_{\Nov}^{R})\] factors through $\Pi\co
\mathbb{Z}[U]\otimes\Lambda^*H_1(D^2\times\Sigma;\mathbb{Z})\twoheadrightarrow
H^*(S^d\Sigma;\mathbb{Z})$ to give a map\break
$H^*(S^d\Sigma;\Lambda_{\Nov}^{R})\to
HF(\frak{o}_+;\Lambda_{\Nov}^{R})$ which is none other than the
Piunikhin--Salamon--Schwarz isomorphism between the cohomology of
$S^d\Sigma$ and its Hamiltonian Floer homology, constructed in
\cite{PSS}.  Note that the monopole and Heegaard Floer homologies of
$S^1\times\Sigma$ in the corresponding spin$^c$ structure are also
known to be given by $H^*(S^d\Sigma)$ when $d<g-1$
\cite{MW2,OS}.  On the other hand, for $d>g-1$ we obtain
$H^*(S^d\Sigma)$ rather than $H^*(S^{2g-2-d}\Sigma)$ as in
\cite{MW2,OS}; this discrepancy results from the fact that
for this range of $d$ we cannot choose $c$ and $\tau$ to have the
property that the corresponding Seiberg--Witten theory perturbation
class $\eta(h,c,\tau)$ mentioned in the introduction is zero.
\vspace{8pt}

Finally, in the case when $Y_-=Y_+=\varnothing$, so that the
morphism $\frak{m}$ corresponds to a Lefschetz fibration
$\tilde{f}\co X\to B$ on a closed manifold, we define
$F_{\frak{m},\theta,\tilde{h}}\co
\mathbb{A}(X)\otimes\Lambda_{\Nov}^{R}\to \Lambda_{\Nov}^{R}$ by once
again counting pseudoholomorphic sections of the relative Hilbert
scheme $\bar{X}_d(\bar{f})$; since $\partial X=\varnothing$ the
relevant classes $\tilde{h}$ belong to $H_2(X;\mathbb{Z})$, and by
\cite[Proposition 4.3]{Smith} the virtual dimension
$\delta(\tilde{h})$ of the space
$\mathcal{M}_{J,\frak{m},\tilde{h}}$ of pseudoholomorphic sections
of $\bar{X}_f(\bar{f})$ in the homotopy class $\gamma_{\tilde{h}}$
is $\langle PD(\tilde{h})-\kappa_X,\tilde{h}\rangle$.  If
$\delta(\tilde{h})\neq 2r+k$ we put
$F_{\frak{m},\theta,\tilde{h}}(U^r\otimes\eta_1\wedge\cdots\eta_k)=0$;
otherwise we put \[
F_{\frak{m},\theta,\tilde{h}}(U^r\otimes\eta_1\wedge\cdots\eta_k)=T^{\langle[\Omega],\gamma_{\tilde{h}}\rangle}e^{\langle\uparrow[\theta],\gamma_{\tilde{h}}\rangle}\#\mathcal{M}_{J,\frak{m},\tilde{h}}(A_1,\ldots,A_r,D_1,\ldots,D_k).\]
$\#\mathcal{M}_{J,\frak{m},\tilde{h}}(A_1,\ldots,A_r,D_1,\ldots,D_k)$
is independent of $J$ and of the choices of $A_i$ and $D_j$, and
indeed is (by definition) the Donaldson--Smith invariant
\[ DS_{(X,f)}(\tilde{h};pt^r,\eta_1,\ldots,\eta_k).\]  Setting
$F_{\frak{m},\theta}=\sum_{\tilde{h}\in H_2(X;\mathbb{Z})}
F_{\frak{m},\theta,\tilde{h}}$ then completes the definition of the
maps $F_{\frak{m},\theta}$ in all cases, and also establishes part
(iv) of \fullref{cobmaps} (noting that
$\langle\uparrow[\theta],\gamma_{\tilde{h}}\rangle=\langle[\theta],\downarrow\gamma_{\tilde{h}}\rangle=\langle[\theta],\tilde{h}\rangle$).
\vspace{8pt}
  
Along the same lines, we could also use Lefschetz fibrations over
surfaces with an arbitrary number $n$ of boundary components to
obtain ``quantum multiplication'' maps
$HF(\frak{o}_1;\Lambda_{\Nov}^{R})\otimes\cdots\otimes
HF(\frak{o}_{n-1};\Lambda_{\Nov}^{R})\to
HF(\frak{o}_n;\Lambda_{\Nov}^{R})$ for suitable objects $\frak{o}_i$,
but we shall not develop this here.
\vspace{8pt}
   
Now that we have defined the maps $F_{\frak{m},\theta}$, the various
parts of \fullref{cobmaps} follow fairly quickly.  Part (i) is
obtained by standard gluing arguments: if
$\frak{m}=(X,\tilde{f},\tau)=\frak{m}_1\circ\frak{m}_0$ with
$\frak{m}_0=(X_0,\tilde{f}_0,\tau)\in
\Mor(\frak{o}_0,\frak{o}_1),\frak{m}_1=(X_1,\tilde{f}_1,\tau)\in
\Mor(\frak{o}_1,\frak{o}_2)$,  writing
$\frak{o}_1=(Y_1,f_1,h_1,c_1,\tau)$ we may choose the form
$\bar{\Omega}$ on $\bar{X}_d(\bar{f})$ to restrict to a neighborhood
of $(Y_1)_{d}(f_1)$ as the pullback of the Salamon form
$\Omega_{d,\omega,\tau}$; varying the complex structure on the base
$\bar{B}$ (and recalling that the almost complex structure $J$ on
$\bar{X}_d(\bar{f})$ is constrained to make the projection to
$\bar{B}$ pseudoholomorphic, so this also varies $J$) metrically
identifies this neighborhood $N$ of $(Y_1)_d(f_1)$ with
$[-T,T]\times (Y_1)_d(f_1)$ for arbitrarily large $T$. Let $J_T$
denote a generic almost complex structure obtained in this way.
Choose representatives $A_1,\ldots, A_k$ of $PD(\uparrow PD[pt])\in
H_{2d-2}(\bar{X}_d(\bar{f}))$ which are contained in the
``$\frak{m}_1$ side'' $(X_1)_d(\tilde{f}_1)\subset
\bar{X}_d(\bar{f})$ (and are outside the neighborhood $N$ of
$(Y_1)_d(f_1)$ mentioned above), and representatives
$A_{k+1},\ldots,A_{k+l}$ of $PD(\uparrow PD[pt])\in
H_{2d-2}(\bar{X}_d(\bar{f}))$ which are contained in the
``$\frak{m}_2$ side'' $(X_2)_d(\tilde{f}_2)\subset
\bar{X}_d(\bar{f})$ and are also disjoint from $N$.
 Once $T$ is large enough,
it follows from gluing theorems that in special cases date back at
least to \cite{Floer} (the argument in \cite[Section 3.3]{Salsurv}
can be rather directly applied to our case, the only essential
difference being that here we glue at two points rather than one)
that the moduli spaces
$\mathcal{M}_{J_T,\frak{m},\tilde{h}}(A_1,\ldots,A_{k+l};\mathbf{x}_0,\mathbf{x}_2)$
of pseudoholomorphic sections will, for generic $J_T$ and
appropriate generic $J^0,J^1$ on
$(\bar{X}_0)_f(\bar{f}_0),(\bar{X}_1)_f(\bar{f}_1),$ be in a
one-to-one, orientation-preserving correspondence (which preserves
the Novikov--ring and local--coefficient weights) with
\begin{align*}
\smash{\bigcup_{\begin{subarray}{c}\mathbf{x}_{1}^{-},\mathbf{x}_{1}^{+}\in\mathcal{A}_1,\\
+(\tilde{h}_0,\tilde{h}_{Y_1}\tilde{h}_1)=\tilde{h}\end{subarray}}}
\bigg(\mathcal{M}_{J^0,\frak{m}_0,\tilde{h}_0}&(A_1,\ldots,A_k;\mathbf{x}_0,\mathbf{x}_{1}^{-})\times
\mathcal{M}_{J,\mathbb{R}\times(Y_1)_d(f),\tilde{h}_{Y_1}}(\mathbf{x}_{1}^{-},\mathbf{x}_{1}^{+})\times\\
&\mathcal{M}_{J^1,\frak{m}_1,\tilde{h}_1}(A_{k+1},\ldots,A_{k+l};\mathbf{x}_{1}^{+},\mathbf{x}_2)\bigg),\end{align*}
where the notation
$+(\tilde{h}_0,\tilde{h}_Y,\tilde{h}_1)=\tilde{h}$ means that there
exist an asymptotically constant section $s_0$ of
$(\bar{X}_0)_f(\bar{f}_0)$ representing $\gamma_{\tilde{h}_0}$ (say
approaching $\gamma_-\in\Gamma((Y_1)_d(f))$ at the positive end of
$B_0$), a section $s_Y$ of $\mathbb{R}_t\times Y_d(f)$ representing
$\gamma_{\tilde{h}_Y}$ which is asymptotic to
$\gamma_{\pm}\in\Gamma((Y_1)_d(f))$ as $t\to\pm\infty$, and an
asymptotically constant section $s_1$ of $(\bar{X}_1)_f(\bar{f}_1)$
representing $\gamma_{\tilde{h}_1}$ which approaches $\gamma_+$ at
the negative end of $B_1$, and that when these sections are glued
along their corresponding asymptotic limits to obtain a section
$s_0\#s_Y\#s_1$ of $\bar{X}_d(\bar{f})$, $s_0\#s_Y\# s_1$ represents
the homotopy class $\gamma_{\tilde{h}}$ (we allow, of course, $s_Y$
to be the trivial section $\mathbb{R}\times \gamma_-$). Summing over
those $\tilde{h}$ with $h|_{X_i}=\tilde{h}_i$ ($i=0,1$), this
one-to-one correspondence then translates into the language of our
cobordism maps as the statement that \begin{align*}
\smash{\sum_{\begin{array}{c}\scriptstyle \tilde{h}\in
H_{h_0,h_2}\\\scriptstyle
\tilde{h}|_{X_i}=\tilde{h}_i\end{array}}}\tilde{F}&_{\frak{m}_1\circ\frak{m}_0,\theta,\tilde{h}}(U^{k+l}\otimes
1\otimes
\mathbf{x})\\&\qquad=\tilde{F}_{\frak{m}_1,\theta|_{X_1},\tilde{h}_0}\left(U^{k}\otimes
1\otimes (1+\partial_{\frak{o}_1})
\tilde{F}_{\frak{m}_0,\theta|_{X_0},\tilde{h}_1}(U^l\otimes 1\otimes
\mathbf{x})\right),\end{align*} and passing to the induced maps on
homology proves part (i) of \fullref{cobmaps}.

  Part (ii) of
\fullref{cobmaps}, which asserts that where $i_-\co
\mathbb{A}(Y_-)\to \mathbb{A}(X)$ is induced by the inclusion
$Y_-\subset X$ we have
\begin{equation} \label{modcob} F_{\frak{m},\theta,\tilde{h}}(1\otimes\lambda\cdot
\mathbf{x})=F_{\frak{m},\theta,\tilde{h}}(i_-(\lambda)\otimes
\mathbf{x}),
\end{equation} follows from the last part of \fullref{basiccob} and a
similar gluing argument.  Here we use the fact that
$\bar{X}_d(\bar{f})$ contains a half--infinite cylinder on
$(Y_-)_d(f_-)$, take the chains $A_i,D_j$ to be contained in some
fixed $\{t_0\}\times (Y_-)_d(f_-)$ in this cylinder, and send the
length $T$ of the cylinder to $\infty$. The
 right hand side of
\eqref{modcob} counts sections satisfying the incidence conditions
given by $A_i,D_j$ for any finite $T$; as $T$ becomes large these
sections approach the broken sections (consisting of a section of
$\mathbb{R}\times (Y_-)_d(f_-)$ followed by a section of
$\bar{X}_d(\bar{f})$) counted by the left hand side of
\eqref{modcob}.

The duality statement comprising part (iii) of \fullref{cobmaps}
follows immediately from the definition of $F_{\frak{m},\theta}$:
the quantities \[\langle
F_{\frak{m},\theta,\tilde{h}}(\mathbf{x}_-),\mathbf{x}_+\rangle_{\frak{o}_+}
\mbox{ and } \langle
\mathbf{x}_-,F_{-\frak{m},\theta,\tilde{h}}(\mathbf{x}_+)\rangle\]
count precisely the same objects (namely holomorphic sections of
$\bar{X}_d(\bar{f})$ asymptotic to $\mathbf{x}_{\pm}$ at the
boundary components $(Y_{\pm})_d(f_{\pm})$ in the relative homotopy
class $\gamma_{\tilde{h}}$), and do so with identical weights (since
the relevant forms $\bar{\Omega}$ and $\uparrow\theta$ are the same
for $\frak{m}$ as for $-\frak{m}$); hence these quantities are
equal.

Since part (iv) has already been established, the proof of \fullref{cobmaps} is now complete.

\section{Periodic points of symplectomorphisms and asymptotics for the parallel translation of vortices}\label{sec.6}

Let us first prove \fullref{fixedpoint} assuming \fullref{pfhmain}. Note that for any diffeomorphism $\phi\co
\Sigma\to\Sigma$ there is a Lefschetz fibration $f\co X\to S^2$ with
fiber $\Sigma$ over some point and all fibers irreducible such that
the preimage of the equator $\gamma$ in $S^2$ is isomorphic as a
smooth fibration to the mapping torus $Y_{\phi}$ of $\phi$ (for by
\cite[Theorem 2.2]{OS2} one can factor the mapping class of $\phi$
as a product of right-handed Dehn twists along nonseparating loops
to get a Lefschetz fibration over the disc such that the monodromy
around the boundary is isotopic to $\phi$, and then complete this
factorization to a factorization of the identity as a product of
right-handed Dehn twists along nonseparating loops in order to
complete the Lefschetz fibration to a Lefschetz fibration over the
whole sphere), such that $b^+(X)>1$ (if the initially--constructed
Lefschetz fibration does not satisfy this property, then its fiber
sum with itself will, by, eg, \cite[Lemma 3.1]{stip}).
Now $T^{vt} X$ is a well-defined complex line bundle on the
complement of a set of codimension four (namely $Crit(f)$) in $X$
and so extends from the complement of a neighborhood of $Crit(f)$ to
a complex line bundle on all of $X$, and then $\langle
c_1(T^{vt}X),[\Sigma]\rangle<0$. Hence there are, by the proof of
\cite[Theorem 10.2.18]{Gompf}, symplectic forms $\beta$ on $X$ in
classes of form $-c_1(T^{vt}X)+Mf^*\omega_{S^2}$ for large $M$.  We
wish to say that $DS(PD(\kappa_X);pt^0)\neq 0$ where $\kappa_X$ is
the canonical class; if our Lefschetz fibration were obtained by
blowing up a high degree Lefschetz pencil on a manifold with
$b^+>b_1+1$ we could deduce this directly from the main result of
\cite{DS}.  For more general Lefschetz fibrations, Taubes' theorems
\cite{Tbook} show that, with respect to the symplectic structure
$\beta$, we have $Gr_X(PD(\kappa_X))=\pm 1$, and hence, by the main
theorem of \cite{U},\footnote{In the statement of the main theorem
of \cite{U}, there is a hypothesis on the area of the fiber of the
Lefschetz fibration.  However, for any Lefschetz fibration $f\co
X\to S^2$ with all fibers irreducible, the main theorem of \cite{U}
still applies to show that $Gr(\alpha;\cdot)=DS(\alpha;\cdot)$ for
any class $\alpha\in H_2(X;\mathbb{Z})$ satisfying
$d=\alpha\cap[\Sigma]>g-1$, since then the virtual dimension
of the space of pseudoholomorphic curves representing any class of
form $\alpha-n[\Sigma]$ with $n>0$ will be smaller than the virtual
dimension of pseudoholomorphic representatives of $\alpha$, and so
the Gromov--Taubes moduli spaces for the class $\alpha$ will, for
generic almost complex structures making $f$ pseudoholomorphic, not
contain any curves with fiber components; ensuring that this be the
case was the only role played by the hypothesis on the area of the
fiber in \cite{U}.} $DS(PD(\kappa_X);pt^0)=\pm 1$; meanwhile all
other classes $\alpha\in H_2(X;\mathbb{Z})$ differing from
$PD(\kappa_X)$ by a torsion element have
$Gr_X(\alpha)=DS(\alpha;pt^0)=0$.  So where
$h=PD(\kappa_X|_{Y_{\phi}})\in H_1(Y_{\phi};\mathbb{Z})$ and
$c=\beta|_{Y_{\phi}}=-c_1(T^{vt}Y_{\phi})$, we conclude that the
composition \begin{align*} \smash{\sum_{\begin{array}{c}\scriptstyle
\tilde{h}\in
 H_2(X;\mathbb{Z}):\\\scriptstyle PD(\tilde{h})|_{Y_{\phi}}=\kappa_X|_{Y_{\phi}}\end{array}}}F_{X,\theta,\tilde{h}}(1\otimes\cdot)\co\mathbb{R}=HF(&\varnothing;\mathbb{R})\to
HF(Y_{\phi},f,h,c,\tau;\Gamma_{\theta|_{Y_{\phi}}})\\
&\to\mathbb{R}=HF(\varnothing;\mathbb{R})\end{align*} is nonzero for
certain choices of $\theta\in \Omega^2(X)$ (and for arbitrary
$\tau$), and so\break
$HF(Y_{\phi},f,h,c,\tau;\Gamma_{\theta}|_{Y_{\phi}})\neq 0$. The
monotonicity assumption on $\phi$ implies that $\omega_{\phi}$
(after rescaling) belongs to the class $c+\frac{2\pi}{\tau}PD(h)$,
and so may be used in the definition of $HF(Y_{\phi},f,h,c,\tau)$.
Hence, for all  $\tau>2\pi d$, noting that $\langle
h,\fiber\rangle=2g-2$, the symplectomorphism
$\Phi_{2g-2,\omega_{\phi},\tau}$ has a fixed point.  But as
$\tau\to\infty$, $\Phi_{2g-2,\omega_{\phi},\tau}\to S^{2g-2}\phi$ by
\fullref{pfhmain}, so the latter map has a fixed point as
well.\quad \qed

As was alluded to in the introduction, the same argument  reveals
that for $d>g-1$ $S^d\phi$ has a fixed point whenever $\phi$ is
monotone and there is a Lefschetz fibration $f\co X\to S^2$ with
irreducible fibers having monodromy around some loop isotopic to
$\phi$ whose total space has the property that, for some homology
class $\tilde{h}\in H_2(X;\mathbb{Z})$ having intersection number
$d$ with the fibers, the sum of the Gromov--Taubes invariants in
classes congruent to $\tilde{h}$ mod torsion and mod restriction to
$Y_{\phi}$ is nonzero (with the slight modification that for the
class $c$ of the fiberwise symplectic form one should use
$c=-c_1(T^{vt}X)-\frac{2\pi}{\tau}PD(\tilde{h})$, with $\tau$ large
enough to ensure that this class is positive on the fibers).  As
mentioned in footnote 4, the requirement that $d>g-1$ along with the
irreducibility of the fibers  suffice to replace the assumption on
the symplectic areas of $\tilde{h}$ and $[\Sigma]$ in the main
theorem of \cite{U}. The same reasoning can also be applied to
certain non-monotone symplectomorphisms $\phi$, provided that there
is a Lefschetz fibration containing the mapping torus of $\phi$ as
the preimage of some circle in the base, and carrying a symplectic
form in a cohomology class which restricts appropriately to this
mapping torus.

We turn finally to the proof of \fullref{pfhmain}.  We consider
a symplectomorphism $\phi\co (\Sigma,\omega)\to(\Sigma,\omega)$ of a
symplectic $2$--manifold. $\omega$ induces on the mapping torus
$Y_{\phi}$ a closed fiberwise symplectic form $\omega_{\phi}$ in the
cohomology class $c\in H^2(Y;\mathbb{R})$.  In \fullref{Define}
we have, for each large enough $\tau$, chosen closed fiberwise
symplectic forms $\omega^{\tau}$ on $Y_{\phi}$ representing the
classes $c+2\pi PD(h)/\tau$; since the homology will be independent
of the particular forms in these classes that we choose we may as
well assume that $\omega^{\tau}\to \omega_{\phi}$ as
$\tau\to\infty$, and (using the Moser trick) that the restriction of
$\omega^{\tau}$ to some fixed base fiber is proportional to
$\omega$.  The monodromies $\phi_{\tau}$ of the $\omega^{\tau}$ then
converge to $\phi$, and so from the definition of the chain complex
$CF$ it follows that \fullref{pfhmain} can be translated into
the statement that the parallel transport map $F_{\{J_t\}}\co
S^d(\Sigma,J_0)\to S^d(\Sigma,J_1)$ of \eqref{transport} approaches
the identity as the parameter $\tau$ tends to $\infty$.

In proving this, we shall make use of the asymptotic properties of
the vortices themselves for large $\tau$.  Recall that the vortex
equations are obtained by fixing a degree $d$ Hermitian line bundle
$L$ on the K\"ahler curve $(\Sigma,\omega,J)$; they read
\begin{align}\label{vortextw} \bar{\partial}_{J,A}\theta&=0\nonumber \\ iF_A&=\tau(1-|\theta|^2)\omega.\end{align}
where the unknown $(A,\theta)$ consists of a connection $A$ in $L$
and a not--identically--zero section $\theta$ of $L$.    In the case
that $(\Sigma,\omega,J)$ is $\mathbb{R}^2$ with its standard
symplectic and complex structure, solutions for general $\tau$ can
be obtained from those from the case $\tau=1$ by pulling back via
the dilation $z\mapsto\sqrt{\tau}z$. The case of the standard plane
with $\tau=1$ was exhaustively analyzed in \cite[Chapter III]{JT};
in particular, according to Theorem III.8.5, the curvature satisfies
an exponential decay condition which translates to the general
$\tau$ case as \[ |\ast iF_A|\leq M\tau e^{-c\sqrt{\tau}|z|},\]
where $c$ can be taken to be any constant smaller than $2$.  We
shall be needing analogous (though somewhat weaker) results for the
vortices on general $(\Sigma,\omega,J)$.  The referee has pointed
out that bounds similar to what we prove (at least for a fixed $J$)
can be deduced from estimates on solutions to the Seiberg--Witten
equations from \cite{swgr} (in particular 1.24 (6)) by specializing
to the case where the four--manifold under consideration is
$\Sigma\times T^2$ with a product metric; however, we shall still
give our proof of these bounds because the proof is simpler (though
similar in spirit) in the purely two-dimensional case, because we
need to see explicitly  that the estimates are uniform when we vary
$J$ in a compact 1--parameter family, and because some of the
necessary ingredients will reappear later when we analyze a certain
Green's function. Readers familiar with the proofs of such bounds
might skip \fullref{basic} through \fullref{expldec}.

Throughout our discussion, we work with a fixed Hermitian line
bundle $L\to \Sigma$ of degree $d>0$ over a fixed compact symplectic
$2$--manifold $(\Sigma,\omega)$.  To connect this to the setup in
\fullref{Define}, we should note that in that section the closed
fiberwise symplectic form $\omega^{\tau}$ restricts to $\Sigma$ as
$\frac{\langle c,[\Sigma]\rangle+2\pi d/\tau}{\langle
c,[\Sigma]\rangle}$ times $\omega$. As such, the parameter $\tau$ in
\eqref{vortextw} would be $\tau+\frac{2\pi d}{\langle
c,[\Sigma]\rangle}$ in the notation of \fullref{Define}.  Since
we are interested here in the behavior of the vortex equations as
$\tau\to\infty$ and since $\langle c,[\Sigma]\rangle>0$, this
distinction is immaterial to our present concerns and we shall
henceforth suppress it.

We will also fix a smooth path $\{J_t\}_{t\in [0,1]}$ of almost
complex structures on $\Sigma$; together with $\omega$, these induce
metrics $g_t$. Consider the vortex equations \eqref{vortextw} where
$J$ is one of the almost complex structures appearing in the path
$J_t$. We shall be making a variety of estimates on some quantities
relating to solutions of these equations, which shall involve
certain constants; these constants may be taken independent of
$\tau$ and of $J$ \emph{provided that $J$ is chosen from within the
fixed smooth $1$--parameter family} $\{J_t\}_{t\in [0,1]}$, but
might not apply to an entirely arbitrary choice of $J$.  More
specifically, where $g$ is the metric induced by $\omega$ and $J$,
the constants may depend on any or all of: the minimal or maximal
curvature of the Riemannian $2$--manifold $(\Sigma, g)$; the
injectivity radius $r_0$ of $(\Sigma,g)$; the diameter of
$(\Sigma,g)$; or the maximum of the Jacobians of the exponential
maps $\exp_{p}^{g}\co B_{r_0/2}(0)\to \Sigma$ for $p\in\Sigma$.

First, we prove a direct analogue for the case of a general Riemann
surface to a pair of properties proven for the case of the flat
plane in \cite{JT}.  Let \[
\kappa=\max_{t\in[0,1],p\in\Sigma}\{0,-\sec_{g_t}(p)\},\] where
$\sec_{g_t}(p)$ is the sectional curvature of $\Sigma$ at $p$ in the
metric $g_t$.

\begin{lemma} \label{basic} Any solution $(A,\theta)$ to \eqref{vortextw} satisfies:
\[ w:=1-|\theta|^2\geq 0  \qquad  |d_A\theta|\leq
 2\tau^{1/2}w+2\kappa\tau^{-1/2},\] provided that $\tau\geq \kappa$.
\end{lemma}

\begin{proof}
First note that, for any section $\phi$ of a holomorphic line bundle
$V$ with unitary connection $A$ over any K\"ahler manifold $M$, one
has
\begin{align*}
(d_A\delbar_A\phi)(v,w)&=\nabla_v\iota_w\delbar_A\phi-\nabla_w\iota_v\delbar_A\phi-(\delbar_A\phi)([v,w])\\
&=\frac{1}{2}F_A\phi(v,w)+\frac{i}{2}(\nabla^{2}_{v,iw}\phi-\nabla^{2}_{w,iv}\phi),\end{align*}
as can be seen by expanding out
$\iota_u\overline{\partial}_A\phi=\nabla_u\phi+i\nabla_{iu}\phi$ and
then using the fact that $M$ is K\"ahler to move various factors of
$i$ past covariant derivatives.

In particular \[
d_A\delbar_A\phi(v,iv)=\frac{1}{2}F_A\phi(v,iv)-\frac{i}{2}(\nabla^{2}_{v,v}\phi+\nabla^{2}_{iv,iv}\phi),\]
so that if $M$ is $1$--complex dimensional we see that \[ \ast
d_A\delbar_A\phi=\frac{1}{2}\ast F_A\phi-\frac{i}{2}\Delta\phi,\]
where $\ast$ is the Hodge star operator induced by the metric and
$\Delta=\ast d_A\ast d_A$ is the (negative) Laplacian on sections of
$V$ induced by $A$. Applying this to our vortex $(A,\theta)$ on
$\Sigma$, since $\delbar_A\theta=0$ and $\ast
iF_A=\tau(1-|\theta|^2)$, we see \begin{equation} \label{Won}
\Delta\theta+\tau(1-|\theta|^2)\theta=0.\end{equation} Now \[
\Delta|\theta|^2=2Re\langle\Delta\theta,\theta\rangle+2|d_A\theta|^2,\]
while \eqref{Won} implies that $\langle\Delta\theta,\theta\rangle$ is
real, so that \[ \langle
\Delta\theta,\theta\rangle=\Delta|\theta|^2/2-|d_A\theta|^2.\] Hence
taking the inner product of \eqref{Won} with $\theta$ and setting
$w=1-|\theta|^2$ yields \begin{equation}\label{lapw} -\Delta
w+2\tau|\theta|^2w=2|d_A\theta|^2.\end{equation}  In particular if
$z_0\in\Sigma$ were such that $w(z_0)<0$, we would have $(-\Delta
w)(z_0)>0$, so that $z_0$ could not be a local minimum for $w$.  So
since $\Sigma$ is compact and $w\co \Sigma\to\mathbb{R}$ cannot
attain a negative local minimum, we have $w\geq 0$ everywhere.
\vspace{-3pt}

Now  $\Delta\theta=\ast d_A\ast d_A\theta=-\ast iF_A\theta=-\tau
w\theta$, where we have used that, since $d_A\theta$ has type
$(1,0)$, $\ast d_A\theta=-id_A\theta$.  So setting $h=d_A\theta\in
\Omega^{1,0}(L)$, we see \begin{align*}d_A\ast d_A\ast h&=d_A(-\tau
w\theta)=-\tau(wh+\theta
dw)=-\tau(wh-\theta(\bar{\theta}h+\theta\bar{h}))\\&=\tau
h(1-2w)+\tau\theta\theta\bar{h},\end{align*} while
\vspace{-3pt}
\begin{align*}
\ast d_A\ast d_A h&=\ast d_A\ast d_Ad_A\theta=-i\tau\ast
d_A(w\theta)\\&=-i\tau\ast(\theta(-\bar{\theta}h-\theta\bar{h})+wh)=-i\tau\ast((2w-1)h-\theta\theta\bar{h})\\
&=\tau (1-2w)h-\tau\theta\theta\bar{h},\end{align*}
where in the
last equality we have used that since $h$ has type $(1,0)$, $\ast
h=-ih$ and $\ast \bar{h}=i\bar{h}$.   So following \cite[Section III.6]{JT} 
by writing $\Delta_A=\ast d_A\ast d_A+d_A\ast d_A \ast$,
we see \[ \Delta_A h=2\tau h(1-2w).\]  Further, on $L$--valued
$1$--forms there is a Weitzenb\"ock formula (see, eg,
\cite[III.6.15]{JT}; \cite[Chapter 7]{Petersen})
\[ tr\nabla_{A}^{2}=\Delta_A+(\ast F_A)\ast +\sec,\] so since $h$
has type $(1,0)$ and so $(\ast F_A)\ast h=(-\ast iF_A)h=-\tau w h$,
we obtain
\[ tr\nabla^{2}_{A}h=\tau(2-5w+\sec/\tau) h,\] from which Kato's inequality
\cite[III.6.20]{JT} provides \[ |h|\Delta |h|\geq \tau
(2-\kappa/\tau-5w)|h|^2.\]  Hence \begin{align*}
|h|\Delta&(2\sqrt{\tau}(w+\kappa/\tau)-|h|)\leq
2\sqrt{\tau}|h|\left(2\tau
w(1-w)-2|h|^2\right)+(5w+\kappa/\tau-2)\tau|h|^2\\
&=\sqrt{\tau}|h|\left(4\tau(1-w)w+\sqrt{\tau}|h|(5w+\kappa/\tau-2)-4|h|^2\right)\\&=\sqrt{\tau}|h|\left((2\sqrt{\tau}(w+\kappa/\tau)-|h|)\left(2\sqrt{\tau}\frac{w(1-w)}{w+\kappa/\tau}+4|h|\right)\right)\\&\quad -\tau|h|^2\left((2-\kappa/\tau-5w)-\left(\frac{2w(1-w)}{w+\kappa/\tau}-8(w+\kappa/\tau)\right)\right)\\
&\leq
\sqrt{\tau}|h|\left((2\sqrt{\tau}(w+\kappa/\tau)-|h|)\left(2\sqrt{\tau}\frac{w(1-w)}{w+\kappa/\tau}+4|h|\right)\right),
\end{align*} where we have used the fact that, since $0\leq w\leq
1$ and $\kappa/\tau\leq 1$, we have
$\frac{2w}{w+\kappa/\tau}\leq\frac{2}{1+\kappa/\tau}\leq
2-\kappa/\tau$.

So  we see that wherever $2\tau^{1/2}w+2\kappa\tau^{-1/2}-|h|$ is
negative (which forces $|h|>0$ since we've already shown that $w\geq
0$ everywhere) we have
$\Delta(2\tau^{1/2}w+2\kappa\tau^{-1/2}-|h|)<0$. But then
$2\tau^{1/2}w+2\kappa\tau^{-1/2}-|h|\co\Sigma\to\mathbb{R}$ cannot
attain a negative local minimum, which by the compactness of
$\Sigma$ forces $|d_A\theta|=|h|\leq
2\tau^{1/2}w+2\kappa\tau^{-1/2}$ everywhere.\end{proof}

\begin{prop} \label{firstdecay}There is a constant $C>0$ with the property that, for all sufficiently large $\tau$ if
$(A,\theta)$ is a solution to \eqref{vortextw} with
$J\in\{J_t\}_{t\in[0,1]}$ and $w=1-|\theta|^2$, we have, for each
$z_0\in \Sigma$,
\[ w(z_0)\min\{d(z_0,p)|\theta(p)=0\}\leq \frac{C}{\sqrt{\tau}},\]
where $d(\cdot,\cdot)$ denotes the distance measured in the metric
$g$ induced by $J$ and $\omega$.\end{prop}

\begin{proof}  First note that the first statement of \fullref{basic} shows that $|\theta|\leq 1$, so
$|d(w+\kappa/\tau)|=|2Re\langle\theta,d_A\theta\rangle|\leq
2|d_A\theta|\leq 4\sqrt{\tau}(w+\kappa/\tau)$.  So if $\gamma$ is an
arc--length parametrized path in $\Sigma$, say from $z$ to $z'$ and
having length $l_{\gamma}$, we have \begin{align*}
\log\left(\frac{w(z')+\kappa/\tau}{w(z)+\kappa/\tau}\right)&=\int_{0}^{l_{\gamma}}\frac{d}{dt}\log(w(\gamma(t))+\kappa/\tau)\,dt\leq\int_{0}^{l_{\gamma}}\frac{|d(w+\kappa/\tau)|}{|w+\kappa/\tau|}\,dt\\&\leq
4\sqrt{\tau}l_{\gamma}.\end{align*}  Thus, for any $z,z'\in\Sigma$,
\[ w(z')+\kappa/\tau \geq
(w(z)+\kappa/\tau)e^{-4\sqrt{\tau}d(z,z')}.\]

We claim now that $w=1-|\theta|^2$ is equal to either $0$ or $1$ at
each of its local maxima. Indeed, note that, again writing
$h=d_A\theta\in \Omega^{1,0}(L)$, we have
$dw=-d|\theta|^2=-(\bar{\theta}h+\theta\bar{h})$, and $\ast
dw=-\ast(\bar{\theta}h+\theta\bar{h})=i\bar{\theta}h-i\theta\bar{h}$,
so that \[ dw+i\ast dw=-2\bar{\theta}h, \] so at a putative local
maximum $z$ of $w$ with $w(z)\notin \{0,1\}$ (so $\theta(z)\neq 0$),
we necessarily have $h(z)=0$.  Meanwhile since $w$ takes values only
in $[0,1]$ we must also have $w(z)(1-w(z))>0$, so
 recalling
the equation \[ -\Delta w+2\tau w(1-w)=|h|^2 \] we see that $\Delta
w(z)>0$, in contradiction with the fact that $z$ was taken to be a
local maximum.

Suppose now that \[ N\tau^{-1/2}\geq
\min\{d(z_0,p)|\theta(p)=0\}\geq (N-1)\tau^{-1/2}\quad
(N\in\mathbb{N}).\]  Since $w$ is everywhere nonnegative and is
strictly less than $1$ on $B_{(N-1)\tau^{-1/2}}(z_0)$, we deduce
that for $k=1,\ldots,N-1,$ $\sup_{B_{k/\sqrt{\tau}}(z_0)}w$ must be
attained at some point $z_k$ with $d(z_k,z_0)=k\tau^{-1/2}$; in
particular $w(z_k)\geq w(z_0)$. This together with the conclusion of
the first paragraph of the proof shows that, on each of the disjoint
balls $B_k=B_{\frac{1}{2\sqrt{\tau}}}(z_k)$ ($k=0,\ldots,N-1$),
\[ w|_{B_k}\geq e^{-2}w(z_0)-\kappa/\tau.\]  Now for some constant $A$ (related to the Jacobian of the exponential map on balls of radius smaller than the injectivity radius, if $\tau$ is large enough) we have
$\vol(B_k)\geq A\tau^{-1}$ for each $k$, and so since $w\geq 0$
throughout $\Sigma$
\[ \int_{\Sigma}w\omega\geq \sum_{k=0}^{N-1}\int_{B_k}w\omega\geq
ANe^{-2}w(z_0)\tau^{-1}-A\kappa\tau^{-2}.\] But the original vortex
equations imply that \[
\int_{\Sigma}w\omega=\tau^{-1}\int_{\Sigma}\ast iF_A=2\pi
d\tau^{-1}.\]  Thus \[ w(z_0)\min\{d(z_0,p)|\theta(p)=0\}\leq
w(z_0)N\tau^{-1/2}\leq C\tau^{-1/2}\] for an appropriate choice of
$C$.
\end{proof}

Having taken these first steps, we can now prove a basic exponential
decay estimate for vortices.

\begin{theorem}\label{expldec} There are constants $R,M,\tau_0>0$ with the property that, if
$(A,\theta)$ is a solution to \eqref{vortextw} with
$J\in\{J_t\}_{t\in[0,1]}$ and $w=1-|\theta|^2$, we have, for each
$z\in \Sigma$ and for $\tau\geq \tau_0$, \[
w(z):=1-|\theta(z)|^2\leq
\kappa\tau^{-1}+M\sum_{\{p:\theta(p)=0\}}\left(e^{-\sqrt{\tau}d(z,p)}+e^{-R\sqrt{\tau}}\right).\]
\end{theorem}

\begin{proof}
According to \fullref{basic} and Equation \ref{lapw}, we have \[
-\Delta w+2\tau w(1-w)=2|d_A\theta|^2\leq 8\tau(w+\kappa/\tau)^2,\]
so that
\[ -\Delta w\leq\tau w(34w-2) \mbox{ wherever $w>\kappa/\tau$}.\]  Meanwhile, according to
\fullref{firstdecay}, where \[ V=\{z\in\Sigma:
d(z,\theta^{-1}(0))\geq 68 C\tau^{-1/2}\},\] at each $z\in V$ we
have $w(z)\leq 1/68$, and so \[ -\Delta w\leq -\frac{3}{2}\tau w
\mbox{ on $V\cap\{w>\kappa/\tau\}$}.\]

Now note that if $q\in \Sigma$ and $\alpha\co \mathbb{R}_{\geq
0}\to\mathbb{R}$, the function
$u_{q,\alpha}(z)=e^{-\sqrt{\tau}\alpha(d(q,z))}$ satisfies (wherever
all terms exist) \[ \Delta
u_{q,\alpha}(z)=\left(\tau\alpha'(d(q,z))^2-\sqrt{\tau}(\alpha''(d(q,z))-\sqrt{\tau}\alpha'(d(q,z))\Delta(d(q,z))\right)u_{q,\alpha}(z).\]

Also, if the curvature of $(\Sigma,g)$ is bounded above by $K>0$,
then it follows from Theorem 6.2.1 and the discussion before Lemma
9.1.1 in \cite{Petersen} that where $R$ is the minimum of
$\frac{\pi}{2\sqrt{K}}$ and the injectivity radius of $(\Sigma,g)$,
we have $\Delta(d(q,z))\geq 0$ as long as $d(q,z)\leq R$.

Take for $\alpha$ a smooth function with the following properties:
\begin{itemize} \item[(i)] $\alpha(0)=0$,\item[(ii)]
$0\leq\alpha'(t)\leq 6/5$, with $\alpha'(t)=6/5$ for
$t<R/2$,\item[(iii)]$\alpha(t)=R$ for $t\geq R$, \item[(iv)]
$-3/R\leq \alpha''(t)\leq 0$.\end{itemize}  Then since
$u_{q,\alpha}$ is positive everywhere and is constant outside the
region that $\Delta d(q,\cdot)$ is known to be nonnegative, we have
\begin{align*} \Delta u_{q,\alpha}&\leq
\left(\frac{36}{25}\tau+\frac{3\sqrt{\tau}}{R}\right)u_{q,\alpha}\\[1ex]
&\leq \frac{3}{2}\tau u_{q,\alpha}\end{align*} provided that
$\sqrt{\tau}\geq 50/R$.
\vspace{3pt}

Also, assuming that $\tau$ is large enough that $68C\tau^{-1/2}\leq
R/2$, if $d(q,z)=68C\tau^{-1/2}$ then $u_{q,\alpha}(z)=e^{-408C/5}$.
Hence setting $M=e^{408C/5}$ and \[
u=M\sum_{p\in\theta^{-1}(0)}u_{p,\alpha},\] \[ u|_{\partial
V}\geq 1> \frac{1}{68}\geq w|_{\partial V}\tag*{\hbox{we have}} \] \[ \Delta
(u-w)\leq\frac{3}{2}\tau (u-w)\mbox{ on $V\cap\{w>\kappa/\tau\}$}.
\tag*{\hbox{and}} \]

But then, as usual, $\Delta (u-w)<0$ anywhere on
$V\cap\{w>\kappa/\tau\}$ that $u-w$ is negative, so that $u-w$
cannot attain a negative local minimum on $V\cap\{w>\kappa/\tau\}$.
In particular, then, $u+\kappa/\tau-w$ also cannot attain a negative
local minimum on $V\cap\{w>\kappa/\tau\}$, so since
$u+\kappa/\tau-w$ is obviously positive where $w\leq\kappa/\tau$,
$u+\kappa/\tau-w$ cannot attain a negative local minimum anywhere on
$V$. So since $u-w>0$ on $\partial V$ we deduce that $w\leq
u+\kappa/\tau$ througout $V$, and indeed throughout $\Sigma$ since
away from $V$ we have $u\geq 1\geq w$. The proof is then completed
by noting that the construction of $\alpha$ ensures that, for each
$p\in \theta^{-1}(0)$, we have $$ u_{p,\alpha}(z)\leq
e^{-\sqrt{\tau}d(p,z)}+e^{-R\sqrt{\tau}}.\proved$$
\end{proof}
\vspace{3pt}

Now the parallel translation inducing the map $F_{\{J_t\}}\co
S^d(\Sigma,J_0)\to S^d(\Sigma,J_1)$ that we are investigating is, by
\cite[Theorem 5.1]{Sal}, given by \[
F_{\{J_t\}}([A,\theta])=[A(1),\theta(1)]\] where $(A(t),\theta(t))$
solves the ODE
\[ i\dot{A}(t)=2\tau Re\langle\theta(t),\eta(t)\rangle\qquad
i\dot{\theta}(t)=\delbar_{A(t)}^{*}\eta(t)
\] with initial condition $(A(0),\theta(0))=(A,\theta)$ where $\eta(t)\in \Omega^{0,1}(L)$ is the unique solution to
\begin{equation} \label{trans}
\delbar_{J_t,A(t)}\delbar_{J_t,A(t)}^{*}\eta(t)+\tau|\theta(t)|^2\eta=\frac{1}{2}(d_{A(t)}\theta(t))\circ
\dot{J}(t).\end{equation}

(Recall that the isomorphism between the set of gauge equivalence
classes of vortices takes $[A,\theta]$ to the vanishing locus of
$\theta$.  Also, to compare to \cite{Sal}, our terms $\theta$ and
$\eta$ are $1/\sqrt{2\tau}$ times the corresponding terms
$\Theta_0,\Theta_1$, respectively, in \cite{Sal}.  The reader may
calculate directly or consult the proof of \cite[Theorem 5.1]{Sal}
to see that $(A(t),\theta(t))$ so defined does indeed satisfy
\eqref{vortextw} with $J=J_t$ for all $t$ and that this recipe is
consistent with the symplectic parallel transport description
discussed in Section 2.)
\vspace{3pt}

Our goal is to show that $F_{\{J_t\}}$ is close to the identity; we
shall accomplish this by obtaining upper bounds on
$|\delbar_{J_t,A(t)}^{*}\eta(t)|$ where $\eta(t)$ solves
\eqref{trans}. Now where \[G(x,p)\co \Omega^{0,1}(L)|_p\to
\Omega^{0,1}(L)|_x\] denotes the Green's kernel for the operator
\[ \delbar_{J,A}\delbar_{J,A}^{*}+\tau|\theta|^2\co
\Omega^{0,1}(L)\to\Omega^{0,1}(L)\] we have \[
\eta(t)(x)=\frac{1}{2}\int_{\Sigma}G(x,p)\left((d_{A(t)}\theta(t,p))\circ
\dot{J}(t,p)\right)\omega_p;\] the desired upper bounds on
$|\delbar_{J_t,A(t)}^{*}\eta(t)|$  will follow from our
already--obtained exponential decay bounds on $w(t)=1-|\theta(t)|^2$
(and hence on $|d_{A(t)}\theta(t)|$ by \fullref{basic}), together
with bounds on the derivatives of the Green's kernel.

We now set about deriving these Green's kernel estimates.  Let
$(A,\theta)$ be an arbitrary solution to \eqref{vortextw} (with $J$
taken from the path $\{J_t\}_{t\in[0,1]}$; with this $J$ understood,
we shall just write $\delbar_A$ for $\delbar_{J,A}$). Note that
$\delbar_{A}\delbar_{A}^{*}+\tau|\theta|^2$ is manifestly positive
definite, and in fact the Weitzenb\"ock formula used in the proof of
\fullref{ineq} below allows us to rewrite this operator as
$\frac{1}{2}(\nabla_A^{*}\nabla_A+\tau(1+|\theta|^2)+\sec)$, and so
as long as $\tau>4\kappa$ (as we shall assume hereinafter) its
spectrum is bounded below by $\tau/4$.

We first obtain estimates on $G(x,p)$ for $p$ close to $x$.  In this
direction, consider the effect of replacing the metric $g$ induced
by $J$ and $\omega$ by $\tilde{g}=\tau g$.  Then, since on
$1$--forms we have
$\delbar_{A}^{*_{\tilde{g}}}=-\ast_{\tilde{g}}\partial_A\ast_{\tilde{g}}=\tau^{-1}\delbar_{A}^{*_{g}}$
we see that $G$ is also the Green's kernel (using the metric
$\tilde{g}$) for the operator \[
\delbar_A\delbar_{A}^{*_{\tilde{g}}}+|\theta|^2\co\Omega^{0,1}(L)\to\Omega^{0,1}(L),\]
where by \fullref{basic} we have $\tau$--independent bounds
\[0\leq|\theta|^2\leq 1, \quad |d|\theta|^2|_{\tilde{g}}\leq 3/2\] for the
potential term $|\theta|^2$.  Furthermore \eqref{lapw} gives \[
\Delta^{\tilde{g}}|\theta|^2+2|\theta|^2(1-|\theta|^2)=|d_A\theta|_{\tilde{g}}^{2};\]
differentiating this and repeatedly using the bounds of \fullref{basic} and the fact that $\ast iF_A=\tau (1-|\theta|^2)$ gives,
for all $k$, $\tau$--independent constants $C_k$ such that
\[ |(\Delta^{\tilde{g}})^k |\theta|^2|\leq C_k, \quad
|d(\Delta^{\tilde{g}})^k|\theta|^2|_{\tilde{g}}\leq C_k.\]  Using
the approach of  \cite[Chapitre III, E.III]{BGM} (adapted from the
case of the Laplacian on functions to that of a more general Laplace
type operator on sections of a vector bundle as in \cite{Gilkey})
one then finds, for a fixed $c<\inf\,injrad(\Sigma,\tilde{g})$,
uniform--in--$\tau$ estimates on the $c$--neighborhood of the
diagonal in $(\Sigma,\tilde{g})\times(\Sigma,\tilde{g})$ for the
$C^1$--accuracy of the third-order asymptotic approximation
$S_3(t,x,y)$ to the heat kernel $S(t,x,y)$ for
$\delbar_A\delbar_{A}^{*_{\tilde{g}}}+|\theta|^2$.  (Note that since
the functions $K_k$ of \cite[Chapitre III, Lemme E.III.6]{BGM}
vanish outside the $c$--neighborhood of the diagonal the term $V$ on
page 212 can be replaced by the maximal volume of a ball of radius $c$
in the $(\Sigma,\tilde{g})$, which is bounded independently of
$\tau$.)  Since the spectrum of
$\delbar_A\delbar_{A}^{*_{\tilde{g}}}+|\theta|^2$ is bounded below
by $1/4$, we may then integrate with respect to $t$ to see that
these estimates imply a uniform bound on the $C^1$--norm of the
difference between a cut-off version of the third-order Hadamard
expansion of the Green's kernel and the actual kernel $G$ (see,
eg, \cite[section II.2]{Av}; the relevant coefficients
in the expansion may be found on page 336 of \cite{Gilkey}, and the
salient point for our purposes is that the $k$th derivatives of the
potential term $|\theta|^2$ only contribute a correction factor
proportional to the $2(k+1)$th power of the distance and so do not
substantially affect the rate at which $G(x,p)$ diverges near the
diagonal). As a result, there is a $\tau$--independent constant
$C>0$ such that, whenever $dist_{\tilde{g}}(x,p)\leq c$, we have
\[ |G(x,p)|\leq C(1+|\log d_{\tilde{g}}(x,p)|),\qquad
|\delbar_{A}^{*_{\tilde{g}}}G(x,p)|_{\tilde{g}}\leq
\frac{C}{d_{\tilde{g}}(x,p)},\] where in the second formula we are
viewing $p$ as fixed, so that $x\mapsto G(x,p)$ is an element of
$(\Omega^{0,1}(L)|_{p})^{*}\otimes\Omega^{0,1}(L)$, and then taking
$\delbar_{A}^{*_{\tilde{g}}}$ of this section (with respect to $x$).
Scaling back, these relations translate to:
\begin{align}\label{localgreen} |G(x,p)|\leq C(1+|\log
\tau^{1/2}d_g(x,p)|)&,\quad
|\delbar_{A}^{*_g}G(x,p)|_g\leq\frac{C}{d_g(x,p)}\quad
\\\mbox{whenever }d_g(x,p)&\leq c\tau^{-1/2}\nonumber\end{align}
since the $\tilde{g}=\tau g$--norm of a given element of
$\Hom(\Omega^{0,1}(L)|_p,L|_x)$ is $\tau^{1/2}$ times its $g$--norm.
Hereinafter we use $g$ to measure all distances and norms and to
take all adjoints, and so we shall drop $g$ from notations such as
$\delbar_{A}^{*_g}$.

\begin{lemma}\label{ineq} Fix $p\in\Sigma$, view $G(p,\cdot)$ as a section of
$\Hom(T^{0,1}_{p}\Sigma\otimes L|_p,\Lambda^{0,1}\Sigma\otimes L)$,
and write
\[ \beta=\delbar_{A}^*G(p,\cdot)\in \Gamma(\Hom(T^{0,1}_{p}\Sigma\otimes
L|_p,L)).\]  Then provided that $\tau\geq\kappa/8$ we have the
differential inequality
\[ \Delta\left(4|\beta|^2+81\tau|G|^2\right)\geq
\frac{3}{2}\tau\left(4|\beta|^2+81\tau|G|^2\right) \mbox{ on
}\Sigma\setminus\{p\}.\]
\end{lemma}

\begin{proof}
First, on $\Sigma\setminus\{p\}$, we have by the definition of $G$
\begin{equation}\label{Gzero} \delbar_A\delbar_{A}^{*}G+\tau|\theta|^2G=0.\end{equation}  Now
\begin{align*}
\langle\delbar_A\delbar_{A}^{*}G,G\rangle&=\langle-\frac{1}{2}\Delta
G,G\rangle=-\frac{1}{2}\left(\langle
tr\nabla_{A}^{2}G,G\rangle-(\ast
iF_A+\sec)|G|^2\right)\\&=-\frac{1}{4}\Delta|G|^2+\frac{1}{2}|\nabla_AG|^2+\frac{\tau}{2}(1+\sec/\tau-|\theta|^2)|G|^2,\end{align*}
where in the second equality we have used the Weitzenb\"ock formula
on $L$--valued $1$--forms $\Delta=tr\nabla_{A}^{2}-(\ast
F_A)\ast-\sec$ and the fact that $\ast G=iG$ since $G$ has type
$(0,1)$.  Hence taking the inner product of \eqref{Gzero} with $G$
gives
\begin{equation} \label{gsquared}
-\frac{1}{4}\Delta|G|^2+\frac{\tau}{2}(1+\sec/\tau+|\theta|^2)|G|^2+\frac{1}{2}|\nabla_A
G|^2=0\end{equation}

Meanwhile, applying $\delbar_{A}^{*}$ to \eqref{Gzero} gives
\begin{equation}\label{firstb}
\delbar_{A}^{*}\delbar_A\beta+\tau|\theta|^2\beta=i{\tau}\ast(\bar{\theta}
\partial_A\theta\wedge G)\end{equation} (we have used that $\partial_A\bar{\theta}=0$ here).  Now since $\Delta=tr\nabla_{A}^{2}$ on
sections of $L$, we have \[
Re\langle\delbar_{A}^{*}\delbar_A\beta,\beta\rangle=-\frac{1}{2}Re\langle
tr\nabla_{A}^{2}\beta,\beta\rangle=-\frac{1}{4}\Delta|\beta|^2+\frac{1}{2}|\nabla_A\beta|^2,\]
while \[ |\ast(\bar{\theta}\partial_A\theta\wedge G)|\leq
2\tau^{1/2}|\theta|(9/8-|\theta|^2)|G| \] by \fullref{basic} and
the assumption $\kappa\tau^{-1}\leq 1/8,$ so taking the real part of
the inner product of \eqref{firstb} with $\beta$ shows
\begin{equation}\label{bsquared}
-\frac{1}{4}\Delta|\beta|^2+\tau|\theta|^2\beta+\frac{1}{2}|\nabla_A\beta|^2\leq
2\tau^{3/2}|\theta|(9/8-|\theta|^2)|G||\beta|.\end{equation} But now
note that $|\nabla_A G|^2\geq |\delbar_{A}^{*}G|^2=|\beta|^2$, while
$|\nabla_A\beta|^2\geq
|\delbar_A\beta|^2=|\delbar_A\delbar_{A}^{*}G|^2=\tau^2|\theta|^4|G|^2$
by \eqref{Gzero}.  Substituting these relations into
\eqref{gsquared} and \eqref{bsquared} and using that, by assumption,
$\sec\geq -\kappa\geq -\tau/8$, yields
\begin{align}\label{bothsquared} \frac{1}{4}\Delta|G|^2&\geq
\frac{\tau}{2}(\frac{7}{8}+|\theta|^2)|G|^2+\frac{1}{2}|\beta|^2 \\
\frac{1}{4}\Delta|\beta|^2&\geq
\frac{\tau^2}{2}|\theta|^4|G|^2+\tau|\theta|^2|\beta|^2-\frac{9}{4}\tau^{3/2}|\theta||G||\beta|
\nonumber
\end{align}  But since \[ \frac{81}{16}\tau^2|G|^2+4\tau|\theta|^2|\beta|^2\geq
9\tau^{3/2}|\theta||G||\beta|,\] adding $81\tau$ times the first
inequality of \eqref{bothsquared} to $4$ times the second gives
\begin{align*}
\Delta\left(\frac{81}{4}\tau|G|^2+|\beta|^2\right)&\geq
\tau^2\left(\frac{243}{8}+\frac{81}{2}|\theta|^2+2|\theta|^4\right)|G|^2+\frac{81}{2}\tau|\beta|^2
\\&\geq \frac{3}{2}\tau\left(\frac{81}{4}\tau|G|^2+|\beta|^2\right),
\end{align*} as desired.
\end{proof}

\begin{cor} Where $R$ is the constant from \fullref{expldec}, there is a constant $K$ such that, for all sufficiently large $\tau$ and for all $p,x\in \Sigma$ with $d(p,x)\geq
c\tau^{-1/2}$ we have
\[
4|\beta(x,p)|^2+81\tau|G(x,p)|^2\leq
K(e^{-\sqrt{\tau}d(x,p)}+e^{-R\sqrt{\tau}})\tau\log\tau.\]\end{cor}
\begin{proof}    Let $u(x)=e^{-\sqrt{\tau}\alpha(d(p,x))}$
where $\alpha\co \mathbb{R}_{\geq 0}\to\mathbb{R}$ is the same
function as in the proof of \fullref{expldec} (so that in
particular, we have $\Delta u\leq \frac{3}{2}\tau u$ and $u(x)\leq
e^{-\sqrt{\tau}d(p,x)}+e^{-R\sqrt{\tau}}$).  Now if
$d(p,x)=c\tau^{-1/2}$, the local estimates \eqref{localgreen} show
that \[ 4|\beta(x,p)|^2+81\tau|G(x,p)|^2\leq A\tau\log\tau\] for an
appropriate constant $A$, so we can choose $K$ independently of
$\tau$, $p$ such that, when $d(p,x)=c\tau^{-1/2}$,
$Ku(x)\tau\log\tau=Ke^{-6c/5}\tau\log\tau\geq
4|\beta(x,p)|^2+81\tau|G(x,p)|^2$.  So since \[
\Delta\left(Ku\tau\log\tau-\left(4|\beta|^2+81\tau|G|^2\right)\right)\leq
\frac{3}{2}\tau\left(Ku\tau\log\tau-\left(
4|\beta|^2+81\tau|G|^2\right)\right) \] on $\Sigma\setminus
B_{c\tau^{-1/2}}(p)$ and \[Ku\tau\log\tau|_{\partial
B_{c\tau^{-1/2}}(p)}\geq
\left.4|\beta|^2+81\tau|G|^2\right|_{\partial
B_{c\tau^{-1/2}}(p)},\] we deduce by the usual argument that
$Ku\tau\log\tau-\left(4|\beta|^2+81\tau|G|^2\right)$ cannot attain a
negative local minimum and hence must be nonnegative throughout
$\Sigma\setminus B_{c\tau^{-1/2}}(p)$.
\end{proof}

In particular, after renaming $K$, we have \[ |\beta(x,p)|\leq
K(e^{-\sqrt{\tau}d(x,p)/2}+e^{-R\sqrt{\tau}/2})\tau^{1/2}\log\tau\]
when $d(x,p)\geq c\tau^{-1/2}$.
\vspace{3pt}

In light of this corollary, together with \fullref{expldec}, the
local Green's kernel bound \eqref{localgreen} and the parallel
transport prescription \eqref{trans}, we can get bounds on \[
|\dot{\theta}(x)|\leq \left(\sup\frac{|\dot{J}|}{2}\right)
\int_{\Sigma}|\beta(x,y)||d_A\theta(y)|\omega_y\] from simple bounds
on integrals over $\Sigma$ of various expressions involving
functions of form $e^{-a\sqrt{\tau}d(q,\cdot)}$ for $a$ a constant
and $q\in\Sigma$.  Namely, note first that if $r<injrad(\Sigma)$ and
$C>0$ are constants such that for each $z\in\Sigma$, $\exp_z\co
\{(x,y)\in\mathbb{R}^2|x^2+y^2<r^2\}\to \Sigma$ is a diffeomorphism
onto its image with Jacobian at most $C$, then we have, for any
$z\in\Sigma$, \begin{align*}
\int_{\Sigma}e^{-a\sqrt{\tau}d(z,x)}\omega_x &\leq
\vol(\Sigma)e^{-a\sqrt{\tau}r}+C\int_{0}^{2\pi}\int_{0}^{r}e^{-a\sqrt{\tau}\rho}\rho
d\rho d\theta \\[2ex]
&\leq \vol(\Sigma)e^{-a\sqrt{\tau}r}+\frac{2\pi
C}{a^2\tau}.\end{align*}

Along the same lines, if $x,z\in\Sigma$ are two given points, for
any $y\in\Sigma$, adding the equations $d(z,y)+d(x,z)\geq d(x,y)$,
and $2(d(x,y)+d(z,y))\geq 2 d(x,z)$ shows that \[
3(d(x,y)+d(z,y))\geq 2 d(x,y)+d(x,z),\] and so
\begin{align}\label{conex} \int_{\Sigma}
e^{-a\sqrt{\tau}d(x,y)}e^{-a\sqrt{\tau}d(z,y)}\omega_y&\leq
e^{-a\sqrt{\tau}d(x,z)/3}\int_\Sigma
e^{-2a\sqrt{\tau}d(x,y)/3}\omega_y\nonumber\\[2ex]
&\leq
e^{-a\sqrt{\tau}d(x,z)/3}(\frac{C'}{a^2\tau}+Be^{-ar\sqrt{\tau}})\end{align}
for certain constants $B$, $C'$.
\vspace{3pt}

Finally, noting that given $x\in\Sigma$, where $C$ is the same
Jacobian bound as earlier, and we assume that $\tau$ is large enough
that $c\tau^{-1/2}<r$, we have
\[ \int_{B_{c\tau^{-1/2}}(x)}\frac{1}{d(x,y)}\omega_y\leq
C\int_{0}^{2\pi}\int_{0}^{c\tau^{-1/2}}\frac{\rho d\rho
d\theta}{\rho}=\frac{2\pi C c}{\sqrt{\tau}}.\]

So if $d(x,z)\geq
2c\tau^{-1/2}$, so that $d(y,z)\geq 2 d(x,z)$ for each $y\in
B_{c\tau^{-1/2}}(x)$, we get \[
\int_{B_{c\tau^{-1/2}}(x)}\frac{e^{-\sqrt{\tau}d(y,z)}}{d(x,y)}\omega_y\leq
D\tau^{-1/2}e^{-\sqrt{\tau}d(x,z)/2} \] for some constant $D$ while
if $d(x,z)\leq 2c\tau^{-1/2}$ then $e^{-\sqrt{\tau}d(x,z)/2}\geq
e^{-c}$, so that \[
\int_{B_{c\tau^{-1/2}}(x)}\frac{e^{-\sqrt{\tau}d(y,z)}}{d(x,y)}\omega_y\leq
D\tau^{-1/2}e^{-\sqrt{\tau}d(x,z)/2} \] still holds, possibly after
increasing the (still $x,z,$ and $\tau$--independent) constant $D$.

So recalling our estimates \[ |\beta(x,y)|\leq \frac{C}{d(x,y)}
\mbox{ when } d(x,y)\leq c\tau^{-1/2}\]\[ |\beta(x,y)|\leq
K(e^{-\sqrt{\tau}d(x,y)/2}+e^{-R\sqrt{\tau}/2})\tau^{1/2}\log\tau
\mbox{ when }d(x,y)\geq c\tau^{-1/2}\]\[ |d_A\theta (y)|\leq
2\tau^{1/2}w(y)+2\kappa\tau^{-1/2}\leq 4\kappa\tau^{-1/2}+2
M\tau^{1/2}\mskip-25mu\sum_{\{p:\theta(p)=0\}}\mskip-8mu\left(e^{-\sqrt{\tau}d(y,p)}{+}e^{-R\sqrt{\tau}}\right)\]

we deduce

\begin{cor}  There are constants $L,b>0$, depending only on the path
of almost complex structures $\{J_t\}_{t\in[0,1]}$, such that if
 $\tau$ is sufficiently large the path $(A(t),\theta(t))$ of $J_t$--vortices is
obtained by \eqref{trans}, we have, for all $t$,\[
|\dot{\theta}(t)(x)|\leq
L\left(\tau^{-1}+(\log\tau)\sum_{p:\theta(t)(p)=0}(e^{-b\sqrt{\tau}d(x,p)}+\tau
e^{-b\sqrt{\tau}})\right).\]\end{cor}

In particular, since $\theta(t)$, being a not--identically--zero
section of a degree $d$ holomorphic line bundle, vanishes at no more
than $d$ points, we have, where $w(t,x)=1-|\theta(t,x)|^2\geq 0$,
\[ \left|\frac{\partial w}{\partial t}\right|\leq
|\theta||\dot{\theta}|\leq 2dL\log\tau \] everywhere (we restrict
here to $\tau$ large enough that $\tau^{-1}+d\tau
e^{-b\sqrt{\tau}}\leq d$). So if $|h|\leq (4dL\log\tau)^{-1}$ and
$x\in\Sigma$ is such that $w(t+h,x)=1$, we must have had $w(t,x)\geq
1/2$. Referring back to the notation in \fullref{expldec},
assuming that $\tau$ is large enough that
$\kappa\tau^{-1}+Mde^{-R\sqrt{\tau}}\leq 1/4$, this implies that one
of the $d$ expressions $e^{-\sqrt{\tau}d(x,p)}$ for $p\in
\theta(t)^{-1}(0)$ must be at least $(4dM)^{-1}$, so that
$d(x,p)\leq B\tau^{-1/2}$ where the constant $B$ is independent of
$t$.  So since the points where $w(t+h,\cdot)$ is equal to $1$ are
those where $\theta(t+h)$ vanishes, we deduce that \emph{for all
$t\in[0,1]$, if $|h|\leq (4dL\log\tau)^{-1}$, then each zero of
$\theta(t+h)$ is a distance at most $B\tau^{-1/2}$ from a zero of
$\theta(t)$, and vice versa}, where the ``vice versa'' part comes
from just replacing $t$ by $t+h$ and $h$ by $-h$.  But then we can
subdivide $[0,1]$ into at most $(5dL\log\tau)$ intervals each of
length at most $(4dL\log\tau)^{-1}$ and apply this fact to the
endpoints of each interval to deduce that \begin{cor} Where
$N=5dLB$, each zero of $\theta(1)$ lies a distance at most\break
$N\tau^{-1/2}\log\tau$ from some zero of $\theta (0)$, and vice
versa.\end{cor}

Thus since the parallel transport map \[ F_{\{J_t\}}\co
S^d(\Sigma,J_0)\to S^d(\Sigma,\phi_{\omega}^{*}J_0) \quad
(\phi_{\omega}^{*}J_0=J_1) \] sends the zero set of $\theta(0)$ to
that of $\theta(1)$, we deduce that, as $\tau\to\infty$,
$F_{\{J_t\}}$ converges in $C^0$ norm to the identity, and so the
$\Omega_{d,\omega,\tau}$--symplectomorphisms
$\Phi_{d,\omega,\tau}=S^d\phi_{\omega}\circ F_{\{J_t\}}$ converge in
$C^0$-norm to $S^d\phi_{\omega}\co S^d\Sigma\to S^d\Sigma$ as
$\tau\to\infty$.

\bibliographystyle{gtart}
\bibliography{link}

\end{document}